\documentclass[10pt,a4paper]{article}
\usepackage{graphics}
\usepackage{amsmath,amssymb,amsthm,graphics,graphicx,epic,eepic,multicol,ascmac}
\usepackage[mathscr]{euscript}
\usepackage[all]{xy}
\usepackage{fancybox}
\usepackage{url}
\usepackage{ytableau}
\usepackage{here}

\ytableausetup{centertableaux,mathmode,boxsize=1em}
\def\young(#1){\ytableaushort{#1}}
\def\yng(#1){\ydiagram{#1}}

\numberwithin{equation}{section}
\theoremstyle{theorem}
\newtheorem{thm}{Theorem}[section]
\newtheorem{prop}[thm]{Proposition}

\theoremstyle{definition}
\newtheorem{defn}[thm]{Definition}
\newtheorem{ex}[thm]{Example}
\newtheorem{rem}[thm]{Remark}

\def\al{\alpha}

\def\wht(#1){\widehat{\ #1\ }}

\newcommand{\ch}{\mathrm{ch}}

\newcommand{\lbr}{\begin{bmatrix}}
\newcommand{\rbr}{\end{bmatrix}}

\newcommand{\cd}{commutative diagram }

\def\al{\alpha}

\def\beneme{\begin{enumerate}}
\def\beq{\begin{equation}}
\def\beqn{\begin{eqnarray}}
\def\beqnn{\begin{eqnarray*}}

\def\bfii0{{\bf i_0}}

\def\bbra#1,#2,#3{\left\{\begin{array}{c}\hspace{-5pt}
#1;#2\\ \hspace{-5pt}#3\end{array}\hspace{-5pt}\right\}}
\def\cd{\cdots}
\def\ci(#1,#2){c_{#1}^{(#2)}}
\def\Ci(#1,#2){C_{#1}^{(#2)}}
\def\mpp(#1,#2,#3){#1^{(#2)}_{#3}}
\def\bCi(#1,#2){\ovl C_{#1}^{(#2)}}
\def\ch(#1,#2){c_{#2,#1}^{-h_{#1}}}
\def\cc(#1,#2){c_{#2,#1}}

\def\di(#1,#2){D_{#1}^{(#2)}}
\def\dbi(#1,#2){\ovl D_{#1}^{(#2)}}

\def\eneme{\end{enumerate}}

\def\eeq{\end{equation}}
\def\eeqn{\end{eqnarray}}
\def\eeqnn{\end{eqnarray*}}

\def\gau#1,#2{\left[\begin{array}{c}\hspace{-5pt}#1\\
\hspace{-5pt}#2\end{array}\hspace{-5pt}\right]}

\def\ify{\infty}
\def\io{\iota}
\def\ji(#1,#2){j_{#1}^{(#2)}}

\def\lan{\langle}

\def\nd{\noindent}

\def\ovl{\overline}

\def\qed{\hfill\framebox[2mm]{}}
\def\QQ{\mathbb Q}

\def\ran{\rangle}

\def\TY(#1,#2,#3){#1^{(#2)}_{#3}}

\def\vp{\varphi}

\def\xxi(#1,#2,#3){\displaystyle {}^{#1}\Xi^{(#2)}_{#3}}
\def\xsi(#1,#2,#3){\displaystyle {}^{#1}\Sigma^{(#2)}_{#3}}
\def\xE(#1,#2,#3){\displaystyle {}^{#1}E_{#2}[#3]}
\def\xF(#1,#2){\displaystyle {}^{#1}F_{#2}}
\def\xx(#1,#2){\displaystyle {}^{#1}\Xi_{#2}}
\def\W1{W(\varpi_1)}

\def\what{\widehat}

\def\ZZ{\mathbb Z}

\def\m@th{\mathsurround=0pt}
\def\fsquare(#1,#2){
\hbox{\vrule$\hskip-0.4pt\vcenter to #1{\normalbaselines\m@th
\hrule\vfil\hbox to #1{\hfill$\scriptstyle #2$\hfill}\vfil\hrule}$\hskip-0.4pt
\vrule}}

\newcommand{\ba}{\begin{array}}
\newcommand{\ea}{\end{array}}

\newcommand{\eq}{\begin{eqnarray}}
\newcommand{\eneq}{\end{eqnarray}}

\title{\textbf{\large{Polyhedral realizations for crystal bases of integrable highest weight modules
and combinatorial objects of type
${\rm A}^{(1)}_{n-1}$, ${\rm C}^{(1)}_{n-1}$, ${\rm A}^{(2)}_{2n-2}$, ${\rm D}^{(2)}_{n}$
}}}
\author{\normalsize{YUKI KANAKUBO\thanks{Faculty of Pure and Applied Sciences, University of Tsukuba,
1-1-1 Tennodai, Tsukuba, Ibaraki 305-8577,
Japan: {y-kanakubo@math.tsukuba.ac.jp}.}}
}
\date{}

\begin{document}

\maketitle
\vspace{-10pt}

\begin{abstract}

In this paper, we consider polyhedral realizations for crystal bases $B(\lambda)$
of irreducible integrable highest weight modules of a
quantized enveloping algebra $U_q(\mathfrak{g})$,
where $\mathfrak{g}$ is a classical affine Lie algebra of type
${\rm A}^{(1)}_{n-1}$, ${\rm C}^{(1)}_{n-1}$, ${\rm A}^{(2)}_{2n-2}$ or ${\rm D}^{(2)}_{n}$.
We will give explicit forms of polyhedral realizations
in terms of extended Young diagrams or Young walls that appear in
the representation theory of
quantized enveloping algebras 
of classical affine type.
As an application, a combinatorial description of $\varepsilon_k^*$
functions on $B(\infty)$ will be given.

\end{abstract}

\section{Introduction}

The crystal bases introduced in
\cite{K1, L} are powerful
tools to study the representation theory of quantized enveloping algebras $U_q(\mathfrak{g})$
for symmetrizable Kac-Moody Lie algebras $\mathfrak{g}$.
Realizing the crystal bases as combinatorial objects like as Young tableaux, LS paths, Nakajima monomials \cite{K, KN, Lit2, Nj}, 
one can reveal skeleton structures of representations or negative part $U_q^-(\mathfrak{g})\subset U_q(\mathfrak{g})$.

We are interested in the crystal bases $B(\lambda)$ for irreducible integrable highest weight modules $V(\lambda)$ with a dominant integral weight $\lambda$
and the crystal bases $B(\infty)$ for the negative part $U_q^-(\mathfrak{g})$.
In \cite{N99, NZ}, {\it polyhedral realizations} of $B(\lambda)$ and $B(\infty)$ are introduced,
which realize $B(\lambda)$ and $B(\infty)$ as sets of integer points
in polytopes and polyhedral cones, respectively.
There are embeddings of crystals $\Psi_{\iota}:B(\infty)\hookrightarrow \mathbb{Z}^{\infty}_{\iota}$
and $\Psi_{\iota}^{(\lambda)}:B(\lambda)\hookrightarrow \mathbb{Z}^{\infty}_{\iota}\otimes R_{\lambda}$
associated with an infinite sequence $\iota$ from the index set $I=\{1,2,\cdots,n\}$ of $\mathfrak{g}$,
where $R_{\lambda}$ is a crystal, which has a single element (Example \ref{r-ex}) and
\[
\mathbb{Z}^{\infty}_{\iota}
=\{(\cd,a_r,\cd,a_2,a_1)| a_r\in\ZZ
\,\,{\rm and}\,\,a_r=0\,\,{\rm for}\,\,r\gg 0\}
\]
has a crystal structure associated with $\iota$.
Then procedures to compute
linear inequalities defining
${\rm Im}(\Psi_{\iota}) (\cong B(\infty))$ 
and ${\rm Im}(\Psi^{(\lambda)}_{\iota}) (\cong B(\lambda))$ are provided.
It is most fundamental problem in the theory of polyhedral realizations
to find explicit forms of 
${\rm Im}(\Psi_{\iota})$ and ${\rm Im}(\Psi^{(\lambda)}_{\iota})$.

The polyhedral realizations for $B(\infty)$ are deeply related to the {\it string cones} in \cite{BZ, Lit}. In fact,
in case of
$\mathfrak{g}$ is a finite dimensional simple Lie algebra and $\iota=(\cdots,i_{N+1},i_N,\cdots,i_2,i_1)$ is a sequence
such that
$(i_N,\cdots,i_2,i_1)$ is a reduced word of the longest element in the Weyl group $W$,
the image ${\rm Im}(\Psi_{\iota})$ coincides with the set of integer points in the string cone associated to the reduced word 
$(i_1,i_2\cdots,i_N)$ as a set, which is a polyhedral convex cone. 
When $\mathfrak{g}$ is of type ${\rm A}_n$, several combinatorial expressions of string cones
are given in
\cite{GKS16, GP}.

If $\mathfrak{g}$ is a finite dimensional simple Lie algebra and $\iota$ is a specific one then
explicit forms of ${\rm Im}(\Psi_{\iota})$ and ${\rm Im}(\Psi^{(\lambda)}_{\iota})$ are
given in \cite{H1, KS, Lit, N99,NZ}.
In case of 
$\mathfrak{g}$ is a classical affine Lie algebra and $\iota$ is specific one,
explicit forms of them are
given in \cite{H2, HN, N99, NZ}.
When $\mathfrak{g}$ is a finite dimensional classical Lie algebra
and $\iota$ is adapted (Definition \ref{adapt}), we explicitly describe inequalities defining polyhedral realizations of $B(\infty)$ and $B(\lambda)$
in terms of rectangular Young tableaux \cite{KaN, KaN2}. 

In the previous paper \cite{Ka}, we consider the case 
$\mathfrak{g}$ is a classical affine Lie algebra of type
${\rm A}^{(1)}_{n-1}$, ${\rm C}^{(1)}_{n-1}$, ${\rm A}^{(2)}_{2n-2}$ or ${\rm D}^{(2)}_{n}$
and $\iota$ is adapted.
Then we give an explicit form of linear inequalities defining polyhedral realizations of $B(\infty)$
in terms of extended Young diagrams and proper Young walls. 
Extended Young diagrams are used
to realize irreducible integrable highest weight modules $V(\Lambda)$ of $U_q(\mathfrak{g})$
as Fock space representations
for an affine Lie algebra $\mathfrak{g}$ of type ${\rm A}^{(1)}_n$, ${\rm C}^{(1)}_n$, ${\rm A}^{(2)}_{2n}$ or ${\rm D}^{(2)}_{n+1}$
and almost all fundamental weights $\Lambda$ \cite{Ha, JMMO, KMM}.
The proper Young walls are introduced in \cite{Kang}
to combinatorially realize the crystal base $B(\lambda)$ of level $1$ representations $V(\lambda)$ of $U_q(\mathfrak{g})$
for several affine Lie algebras $\mathfrak{g}$ (see also \cite{KK}).
The detail of result in \cite{Ka} is as follows:
Defining ${\rm EYD}_k$ as the set of
extended Young diagrams with $y_{\infty}=k\in I$ (see Definition \ref{EYD-sub}),
one can assign linear functions $L^{{\rm A}^{(1)}}_{s,k,\iota}(T)$
and $L^{{\rm C}^{(1)}}_{s,k,\iota}(T)$
for each $T\in {\rm EYD}_k$ and $s\in\mathbb{Z}$ (see (\ref{ovlL2})).
We also define ${\rm REYD}^{{\rm A}^{(2)}}_{k}$ and ${\rm REYD}^{{\rm D}^{(2)}}_{k}$
as sets of revised extended Young diagrams (Definition \ref{ad-rem-pt2}, \ref{ad-rem-pt3})
and ${\rm YW}^{{\rm A}^{(2)}}_{k}$ and ${\rm YW}^{{\rm D}^{(2)}}_{k}$
as the sets of proper Young walls (Definition \ref{def-YW}, \ref{def-YW2}).
One can assign linear functions $L^{{\rm A}^{(2)}}_{s,k,\iota}(T)$ and 
$L^{{\rm D}^{(2)}}_{s,k,\iota}(T)$ for each revised extended Young diagram or proper Young wall $T$ ((\ref{L1kdef}), (\ref{L11-def}))
and $s\in\mathbb{Z}$.
Setting
${\rm Comb}_{\iota}^{X}[\infty]$ as in (\ref{comb-inf-1}) and (\ref{comb-inf-2}),
we obtain 
\[
{\rm Im}(\Psi_{\iota})=
\{
\mathbf{a}\in\mathbb{Z}^{\infty}_{\iota}|
\varphi(\mathbf{a})\geq0\text{ for any }\varphi\in {\rm Comb}_{\iota}^{X^L}[\infty]
\},
\]
where, $X$ is the type of $\mathfrak{g}$ and $X^L={\rm A}^{(1)}, {\rm D}^{(1)}, {\rm A}^{(2)}$ or ${\rm C}^{(2)}$
when $X={\rm A}^{(1)}_{n-1}$, ${\rm C}^{(1)}_{n-1}$, ${\rm A}^{(2)}_{2n-2}$ or ${\rm D}^{(2)}_{n}$, respectively.
That is, the inequalities defining ${\rm Im}(\Psi_{\iota})$
can be expressed in terms of combinatorial objects that deeply related
to fundamental representations of $U_q( ^L\mathfrak{g})$, where
$^L\mathfrak{g}$ is the classical affine Lie algebra whose generalized Cartan matrix is
the transposed matrix of that of $\mathfrak{g}$.

In this paper, 
assuming that
$\mathfrak{g}$ is a classical affine Lie algebra of type
${\rm A}^{(1)}_{n-1}$, ${\rm C}^{(1)}_{n-1}$, ${\rm A}^{(2)}_{2n-2}$ or ${\rm D}^{(2)}_{n}$
and $\iota$ is adapted,
we give an explicit form of inequalities defining 
${\rm Im}(\Psi_{\iota}^{(\lambda)})$ in terms of 
combinatorial objects like as extended Young diagrams and Young walls.
Depending on $\iota$ and $k\in I$,
we will define
${\rm Comb}^X_{k,\iota}[\lambda]$
as a set
of
(i) a single function,
(ii) functions parametrized by boxes with integers $\mathbb{Z}_{\leq k}$,
(iii) functions parametrized by boxes with integers $\mathbb{Z}_{\geq k+1}$,
or 
(iv) functions parametrized by (revised) extended Young diagrams or Young walls
((\ref{comb-lam-1}), (\ref{comb-lam-2}), (\ref{comb-lam-3})). 
The functions in (ii), (iii) are related to crystal bases for extremal weight modules of $U_q( ^L\mathfrak{g})$ (see Remark \ref{rem-1}).
Then our main results Theorem \ref{thmA1} and \ref{thmA2} say
\[
{\rm Im}(\Psi_{\iota}^{(\lambda)})
=\left\{ \mathbf{a}\in\mathbb{Z}^{\infty}
\left|
\varphi(\mathbf{a})\geq0\\
\text{ for any }
\varphi\in
{\rm Comb}^{X^L}_{\iota}[\infty]\cup
\bigcup_{k\in I}
{\rm Comb}^{X^L}_{k,\iota}[\lambda]
\right.
\right\}.
\]
In this way, inequalities defining ${\rm Im}(\Psi_{\iota}^{(\lambda)})$
can be described by
combinatorial objects related
to representations of $U_q( ^L\mathfrak{g})$,
As an application, we will give a combinatorial
description of $\varepsilon_k^*$-functions on $B(\infty)$ (Theorem \ref{thm3}).
In case of other classical affine types ${\rm B}_{n-1}^{(1)}$, ${\rm D}_{n-1}^{(1)}$ and ${\rm A}_{2n-3}^{(2)}$,
we do not describe the inequalities in this paper
since
the combinatorial objects related to representations of $U_q( ^L\mathfrak{g})$ are complicated
comparing with the four types in the title.

The organization
of this paper is as follows:
In Sect.2, we will review the notion of crystals and polyhedral realizations.
A procedure to compute explicit forms of ${\rm Im}(\Psi_{\iota}^{(\lambda)})$ is
given in this section.
In Sect.3, we recall several combinatorial objects like as (revised) extended Young diagrams
and Young walls.
Sect.4 is devoted to present
our main results, which describe explicit forms of ${\rm Im}(\Psi_{\iota}^{(\lambda)})$
in terms of combinatorial objects seen in Sect.3.
In Sect.5 and Sect.6, we give proofs of main theorems.

\vspace{2mm}

\nd
{\bf Acknowledgements}
This work was supported by JSPS KAKENHI Grant Number JP20J00186.

\section{Polyhedral realizations and a procedure}

In this section, we review crystals and polyhedral realizations.
We also give a modified procedure computing explicit forms of polyhedral realizations for $B(\lambda)$.

\subsection{Notations}

Let $\mathfrak{g}$ denote a symmetrizable Kac-Moody algebra over $\mathbb{Q}$ with the index set $I=\{1,2,\cdots,n\}$ and
a generalized Cartan matrix $A=(a_{i,j})_{i,j\in I}$. Let  $\mathfrak{h}$ be a Cartan subalgebra, 
$\langle \cdot,\cdot \rangle : \mathfrak{h} \times \mathfrak{h}^*\rightarrow \mathbb{Q}$
be the canonical pairing,
$P\subset \mathfrak{h}^*$ a weight lattice,
$\{\alpha_i\}_{i\in I}$ a set of simple roots and $\{h_i\}_{i\in I}$ a set of simple coroots. 
It holds $\langle h_{i},\alpha_j \rangle=a_{i,j}$. 
We also define
$P^*:=\{h\in\mathfrak{h} | \langle h,P \rangle\subset\mathbb{Z}\}$ and
$P^+:=\{\lambda\in P | \langle h_i,\lambda \rangle \in\mathbb{Z}_{\geq0} \text{ for all }i\in I \}$.
For each $i\in I$, the $i$-th fundamental weight $\Lambda_i\in P^+$ is defined as
$\langle h_j,\Lambda_i \rangle=\delta_{i,j}$.
Let $U_q(\mathfrak{g})$ be the quantized universal enveloping algebra
with generators $e_i$, $f_i$ ($i\in I$) and $q^h$ ($h\in P^*$) and $U_q^-(\mathfrak{g})$ be
the subalgebra of $U_q(\mathfrak{g})$ generated by $f_i$ ($i\in I$).
It is known that the algebra $U_q^-(\mathfrak{g})$ has
a crystal base $(L(\infty),B(\infty))$
and the irreducible integrable highest weight module $V(\lambda)$ of $U_q(\mathfrak{g})$
has a crystal base $(L(\lambda),B(\lambda))$ for $\lambda\in P^+$ (see \cite{K0,K1}).

\begin{rem}\label{rem1}
In this paper, we mainly treat Lie algebras of type
${\rm A}^{(1)}_{n-1}$, ${\rm C}^{(1)}_{n-1}$, ${\rm A}^{(2)}_{2n-2}$ and ${\rm D}^{(2)}_{n}$.
For the sake of simplicity, we often write them by
${\rm A}^{(1)}$, ${\rm C}^{(1)}$, ${\rm A}^{(2)}$ and ${\rm D}^{(2)}$.
\end{rem}
\nd
We will use the following
numbering of vertices in affine Dynkin diagrams:
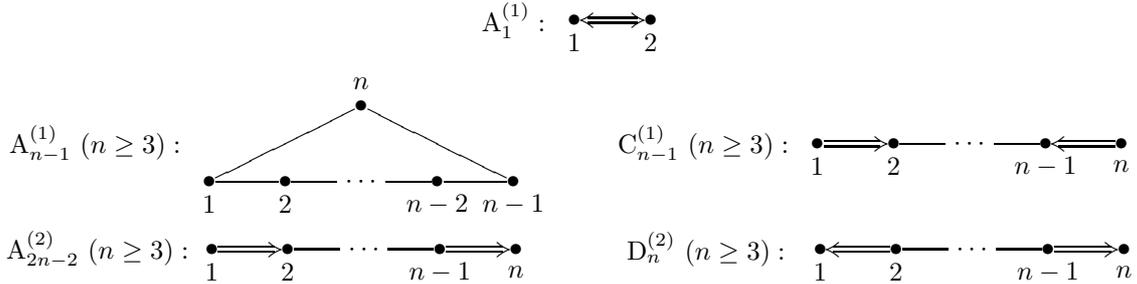
\[
\begin{xy}
(-8,0) *{{\rm A}_{1}^{(1)} : }="A1",
(0,0) *{\bullet}="1",
(0,-3) *{1}="1a",
(10,0)*{\bullet}="2",
(10,-3) *{2}="2a",
\ar@{<=>} "1";"2"^{}
\end{xy}
\]
\[
\begin{xy}
(-15,5) *{{\rm A}_{n-1}^{(1)} \ (n\geq 3) : }="A1",
(20,10) *{\bullet}="n",
(20,13) *{n}="na",
(0,0) *{\bullet}="1",
(0,-3) *{1}="1a",
(10,0)*{\bullet}="2",
(10,-3)*{2}="2a",
(20,0)*{\ \cdots\ }="3",
(30,0)*{\bullet}="4",
(30,-3)*{n-2}="4a",
(40,0)*{\bullet}="5",
(40,-3)*{n-1}="5a",
(65,5) *{{\rm C}_{n-1}^{(1)} \ (n\geq 3) : }="C1",
(80,5) *{\bullet}="11",
(80,2) *{1}="11a",
(90,5)*{\bullet}="22",
(90,2)*{2}="22a",
(100,5)*{\ \cdots\ }="33",
(110,5)*{\bullet}="44",
(110,2)*{n-1}="44a",
(120,5)*{\bullet}="55",
(120,2)*{n}="55a",
\ar@{-} "1";"2"^{}
\ar@{-} "2";"3"^{}
\ar@{-} "3";"4"^{}
\ar@{-} "4";"5"^{}
\ar@{-} "1";"n"^{}
\ar@{-} "n";"5"^{}
\ar@{=>} "11";"22"^{}
\ar@{-} "22";"33"^{}
\ar@{-} "33";"44"^{}
\ar@{<=} "44";"55"^{}
\end{xy}
\]
\[
\begin{xy}
(-15,5) *{{\rm A}_{2n-2}^{(2)}\ (n\geq 3) : }="A2",
(0,5) *{\bullet}="1",
(0,2) *{1}="1a",
(10,5)*{\bullet}="2",
(10,2)*{2}="2a",
(20,5)*{\ \cdots\ }="3",
(30,5)*{\bullet}="4",
(30,2)*{n-1}="4a",
(40,5)*{\bullet}="5",
(40,2)*{n}="5a",
(65,5) *{{\rm D}_{n}^{(2)} \ (n\geq 3) : }="D2",
(80,5) *{\bullet}="11",
(80,2) *{1}="11a",
(90,5)*{\bullet}="22",
(90,2)*{2}="22a",
(100,5)*{\ \cdots\ }="33",
(110,5)*{\bullet}="44",
(110,2)*{n-1}="44a",
(120,5)*{\bullet}="55",
(120,2)*{n}="55a",
\ar@{=>} "1";"2"^{}
\ar@{-} "2";"3"^{}
\ar@{-} "3";"4"^{}
\ar@{=>} "4";"5"^{}
\ar@{<=} "11";"22"^{}
\ar@{-} "22";"33"^{}
\ar@{-} "33";"44"^{}
\ar@{=>} "44";"55"^{}
\end{xy}
\]
Replacing our numbering $1,2,\cdots,n-1,n,n+1$ of ${\rm A}_{n}^{(1)}$, ${\rm C}_{n}^{(1)}$ and ${\rm D}_{n+1}^{(2)}$
with $0,1,2,\cdots,n$, one obtains the numbering in \cite{JMMO, Kang, KMM}.
Replacing our numbering $1,2,\cdots,n-1,n,n+1$ of ${\rm A}_{2n}^{(2)}$ with $n,n-1,\cdots,1,0$, we get the numbering in \cite{Kang, KMM}.

\subsection{Crystals}

We review the notion of {\it crystals} following \cite{K3}:

\begin{defn}
A set $\mathcal{B}$ together with the maps
${\rm wt}:\mathcal{B}\rightarrow P$,
$\varepsilon_k,\varphi_k:\mathcal{B}\rightarrow \mathbb{Z}\sqcup \{-\infty\}$
and $\tilde{e}_k$,$\tilde{f}_k:\mathcal{B}\rightarrow \mathcal{B}\sqcup\{0\}$
($k\in I$)
is said to be a {\it crystal} if 
the following relation holds: For $b,b'\in\mathcal{B}$ and $k\in I$,
\begin{itemize}
\item $\varphi_k(b)=\varepsilon_k(b)+\langle h_k,{\rm wt}(b)\rangle$,
\item ${\rm wt}(\tilde{e}_kb)={\rm wt}(b)+\alpha_k$ if $\tilde{e}_k(b)\in\mathcal{B}$,
\quad ${\rm wt}(\tilde{f}_kb)={\rm wt}(b)-\alpha_k$ if $\tilde{f}_k(b)\in\mathcal{B}$,
\item $\varepsilon_k(\tilde{e}_k(b))=\varepsilon_k(b)-1,\ \ 
\varphi_k(\tilde{e}_k(b))=\varphi_k(b)+1$\ if $\tilde{e}_k(b)\in\mathcal{B}$, 
\item $\varepsilon_k(\tilde{f}_k(b))=\varepsilon_k(b)+1,\ \ 
\varphi_k(\tilde{f}_k(b))=\varphi_k(b)-1$\ if $\tilde{f}_k(b)\in\mathcal{B}$, 
\item $\tilde{f}_k(b)=b'$ if and only if $b=\tilde{e}_k(b')$,
\item if $\varphi_k(b)=-\infty$ then $\tilde{e}_k(b)=\tilde{f}_k(b)=0$.
\end{itemize}
Here, $0$ and $-\infty$ are additional elements which do not belong to $\mathcal{B}$ and $\mathbb{Z}$, respectively.
The maps $\tilde{e}_k$,$\tilde{f}_k$ are called {\it Kashiwara operators}.
\end{defn}
The sets $B(\lambda)$ $(\lambda\in P^+)$, $B(\infty)$
have crystal structures. Let us review another crystal $R_{\lambda}$:

\begin{ex}\label{r-ex}
Let $R_{\lambda}:=\{r_{\lambda}\}$ ($\lambda\in P$) be the set consisting of a single element $r_{\lambda}$ and
defining
\[
{\rm wt}(r_{\lambda})=\lambda,\ \ \varepsilon_k(r_{\lambda})= -\lan h_k,\lambda\ran,\ \ \varphi_k(r_{\lambda})=0,\ \ 
\tilde{e}_k(r_{\lambda})=\tilde{f}_k(r_{\lambda})=0.
\]
Then $R_{\lambda}$ is a crystal.
\end{ex}

\begin{ex}\label{star-ex}
Let $\{\tilde{e}_i\}_{i\in I},\{\tilde{f}_i\}_{i\in I},\{\varepsilon_i\}_{i\in I},\{\varphi_i\}_{i\in I}, {\rm wt}$ be the crystal
structure on $B(\infty)$ and
$*:U_q(\frak{g})\rightarrow U_q(\frak{g})$ the antiautomorphism
such that $e_i^*=e_i$, $f_i^*=f_i$ and $(q^h)^*=q^{-h}$ in \cite{K1}.
It is known that the antiautomorphism $*$ induces a bijection $*:B(\infty)\rightarrow B(\infty)$ satisfying $*\circ *=id$.
One can define a crystal $B(\infty)^{*}$ as follows : The underlying set $B(\infty)^{*}$ equals $B(\infty)$,
the maps on $B(\infty)^{*}$ are defined as $\tilde{e}_i^*:=*\circ \tilde{e}_i\circ *$, 
$\tilde{f}_i^*:=*\circ \tilde{f}_i\circ *$,
$\varepsilon_i^*:=\varepsilon_i\circ *$,
$\varphi_i^*:=\varphi_i\circ *$ and ${\rm wt}^*:={\rm wt}$.
\end{ex}

\begin{defn}
Let $\mathcal{B}_1$, $\mathcal{B}_2$ be two crystals.
\begin{enumerate}
\item[$(1)$]
A map $\psi : \mathcal{B}_1\sqcup\{0\}\rightarrow \mathcal{B}_2\sqcup\{0\}$
is said to be a strict morphism 
and denoted by $\psi : \mathcal{B}_1\rightarrow \mathcal{B}_2$
if $\psi(0)=0$ and the following holds:
\begin{itemize}
\item For $k\in I$, if $b\in\mathcal{B}_1$ and $\psi(b)\in \mathcal{B}_2$ then
\[
{\rm wt}(\psi(b))={\rm wt}(b),\quad
\varepsilon_k(\psi(b))=\varepsilon_k(b),\quad
\varphi_k(\psi(b))=\varphi_k(b),
\]
\item $\tilde{e}_k(\psi(b))=\psi(\tilde{e}_k(b))$ and $\tilde{f}_k(\psi(b))=\psi(\tilde{f}_k(b))$ for $k\in I$ and $b\in\mathcal{B}_1$, where
we understand $\tilde{e}_k(0)=\tilde{f}_k(0)=0$.
\end{itemize}
\item[$(2)$] If a strict morphism 
$\psi : \mathcal{B}_1\sqcup\{0\}\rightarrow \mathcal{B}_2\sqcup\{0\}$ is injective then
$\psi$ is said to be a strict embedding and denoted by $\psi:\mathcal{B}_1 \hookrightarrow \mathcal{B}_2$.
\item[$(3)$] If a strict morphism 
$\psi : \mathcal{B}_1\sqcup\{0\}\rightarrow \mathcal{B}_2\sqcup\{0\}$ is bijective then
$\psi$ is said to be an isomorphism.
\end{enumerate}
\end{defn}

\begin{defn}
For two crystals $\mathcal{B}_1$, $\mathcal{B}_2$,
the tensor product $\mathcal{B}_1\otimes\mathcal{B}_2$
is defined to be the set
$\mathcal{B}_1\times\mathcal{B}_2$ whose crystal structure is defined as follows: For $k\in I$, $b_1\in\mathcal{B}_1$,
and $b_2\in\mathcal{B}_2$,
\begin{enumerate}
\item[$(1)$] ${\rm wt}(b_1\otimes b_2)={\rm wt}(b_1)+{\rm wt}(b_2)$,
\item[$(2)$] $\varepsilon_k(b_1\otimes b_2)={\rm max}\{\varepsilon_k(b_1), \varepsilon_k(b_2)-\langle h_k, {\rm wt}(b_1)\rangle\}$,
\item[$(3)$] $\varphi_k(b_1\otimes b_2)={\rm max}\{\varphi_k(b_2), \varphi_k(b_1)+\langle h_k, {\rm wt}(b_2)\rangle\}$,
\item[$(4)$]$\tilde{e}_k (b_1\otimes b_2)=
\begin{cases}
\tilde{e}_k b_1\otimes b_2 & {\rm if}\ \varphi_k(b_1)\geq \varepsilon_k(b_2),\\
b_1\otimes \tilde{e}_kb_2 & {\rm if}\ \varphi_k(b_1)< \varepsilon_k(b_2),
\end{cases}$
\item[$(5)$]$\tilde{f}_k (b_1\otimes b_2)=
\begin{cases}
\tilde{f}_k b_1\otimes b_2 & {\rm if}\ \varphi_k(b_1)> \varepsilon_k(b_2),\\
b_1\otimes \tilde{f}_kb_2 & {\rm if}\ \varphi_k(b_1)\leq \varepsilon_k(b_2).
\end{cases}$
\end{enumerate}
Here, we understand $0\otimes b_2=b_1\otimes0=0$.
\end{defn}

\subsection{Embeddings of crystal bases}\label{def-z}

One defines
\[
\ZZ^{\ify}
:=\{(\cd,a_r,\cd,a_2,a_1)| a_r\in\ZZ\text{ for }r\in\mathbb{Z}_{\geq1}
\text{ and it holds }a_r=0\,\,{\rm for}\ r\gg 0\}
\]
and takes an infinite sequence
$\io=(\cd,i_r,\cd,i_2,i_1)$ of indices from $I$ such that
\begin{equation}
{\hbox{
$i_r\ne i_{r+1}$ for $r\in\mathbb{Z}_{\geq1}$ and $\sharp\{r\in\mathbb{Z}_{\geq1}| i_r=k\}=\ify$ for all $k\in I$.}}
\label{seq-con}
\end{equation}
The set $\ZZ^{\ify}$ has a crystal structure associated with $\iota$ (see subsection 2.4 of \cite{NZ}).
Let $\ZZ^{\ify}_{\iota}$ denote its crystal.

\begin{thm}\cite{K3,NZ}
\label{emb}
There uniquely exists a strict embedding of crystals
\begin{equation}
\Psi_{\io}:B(\ify)\hookrightarrow \ZZ^{\ify}_{\io}
\label{psi}
\end{equation}
such that 
$\Psi_{\io} (u_{\ify}) = \textbf{0}$, where
$u_{\ify}\in B(\ify)$ is
the highest weight vector, $\textbf{0}:=(\cd,0,\cd,0,0)$.
\end{thm}

Recall the crystal $R_{\lambda}$ defined in Example \ref{r-ex}.
We can explicitly write the crystal structure on $\mathbb{Z}^{\infty}[\lambda]:=\mathbb{Z}^{\infty}_{\iota}\otimes R_{\lambda}$ as follows:
Since $R_{\lambda}$ has only one element, one can identify $\mathbb{Z}^{\infty}[\lambda]$ with
$\mathbb{Z}^{\infty}$ as a set.
We consider the infinite dimensional vector space
\[
\QQ^{\ify}:=\{\textbf{a}=
(\cd,a_r,\cd,a_2,a_1)| a_r \in \QQ\text{ for }r\in\mathbb{Z}_{\geq1}
\text{ and it holds }
a_r = 0\,\,{\rm for\ } r \gg 0\},
\]
and its dual space $(\QQ^{\ify})^*:={\rm Hom}(\QQ^{\ify},\QQ)$.
Let $x_r\in (\QQ^{\ify})^*$ denote the linear function defined as $x_r((\cd,a_r,\cd,a_2,a_1)):=a_r$
for $r\in\mathbb{Z}_{\geq1}$ and
\[
\sigma_r:=x_r+\sum_{j>r}\langle h_{i_r},\alpha_{i_j}\rangle x_j\in(\QQ^{\ify})^*\quad (r\in\mathbb{Z}_{\geq1}),\qquad
\sigma_0^{(k)}:=-\langle h_k,\lambda \rangle+ \sum_{j\geq1} \langle h_k,\alpha_{i_j}\rangle x_j\quad (k\in I),
\]
where 
$\sigma_0^{(k)}$ is a map from $\QQ^{\ify}$ to $\mathbb{Q}$ defined as
$\sigma_0^{(k)}(\mathbf{a})=-\langle h_k,\lambda \rangle+ \sum_{j\geq1} \langle h_k,\alpha_{i_j}\rangle a_j$
for $\mathbf{a}=(\cd,a_j,\cd,a_2,a_1)\in\QQ^{\ify}$.
Since $(\cd,a_j,\cd,a_2,a_1)\in\QQ^{\ify}$ satisfies $a_j=0$ for $j\gg0$,
these $\sigma_r$ and $\sigma_0^{(k)}$ are well-defined. 
For $\mathbf{a}=(\cd,a_j,\cd,a_2,a_1)\in\QQ^{\ify}$ and $k\in I$,
let $\sigma^{(k)}(\mathbf{a}):={\rm max} \{\sigma_r(\mathbf{a}) | r\in\mathbb{Z}_{\geq1},\ i_r=k\}$
and
\[
M^{(k)}(\mathbf{a}):=\{r\in\mathbb{Z}_{\geq1} | i_r=k,\ \sigma_r(\mathbf{a})=\sigma^{(k)}(\mathbf{a})\}.
\]
We see that $\sigma^{(k)}(\mathbf{a})\geq0$ and the set
$M^{(k)}(\mathbf{a})$ is finite if and only if 
$\sigma^{(k)}(\mathbf{a})>0$. Then ${\rm wt} : \mathbb{Z}^{\infty}[\lambda]\rightarrow P$, 
$\varepsilon_k, \varphi_k : \mathbb{Z}^{\infty}[\lambda]\rightarrow \mathbb{Z}$ $(k\in I)$ are defined by
\[
{\rm wt}(\mathbf{a}):=\lambda-\sum_{j=1}^{\infty} a_j\alpha_{i_j},\quad
\varepsilon_k(\mathbf{a}):={\rm max}\{\sigma^{(k)}(\mathbf{a}),\sigma_0^{(k)}(\mathbf{a})\},\quad
\varphi_k(\mathbf{a}):=\langle h_k,{\rm wt}(\mathbf{a})\rangle + \varepsilon_k(\mathbf{a})
\]
for $\mathbf{a}=(\cd,a_j,\cd,a_2,a_1)\in\mathbb{Z}^{\infty}[\lambda]$. Next, let us consider
Kashiwara operators $\tilde{f}_k, \tilde{e}_k: \mathbb{Z}^{\infty}[\lambda]\sqcup\{0\}\rightarrow
\mathbb{Z}^{\infty}[\lambda]\sqcup\{0\}$. They are defined by $\tilde{f}_k(0)=\tilde{e}_k(0)=0$ and
\begin{equation}\label{Kashi-def1}
(\tilde{f}_k(\mathbf{a}))_r : = 
a_r +\delta_{r,{\rm min} M^{(k)}(\mathbf{a})} \text{ if }\sigma^{(k)}(\mathbf{a})>\sigma_0^{(k)}(\mathbf{a}) ; \text{ otherwise }\tilde{f}_k(\mathbf{a})=0,
\end{equation}
\begin{equation}\label{Kashi-def2}
(\tilde{e}_k(\mathbf{a}))_r : = 
a_r -\delta_{r,{\rm max} M^{(k)}(\mathbf{a})} \text{ if }\sigma^{(k)}(\mathbf{a})>0\text{ and }\sigma^{(k)}(\mathbf{a})\geq\sigma_0^{(k)}(\mathbf{a})
 ; \text{ otherwise }\tilde{e}_k(\mathbf{a})=0
\end{equation}
for $\mathbf{a}=(\cd,a_j,\cd,a_2,a_1)\in\mathbb{Z}^{\infty}[\lambda]$.
Then above $\{\tilde{e}_k\}_{k\in I}$, $\{\tilde{f}_k\}_{k\in I}$,
$\{\varepsilon_k\}_{k\in I}$, $\{\varphi_k\}_{k\in I}$ and ${\rm wt}$ give
the crystal structure of the set $\mathbb{Z}^{\infty}[\lambda]$ \cite{N99}.
Let $\mathbb{Z}^{\infty}_{\iota}[\lambda]$ denote this crystal.

\vspace{2mm}

For $\lambda\in P^+$,
setting
\begin{equation}\label{Btilde}
\tilde{B}(\lambda)
:=\{
b\otimes r_{\lambda}\in B(\infty)\otimes R_{\lambda} |
\varepsilon_i^*(b)\leq \lan h_i,\lambda\ran \text{ for all }i\in I
\},
\end{equation}
we obtain $B(\lambda)\cong \tilde{B}(\lambda)$ \cite{N99},
where $\varepsilon_i^*$ is the function in Example \ref{star-ex}.
Hence, there exists
a strict embedding of crystals $\Omega : B(\lambda)\hookrightarrow B(\infty)\otimes R_{\lambda}$.
Using the embedding $\Psi_{\iota}$ in Theorem \ref{emb},
the following theorem holds:
\begin{thm}\cite{N99}\label{emb2}
The map
\[
\Psi^{(\lambda)}_{\io}:B(\lambda)\overset{\Omega}{\hookrightarrow} B(\infty)\otimes R_{\lambda}
\overset{\Psi_{\iota}\otimes \text{id}}{\hookrightarrow} \mathbb{Z}^{\infty}_{\iota}\otimes R_{\lambda}
=\ZZ^{\ify}_{\iota}[\lambda]
\]
is the unique strict embedding of crystals
such that $\Psi^{(\lambda)}_{\io}(u_{\lambda})=(\cdots,0,0,0)$. Here,
$u_{\lambda}$ is the highest weight vector in $B(\lambda)$.
\end{thm}

\subsection{Modified Nakashima-Zelevinsky's procedure}\label{poly-uqm}

Let us review a procedure to compute an explicit form of ${\rm Im}(\Psi_{\iota})\subset\mathbb{Z}^{\infty}_{\iota}$ gave in 
the previous paper \cite{Ka}.
As in the previous subsection,
we fix an infinite sequence 
\[
\io=(\cd,i_r,\cd,i_2,i_1)
\]
of indices
satisfying (\ref{seq-con})
and for $r\in \mathbb{Z}_{\geq1}$, define $x_r\in(\QQ^{\ify})^*$ as
$x_r(\cd,a_r,\cd,a_2,a_1)=a_r$. For $r\in \mathbb{Z}_{<1}$, we set $x_r:=0$.
For $r\in \mathbb{Z}_{\geq1}$, let
\[
r^{(+)}:={\rm min}\{\ell\in\mathbb{Z}_{\geq1}\ |\ \ell>r\,\,{\rm and }\,\,i_r=i_{\ell}\},\ \ 
r^{(-)}:={\rm max} \{\ell\in\mathbb{Z}_{\geq1}\ |\ \ell<r\,\,{\rm and }\,\,i_r=i_{\ell}\}\cup\{0\},
\]
and for $r\in\mathbb{Z}_{\geq1}$, let
\begin{equation}
\beta_r:=
\sigma_r-\sigma_{r^{(+)}}
=x_r+\sum_{r<j<r^{(+)}}\lan h_{i_r},\al_{i_j}\ran x_j+x_{r^{(+)}}\in (\QQ^{\ify})^*
\label{betak}
\end{equation}
and $\beta_0:=0$. For $r\in\mathbb{Z}_{\geq1}$,
one defines the map
$S_r'=S_{r,\io}':(\QQ^{\ify})^*\rightarrow (\QQ^{\ify})^*$ as follows:
For $\vp=\sum_{\ell\in\mathbb{Z}_{\geq1}}c_{\ell}x_{\ell}\in(\QQ^{\ify})^*$,
\begin{equation}
S_r'(\vp):=
\begin{cases}
\vp-\beta_r & {\mbox{ if }}\ c_r>0,\\
 \vp+\beta_{r^{(-)}} & {\mbox{ if }}\ c_r< 0,\\
 \vp &  {\mbox{ if }}\ c_r= 0.
\end{cases}
\label{Sk}
\end{equation}
We often write $S_r'(\vp)$ as $S_r'\vp$. Let us define
\begin{eqnarray}
\Xi_{\io}' &:=  &\{S_{j_{\ell}}'\cd S_{j_2}'S_{j_1}'x_{j_0}\,|\,
\ell\in\mathbb{Z}_{\geq0},j_0,j_1,\cd,j_{\ell}\in\mathbb{Z}_{\geq1}\}, \label{xiprime}\\
\Sigma_{\io}' & := &
\{\textbf{a}\in \ZZ^{\ify}\subset \QQ^{\ify}\,|\,\vp(\textbf{a})\geq0\,\,{\rm for}\,\,
{\rm any}\,\,\vp\in \Xi_{\io}'\}.\nonumber
\end{eqnarray}
We say the sequence $\io$ satisfies the $\Xi'$-{\it positivity condition}
when it holds
\begin{equation}
{\hbox{for any 
$\vp=\sum_{\ell\in\mathbb{Z}_{\geq1}} c_{\ell}x_{\ell}\in \Xi_{\io}'$, if $\ell^{(-)}=0$ then $c_{\ell}\geq0$}}.
\label{posi}
\end{equation}
\begin{thm}\label{inf-thm}\cite{Ka}
Let $\io$ be a sequence of indices satisfying $(\ref{seq-con})$ 
and $(\ref{posi})$. Then we have
${\rm Im}(\Psi_{\io})=\Sigma_{\io}'$.
\end{thm}

\subsection{Modified Nakashima's procedure}\label{poly-uql} 

In this subsection, we give a procedure to compute an explicit form of the image
${\rm Im}(\Psi^{(\lambda)}_{\io})$ ($\cong B(\lambda)$) of 
$\Psi^{(\lambda)}_{\io}$ defined in Theorem \ref{emb2}.
Recall that we defined $\beta_r$ in $(\ref{betak})$.
Let $\beta^{(\pm)}_r$ be the functions given by $\beta^{(+)}_r=\beta_{r}$ and
\begin{equation*}
\beta^{(-)}_r = 
\begin{cases}
x_{r^{(-)}}+\sum_{r^{(-)}<j<r} \lan h_{i_r},\alpha_{i_j}\ran x_j +x_r=\sigma_{r^{(-)}}-\sigma_r & {\rm if}\ r^{(-)}>0,\\
- \lan h_{i_r},\lambda\ran + \sum_{1\leq j<r} \lan h_{i_r},\alpha_{i_j}\ran x_j + x_r
=\sigma_0^{(i_r)}-\sigma_r
 & {\rm if}\ r^{(-)}=0.
\end{cases}
\end{equation*}
Note that $\beta^{(-)}_r= \beta_{r^{(-)}}$ if $r^{(-)}>0$.

Using this notation, for each $r\in\mathbb{Z}_{\geq1}$
and $\varphi=c+\sum_{r\geq1} \varphi_rx_r$ with $c$, $\varphi_r\in\mathbb{Q}$,
we define a new function $\what{S}_r'(\varphi)$ as follows:
\begin{equation}\label{Shat}
\what{S}_r'(\varphi):=
\begin{cases}
\varphi - \beta^{(+)}_r & {\rm if}\ \varphi_r>0, \\
\varphi + \beta^{(-)}_r & {\rm if}\ \varphi_r<0, \\
\varphi & {\rm if}\ \varphi_r=0.
\end{cases}
\end{equation}
Note that except for the case $\varphi_r<0$ and $r^{(-)}=0$,
it holds 
\begin{equation}\label{twoS}
\what{S}_r'(\varphi)=S_r'(\varphi-c)+c,
\end{equation}
where $S_r'$ is defined in (\ref{Sk}).
For $k\in I$, one defines
\begin{equation}\label{lmi-def}
\iota^{(k)}:={\rm min}\{r\in\mathbb{Z}_{\geq1} | i_r=k\},\qquad
\lambda^{(k)} :=\lan h_k, \lambda \ran - \sum_{1\leq j<\iota^{(k)}} \lan h_k, \alpha_{i_j} \ran x_j - x_{\iota^{(k)}}.
\end{equation}
For $\iota$ and $\lambda\in P^+$, let $\Xi'_{\iota}[\lambda]$ be the set of all
functions generated by applying $\what{S}'_r$ on the functions $x_j$ ($j \geq 1$) and $\lambda^{(k)}$
($k\in I$), that is,
\begin{equation}\label{xilamdef}
\begin{array}{l}
\Xi'_{\iota}[\lambda]: =\{\what{S}'_{j_{\ell}}\cdots \what{S}'_{j_1}x_{j_0} | \ell\in\mathbb{Z}_{\geq0},\ j_0,\cdots,j_{\ell}\in \mathbb{Z}_{\geq1} \} \\
\qquad\qquad \cup \{\what{S}'_{j_{\ell}}\cdots \what{S}'_{j_1} \lambda^{(k)} | \ell\in\mathbb{Z}_{\geq0},\ k\in I,\ j_1,\cdots,j_{\ell}\in \mathbb{Z}_{\geq1} \}.
\end{array}
\end{equation}
By (\ref{twoS}), if $\iota$ satisfies the $\Xi'$-positivity condition (\ref{posi}) then
we see that
\begin{equation}\label{twoS2}
\{\what{S}'_{j_{\ell}}\cdots \what{S}'_{j_1}x_{j_0} | \ell\in\mathbb{Z}_{\geq0}, j_0,\cdots,j_{\ell}\in \mathbb{Z}_{\geq1} \}
=\{ S'_{j_{\ell}}\cdots S'_{j_1}x_{j_0} | \ell\in\mathbb{Z}_{\geq0}, j_0,\cdots,j_{\ell}\in \mathbb{Z}_{\geq1} \}.
\end{equation}
Now we set
\begin{equation}\label{siglamdef}
\Sigma_{\iota}'[\lambda]:=
\{
\mathbf{a}\in\mathbb{Z}^{\infty}_{\iota}[\lambda] | \varphi(\mathbf{a}) \geq0\ \ {\rm for\ any\ }\varphi\in \Xi'_{\iota}[\lambda] 
\}.
\end{equation}
\begin{defn}\label{ample1}
The pair $(\iota,\lambda)$ is said to be $\Xi'$-{\it ample}
if $\textbf{0}:=(\cdots,0,0,0)\in \Sigma'_{\iota}[\lambda]$.
\end{defn}

\begin{thm}\label{Nthm1}
For $\lambda\in P^+$,
we suppose that $(\iota,\lambda)$ is $\Xi'$-ample. 
Let $\Psi^{(\lambda)}_{\iota}:B(\lambda)\hookrightarrow \mathbb{Z}^{\infty}_{\iota}[\lambda]$
be the embedding of crystals in Theorem \ref{emb2}. Then
the image ${\rm Im}(\Psi^{(\lambda)}_{\iota}) (\cong B(\lambda))$ coincides with
$\Sigma'_{\iota}[\lambda]$.
\end{thm}

\nd
{\it Proof.} The proof is a slight modification of that of Theorem 4.1 in \cite{N99}.
Since ${\rm Im}(\Psi^{(\lambda)}_{\iota})$ is isomorphic to $B(\lambda)$ and $\Psi^{(\lambda)}_{\iota}$
is a strict embedding of crystals, it holds
\[
{\rm Im}(\Psi^{(\lambda)}_{\iota})=
\{
\tilde{f}_{j_{\ell}}\cdots\tilde{f}_{j_1}\Psi_{\iota}^{(\lambda)}(u_{\lambda}) | j_1,\cdots,j_{\ell}\in I,\ \ell\in\mathbb{Z}_{\geq0} \}\setminus
\{0\}
\]
and $\Psi_{\iota}^{(\lambda)}(u_{\lambda})=\mathbf{0}=(\cdots,0,0,0)$.
By the definition of $\tilde{f}_j$ $(j\in I)$ on $\mathbb{Z}^{\infty}_{\iota}[\lambda]$ in subsection \ref{def-z},
one can easily check
\[
{\rm Im}(\Psi_{\iota}^{(\lambda)})\subset 
\mathbb{Z}^{\infty}_{\geq0,\iota}[\lambda]:=\{ (\cdots,a_3,a_2,a_1)\in\mathbb{Z}^{\infty}_{\iota}[\lambda] | a_r\geq0\ (r\in\mathbb{Z}_{\geq1}) \}.
\]
First, let us show
\[
{\rm Im}(\Psi_{\iota}^{(\lambda)})\subset \Sigma_{\iota}'[\lambda].
\]
For any $k\in I$, let us show
\begin{equation}\label{Nthm-pr-1}
\tilde{f}_k \Sigma'_{\iota}[\lambda] \subset \Sigma'_{\iota}[\lambda]\cup\{0\}.
\end{equation}
For $\textbf{a}=(\cdots,a_3,a_2,a_1)\in \Sigma'_{\iota}[\lambda]$, we assume
\[
\tilde{f}_k\textbf{a}=(\cdots,a_{r+1},a_r+1,a_{r-1},\cdots,a_2,a_1)
\]
with some $r\in\mathbb{Z}_{\geq1}$. Note that $i_r=k$.
For any $\varphi=c+\sum_{j\in\mathbb{Z}_{\geq1}}\varphi_jx_j\in \Xi'_{\iota}[\lambda]$ ($c,\varphi_j\in\mathbb{Z}$),
we need to show
\begin{equation}\label{Nthm-pr-2}
\varphi(\tilde{f}_k \mathbf{a})\geq0.
\end{equation}
Since 
$\varphi(\tilde{f}_k \mathbf{a})=\varphi(\mathbf{a})+\varphi_r \geq \varphi_r$,
we may assume that $\varphi_r<0$. It follows from the definition of $\tilde{f}_k$ that
\[
\sigma_r(\mathbf{a})>\sigma_{r^{(-)}}(\mathbf{a})\ \text{if }r^{(-)}>0,\qquad
\sigma_r(\mathbf{a})>\sigma_{0}^{(k)}(\mathbf{a})\ \text{if }r^{(-)}=0
\]
so that $\beta^{(-)}_{r}(\mathbf{a})\leq -1$.
Hence,
\[
\varphi(\tilde{f}_k(\mathbf{a}))=
\varphi(\mathbf{a})+\varphi_r
\geq 
\varphi(\mathbf{a})-\varphi_r\beta^{(-)}_{r}(\mathbf{a})
=(\what{S}_r')^{-\varphi_r} (\varphi)(\mathbf{a})\geq0.
\]
Thus, we get (\ref{Nthm-pr-2}) and (\ref{Nthm-pr-1}).
By the assumption $\mathbf{0}\in \Sigma_{\iota}'[\lambda]$ and
(\ref{Nthm-pr-1}) imply ${\rm Im}(\Psi_{\iota}^{(\lambda)})\subset \Sigma_{\iota}'[\lambda]$.
Next, to prove
\[
{\rm Im}(\Psi_{\iota}^{(\lambda)})\supset \Sigma_{\iota}'[\lambda],
\]
let us show for any $k\in I$,
\begin{equation}\label{Nthm-pr-3}
\tilde{e}_k \Sigma'_{\iota}[\lambda] \subset \Sigma'_{\iota}[\lambda]\cup\{0\}.
\end{equation}
For $\textbf{a}=(\cdots,a_3,a_2,a_1)\in \Sigma'_{\iota}[\lambda]$, we assume
\[
\tilde{e}_k\textbf{a}=(\cdots,a_{r+1},a_r-1,a_{r-1},\cdots,a_2,a_1)
\]
with some $r\in\mathbb{Z}_{\geq1}$. Note that $i_r=k$.
For any $\varphi=c+\sum_{j\in\mathbb{Z}_{\geq1}}\varphi_jx_j\in \Xi'_{\iota}[\lambda]$ ($c,\varphi_j\in\mathbb{Z}$),
we need to show
\begin{equation}\label{Nthm-pr-4}
\varphi(\tilde{e}_k \mathbf{a})\geq0.
\end{equation}
Since $\varphi(\tilde{e}_k\mathbf{a})=\varphi(\mathbf{a})-\varphi_r\geq -\varphi_r$,
we may assume $\varphi_r>0$. By the definition of $\tilde{e}_k$, it follows $\sigma_r(\mathbf{a})-\sigma_{r^{(+)}}(\mathbf{a})>0$
so that $\beta^{(+)}_r(\mathbf{a})\geq1$.
Thus,
\[
\varphi(\tilde{e}_k(\mathbf{a}))=
\varphi(\mathbf{a})-\varphi_r
\geq 
\varphi(\mathbf{a})-\varphi_r\beta^{(+)}_{r}(\mathbf{a})
=(\what{S}_r')^{\varphi_r} (\varphi)(\mathbf{a})\geq0.
\]
Hence, one obtains (\ref{Nthm-pr-4}) and (\ref{Nthm-pr-3}).
Considering
$\Sigma'_{\iota}[\lambda]\subset 
\mathbb{Z}^{\infty}_{\geq0,\iota}[\lambda]$ and
(\ref{Nthm-pr-3}), for any $\textbf{a}=(\cdots,a_3,a_2,a_1)\in \Sigma'_{\iota}[\lambda]$,
taking $\ell\in\mathbb{Z}_{>0}$ sufficiently larger than $0$, it holds
\[
\tilde{e}_{j_1}\cdots \tilde{e}_{j_{\ell}}(\mathbf{a})=0
\]
for any $j_1,\cdots,j_{\ell}\in I$ (indeed, one can take $\ell=\sum_{j\geq1}a_j+1$).
Hence,
all we need to show is that if $\mathbf{a}\in \Sigma_{\iota}'[\lambda]$ satisfies
$\tilde{e}_k(\mathbf{a})=0$ for all $k\in I$ then $\mathbf{a}=\mathbf{0}$. In fact,
it means for any $\mathbf{y}\in\Sigma_{\iota}'[\lambda]$, there exist $j_1,\cdots,j_{\ell}\in I$
such that $\tilde{e}_{j_1}\cdots \tilde{e}_{j_{\ell}}(\mathbf{y})=\mathbf{0}$, which implies
\[
\mathbf{y}=\tilde{f}_{j_{\ell}}\cdots \tilde{f}_{j_{1}}(\mathbf{0})\in {\rm Im}(\Psi_{\iota}^{(\lambda)}).
\]
We assume $\mathbf{a}\in \Sigma_{\iota}'[\lambda]$, $\mathbf{a}\neq\mathbf{0}$ and $\tilde{e}_k(\mathbf{a})=0$ for all $k\in I$
and deduce a contradiction from this assumption.
Considering the definition of $\tilde{e}_k$, it holds either
\begin{equation}\label{Nthm-pr-5}
\sigma^{(k)}(\mathbf{a})\leq 0\quad \text{or}\quad \sigma^{(k)}(\mathbf{a})< \sigma^{(k)}_0(\mathbf{a}).
\end{equation}
The assumption
$\mathbf{a}=(\cdots,a_3,a_2,a_1)\neq\mathbf{0}$ means
there exists $j\in\mathbb{Z}_{\geq1}$ such that $a_j>0$ and $a_r=0$ when $r>j$. Thus,
$\sigma_j(\mathbf{a})=a_j>0$ and $\sigma^{(i_j)}(\mathbf{a})\geq a_j>0$.
The relation $\sigma^{(i_j)}(\mathbf{a})<\sigma^{(i_j)}_0(\mathbf{a})$ is an easy consequence of
(\ref{Nthm-pr-5}). Putting $k:=i_j$, 
\[
0<\sigma^{(k)}_0(\mathbf{a})-\sigma^{(k)}(\mathbf{a})
\leq \sigma^{(k)}_0(\mathbf{a})-\sigma_{\iota^{(k)}}(\mathbf{a})=\beta^{(-)}_{\iota^{(k)}}(\mathbf{a}).
\]
Hence,
\[
\lambda^{(k)}(\mathbf{a})=-\beta^{(-)}_{\iota^{(k)}}(\mathbf{a})<0,
\]
which contradicts $\mathbf{a}\in \Sigma_{\iota}'[\lambda]$. Therefore, one obtains $\mathbf{a}=\mathbf{0}$.\qed

\subsection{A description of $\varepsilon_k^*$}

In general, it seems to be difficult to find a value of $\varepsilon_k^*$ in Example \ref{star-ex}
by a direct calculation.
As an application of polyhedral realizations, we give a combinatorial description of 
$\varepsilon_k^*$, which is an analog of subsection 4.3 of \cite{N99}.

For $k\in I$, let
\[
\xi^{(k)}:=- \sum_{1\leq j<\iota^{(k)}} \lan h_k, \alpha_{i_j} \ran x_j - x_{\iota^{(k)}}\in (\mathbb{Q}^{\infty})^*.
\]
Note that for any $\lambda\in P^+$, one can write
$\xi^{(k)}=\lambda^{(k)}-\lan h_k, \lambda \ran$ by using $\lambda^{(k)}$ in
(\ref{lmi-def}). Using the maps in (\ref{Sk}), we set
\[
\Xi'^{(k)}_{\iota}:=\{S'_{j_{\ell}}\cdots S'_{j_1}\xi^{(k)} | \ell\in\mathbb{Z}_{\geq0},\ j_1,\cdots,j_{\ell}\in\mathbb{Z}_{\geq1}\}.
\]
We say the sequence $\iota$ satisfies the $\Xi'$-{\it strict positivity condition} when it holds the following:
\begin{equation}\label{strict-cond}
\text{for any }
\vp=\sum_{\ell\geq1} \varphi_{\ell}x_{\ell}\in \left(\bigcup_{k\in I}\Xi'^{(k)}_{\iota}\setminus\{\xi^{(k)}\}\right)\cup\Xi'_{\iota},
\text{ if }\ell^{(-)}=0\text{ then }\varphi_{\ell}\geq0,
\end{equation}
where $\Xi'_{\iota}$ is defined in (\ref{xiprime}).
\begin{thm}\label{thm2}
If $\iota$ satisfies the $\Xi'$-strict positivity condition then for $k\in I$ and $x\in {\rm Im}(\Psi_{\iota})$, we have
\[
\varepsilon^*_k(x)={\rm max}\{-\varphi(x)| \varphi\in \Xi'^{(k)}_{\iota}\}.
\]
\end{thm}
\nd
{\it Proof.} The proof is a slight modification of Theorem 4.2 of \cite{N99}. We take an arbitrary $\lambda\in P^+$.
Since $\iota$ satisfies the $\Xi'$-strict positivity condition, it also satisfies
$\Xi'$-positivity condition in (\ref{posi}) and it holds
\begin{equation}\label{thm2-pr0}
{\rm Im}(\Psi_{\iota})=
\{
\mathbf{a}\in\mathbb{Z}_{\iota}^{\infty} | \varphi(\mathbf{a})\geq0\ \text{for any }\varphi\in
\{S'_{j_{\ell}}\cdots S'_{j_1}x_{j_0} | \ell\in\mathbb{Z}_{\geq0}, j_0,\cdots,j_{\ell}\in \mathbb{Z}_{\geq1} \}
\}
\end{equation}
by Theorem \ref{inf-thm}. As checked in (\ref{twoS2}), we have
\begin{equation}\label{thm2-pr1}
\{\what{S}'_{j_{\ell}}\cdots \what{S}'_{j_1}x_{j_0} | \ell\in\mathbb{Z}_{\geq0}, j_0,\cdots,j_{\ell}\in \mathbb{Z}_{\geq1} \}
=\{ S'_{j_{\ell}}\cdots S'_{j_1}x_{j_0} | \ell\in\mathbb{Z}_{\geq0}, j_0,\cdots,j_{\ell}\in \mathbb{Z}_{\geq1} \}.
\end{equation}
Next, for $\ell\in\mathbb{Z}_{\geq0}$ and $j_1,\cdots,j_{\ell}\in \mathbb{Z}_{\geq1}$,
let us show if $\what{S}'_{j_{\ell}}\cdots \what{S}'_{j_1} \lambda^{(k)} \neq0$ then
\begin{equation}\label{thm2-pr2}
\what{S}'_{j_{\ell}}\cdots \what{S}'_{j_1} \lambda^{(k)} = \lan h_k,\lambda\ran + S'_{j_{\ell}}\cdots S'_{j_1} \xi^{(k)}
\end{equation}
by induction on $\ell$. When $\ell=0$ the claim is clear so that we assume $\ell>0$
and (\ref{thm2-pr2}) holds. We put $\varphi:=\lan h_k,\lambda\ran+\sum_{r\in\mathbb{Z}_{\geq1}}
\varphi_rx_r=\what{S}'_{j_{\ell}}\cdots \what{S}'_{j_1} \lambda^{(k)} = \lan h_k,\lambda\ran + S'_{j_{\ell}}\cdots S'_{j_1} \xi^{(k)}$
with some $\varphi_r\in\mathbb{Z}$.
For any $r\in\mathbb{Z}_{\geq1}$, except for the case $\varphi_r<0$ and $r^{(-)}=0$, it holds
\[
\what{S}'_r\what{S}'_{j_{\ell}}\cdots \what{S}'_{j_1} \lambda^{(k)} = \lan h_k,\lambda\ran + S'_rS'_{j_{\ell}}\cdots S'_{j_1} \xi^{(k)}
\]
by (\ref{twoS}). Thus, we may assume $\varphi_r<0$ and $r^{(-)}=0$. Then
$\Xi'$-strict positivity condition implies $S'_{j_{\ell}}\cdots S'_{j_1} \xi^{(k)}=\xi^{(k)}$ so that
$\what{S}'_{j_{\ell}}\cdots \what{S}'_{j_1} \lambda^{(k)} =\lambda^{(k)}$.
It follows from $\varphi_r<0$ and the definition of $\xi^{(k)}$ that $r=\iota^{(k)}$
and $\what{S}'_r\what{S}'_{j_{\ell}}\cdots \what{S}'_{j_1} \lambda^{(k)}=0$, this is not the case of (\ref{thm2-pr2}).
In this way, we have proved (\ref{thm2-pr2}). 

Since $S'_{j_{\ell}}\cdots S'_{j_1} \xi^{(k)}$ and
$S'_{j_{\ell}}\cdots S'_{j_1}x_{j_0}$ are linear maps and $\lan h_k,\lambda\ran\geq0$,  
we see that $(\iota,\lambda)$ is $\Xi'$-ample.
Combining Theorem \ref{Nthm1} with (\ref{thm2-pr0}), (\ref{thm2-pr1}) and (\ref{thm2-pr2}),
it holds
\begin{eqnarray}
{\rm Im}(\Psi_{\iota}^{(\lambda)})
&=&
\{
\mathbf{a}\otimes r_{\lambda}\in
{\rm Im}(\Psi_{\iota})\otimes R_{\lambda}
 | \varphi(\mathbf{a})\geq0\ \text{for any }\varphi\in
\{\what{S}'_{j_{\ell}}\cdots \what{S}'_{j_1}\lambda^{(k)} | k\in I,\ \ell\in\mathbb{Z}_{\geq0}, j_1,\cdots,j_{\ell}\in \mathbb{Z}_{\geq1} \}
\}\nonumber\\
&=&
\{
\mathbf{a}\otimes r_{\lambda}\in
{\rm Im}(\Psi_{\iota})\otimes R_{\lambda}
 | -\varphi(\mathbf{a})\leq\lan h_k,\lambda\ran\ \text{for any }k\in I\text{ and }\varphi\in
\Xi_{\iota}'^{(k)}\}
\}\nonumber\\
&=&
\{
\mathbf{a}\otimes r_{\lambda}\in
{\rm Im}(\Psi_{\iota})\otimes R_{\lambda}
 | {\rm max}\{-\varphi(\mathbf{a}) |\varphi\in
\Xi_{\iota}'^{(k)}\} \leq\lan h_k,\lambda\ran\ \text{for any }k\in I\}
\}.
\label{Btilde2}
\end{eqnarray}
Note that
${\rm max}\{-\varphi(\mathbf{a}) |\varphi\in
\Xi_{\iota}'^{(k)}\}\geq0$. In fact,
taking $r\in\mathbb{Z}_{\geq1}$
such that $a_s=0$ for $s\geq r$
and
applying maps $S'_{j_{\ell}},\cdots, S'_{j_1}$ to $\xi^{(k)}$ properly,
$S'_{j_{\ell}},\cdots, S'_{j_1}\xi^{(k)}\in
\Xi_{\iota}'^{(k)}$ becomes a sum of
$-x_t$ ($t\geq1$) and $x_{s}$ ($s\geq r$) by the definition of $S'$
so that $-S'_{j_{\ell}},\cdots, S'_{j_1}\xi^{(k)}(\mathbf{a})\geq0$.
Comparing (\ref{Btilde2}) with (\ref{Btilde}), one obtains
$\varepsilon^*_k(x)={\rm max}\{-\varphi(x)| \varphi\in \Xi'^{(k)}_{\iota}\}$ on $x\in{\rm Im}(\Psi_{\iota})$.
In fact, for a fixed $x\in{\rm Im}(\Psi_{\iota})$,
 because $\lambda\in P^+$ is arbitrary, one can take $\lambda$ as 
$\varepsilon_k^*(x)=\lan h_k,\lambda\ran$ for all $k\in I$.
Using (\ref{Btilde}), it follows $x\otimes r_{\lambda}\in {\rm Im}(\Psi_{\iota}^{(\lambda)})$
so that $\varepsilon_k^*(x)=\lan h_k,\lambda\ran\geq {\rm max}\{-\varphi(x) |\varphi\in
\Xi_{\iota}'^{(k)}\}$ holds by (\ref{Btilde2}). Similarly, taking $\lambda$ as 
$ {\rm max}\{-\varphi(x) |\varphi\in
\Xi_{\iota}'^{(k)}\}=\lan h_k,\lambda\ran$ for all $k\in I$, the expression (\ref{Btilde2})
yields $x\otimes r_{\lambda}\in {\rm Im}(\Psi_{\iota}^{(\lambda)})$
and
$ {\rm max}\{-\varphi(x) |\varphi\in
\Xi_{\iota}'^{(k)}\}=\lan h_k,\lambda\ran\geq \varepsilon_k^*(x)$ holds by (\ref{Btilde}). \qed

\section{Combinatorial objects}

In this section, we recall several combinatorial objects ((revised) extended Young diagrams, Young walls), which express explicit forms of ${\rm Im}(\Psi_{\iota}^{(\lambda)})$.

\subsection{Extended Young diagrams}\label{EYD-sub}

\begin{defn}\cite{Ha,JMMO}
Let $y_{\infty}$ be a fixed integer. A sequence $(y_r)_{r\in\mathbb{Z}_{\geq0}}$
is called an {\it extended Young diagram} of charge $y_{\infty}$ if the following conditions hold:
\begin{itemize}
\item $y_r\in\mathbb{Z}$, $y_r\leq y_{r+1}$ for all $r\in\mathbb{Z}_{\geq0}$,
\item $y_r=y_{\infty}$ for $r\gg0$.
\end{itemize}
\end{defn}

Each extended Young diagram is described as a Young diagram
drawn in $\mathbb{R}_{\geq0}\times \mathbb{R}_{\leq y_{\infty}}$ as follows:
For $(y_r)_{r\in\mathbb{Z}_{\geq0}}$, we draw a line between the points $(r,y_r)$ and $(r+1,y_r)$ and when
$y_r<y_{r+1}$ draw
a line
between $(r+1,y_r)$ and $(r+1,y_{r+1})$
for each $r\in\mathbb{Z}_{\geq0}$.
\begin{ex}\label{ex-1}
Let $T=(y_r)_{r\in\mathbb{Z}_{\geq0}}$ be
the extended Young diagram of charge
$y_{\infty}=1$ defined by
$y_0=-3$, $y_1=-2$, $y_2=y_3=-1$, $y_4=0$, $y_5=1, y_6=1,\cdots$.
Then $T$ is drawn in $\mathbb{R}_{\geq0}\times \mathbb{R}_{\leq 1}$ as
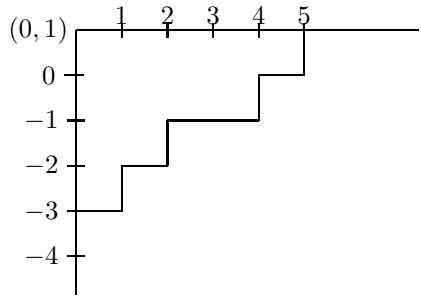
\[
\begin{xy}
(0,0) *{}="1",
(45,0)*{}="2",
(0,-35)*{}="3",
(-5,0)*{(0,1)}="4",
(6,2) *{1}="10",
(6,-1) *{}="1010",
(12,2) *{2}="11",
(12,-1) *{}="1111",
(18,2) *{3}="12",
(18,-1) *{}="1212",
(24,2) *{4}="13",
(24,-1) *{}="1313",
(30,2) *{5}="14",
(30,-1) *{}="1414",
(-3,-6)*{0\ }="5",
(1,-6)*{}="55",
(-4,-12)*{-1\ }="6",
(1,-12)*{}="66",
(-4,-18)*{-2\ }="7",
(1,-18)*{}="77",
(-4,-24)*{-3\ }="8",
(1,-24)*{}="88",
(6,-24)*{}="810",
(6,-18)*{}="8107",
(12,-18)*{}="810711",
(12,-12)*{}="8107116",
(24,-12)*{}="810711613",
(24,-6)*{}="8107116130",
(30,-6)*{}="81071161300",
(30,0)*{}="810711613000",
(-4,-30)*{-4\ }="9",
(1,-30)*{}="99",
\ar@{-} "1";"2"^{}
\ar@{-} "1";"3"^{}
\ar@{-} "5";"55"^{}
\ar@{-} "6";"66"^{}
\ar@{-} "7";"77"^{}
\ar@{-} "8";"810"^{}
\ar@{-} "810";"8107"^{}
\ar@{-} "8107";"810711"^{}
\ar@{-} "810711";"8107116"^{}
\ar@{-} "8107116";"810711613"^{}
\ar@{-} "810711613";"8107116130"^{}
\ar@{-} "8107116130";"81071161300"^{}
\ar@{-} "81071161300";"810711613000"^{}
\ar@{-} "9";"99"^{}
\ar@{-} "10";"1010"^{}
\ar@{-} "11";"1111"^{}
\ar@{-} "12";"1212"^{}
\ar@{-} "13";"1313"^{}
\ar@{-} "14";"1414"^{}
\end{xy}
\]
\end{ex}
\nd
Note that drawing an extended Young diagram
$(y_r)_{r\in\mathbb{Z}_{\geq0}}$ in
$\mathbb{R}_{\geq0}\times \mathbb{R}_{\leq y_{\infty}}$,
if $y_r< y_{r+1}$ then
the points $(r+1,y_r)$ and $(r+1,y_{r+1})$ are corners. 
\begin{defn}\cite{JMMO}
For an extended Young diagram
$(y_r)_{r\in\mathbb{Z}_{\geq0}}$, if $y_r< y_{r+1}$ then
the point $(r+1,y_r)$ is said to be a {\it convex corner} and $(r+1,y_{r+1})$
is said to be a {\it concave corner}. The point $(0,y_0)$ is also called a concave corner.
\end{defn}

\begin{ex}
Let $T$ be the
extended Young diagram as in Example \ref{ex-1}.
The points $(1,-3)$, $(2,-2)$, $(4,-1)$ and $(5,0)$ are convex corners and
$(0,-3)$, $(1,-2)$, $(2,-1)$, $(4,0)$ and $(5,1)$ are concave corners.
\end{ex}

Next, we consider a coloring of each corner in extended Young diagrams
by $I$.

\begin{defn}\cite{JMMO,KMM}\label{pi1-def}
\begin{enumerate}
\item
We define a map $\pi_{{\rm A}^{(1)}}:\mathbb{Z}\rightarrow \{1,2,\cdots,n\}=I$ as
\[
\pi_{{\rm A}^{(1)}}(\ell+rn)=\ell
\]
for any $r\in\mathbb{Z}$ and $\ell\in\{1,2,\cdots,n\}$.
\item
Defining a map
$\{1,2,\cdots,2n-2\}\rightarrow \{1,2,\cdots,n\}$
as
\[
\ell\mapsto \ell,\ 2n-\ell\mapsto \ell \qquad (2\leq \ell\leq n-1),
\]
\[
1\mapsto1,\ n\mapsto n
\]
and extend it to a map $\pi_{{\rm C}^{(1)}}:\mathbb{Z}\rightarrow\{1,2,\cdots,n\}=I$ by periodicity $2n-2$.
\item Defining a map $\{1,2,3,\cdots,2n-1\}\rightarrow \{1,2,\cdots,n\}$ as
\[
\ell\mapsto \ell,\ \ 2n-\ell\mapsto \ell \quad (1\leq \ell\leq n-1),\ \ n\mapsto n
\]
and extend it to the map $\pi_{{\rm A}^{(2)}}:\mathbb{Z}\rightarrow \{1,2,\cdots,n\}=I$ with periodicity $2n-1$.
\item Defining a map $\{1,2,3,\cdots,2n\}\rightarrow \{1,2,\cdots,n\}$ as
\[
\ell\mapsto \ell,\ \ 2n+1-\ell\mapsto \ell \quad (1\leq \ell\leq n),
\]
and extend it to the map $\pi_{{\rm D}^{(2)}}:\mathbb{Z}\rightarrow \{1,2,\cdots,n\}=I$ with periodicity $2n$.
\end{enumerate}
\end{defn}
Originally,
the maps $\pi_{{\rm A}^{(1)}}$, $\pi_{{\rm C}^{(1)}}$, $\pi_{{\rm A}^{(2)}}$ and $\pi_{{\rm D}^{(2)}}$ in this definition were
introduced to define actions of Chevalley generators (or Kashiwara operators) of type ${\rm A}^{(1)}_n$, ${\rm C}^{(1)}_n$, ${\rm A}^{(2)}_{2n}$
and ${\rm D}^{(2)}_{n+1}$ on extended Young diagrams. Each corner $(i,j)$ is colored by
$\pi_{{\rm A}^{(1)}}(i+j)$, $\pi_{{\rm C}^{(1)}}(i+j)$, $\pi_{{\rm A}^{(2)}}(i+j)$ and $\pi_{{\rm D}^{(2)}}(i+j)\in I$, respectively
and roughly speaking, a concave corner colored by $k\in I$
becomes a convex corner by an action of $f_k$ or $\tilde{f}_k$ \cite{JMMO, KMM}.

\subsection{Revised extended Young diagrams}

Next, we reviewed the notion of revised extended Young diagrams (REYD) introduced in \cite{Ka}, which are
used to express inequalities of ${\rm A}^{(2)}_{2n-2}$ and ${\rm C}^{(1)}_{n-1}$.

\begin{defn}\label{AEYD}\cite{Ka}
For $k\in I\setminus\{1\}$, we define ${\rm REYD}^{{\rm A}^{(2)}}_{k}$ as the set of sequences $(y_t)_{t\in\mathbb{Z}}$ such that
\begin{itemize}
\item[$(1)$] $y_t\in\mathbb{Z}$ for all $t\in\mathbb{Z}$,  
\item[$(2)$] $y_{t}=k$ for $t\gg0$ and $y_t=k+t$ for $t\ll0$,
\item[$(3)$] for $t\in\mathbb{Z}$ such that $k+t\not\equiv 0$ (mod $2n-1$), it holds either $y_{t+1}=y_t$ or $y_{t+1}=y_t+1$,
\item[$(4)$] for $t\in\mathbb{Z}_{>0}$ such that $k+t\equiv 0$ (mod $2n-1$), we have $y_{t+1}\geq y_t$,
\item[$(5)$] for $t\in\mathbb{Z}_{<0}$ such that $k+t\equiv 0$ (mod $2n-1$), we have $y_{t+1}\leq y_t+1$.
\end{itemize}
\end{defn}

\begin{defn}\label{DEYD}\cite{Ka}
For $k\in I\setminus\{1,n\}$, we define ${\rm REYD}^{{\rm D}^{(2)}}_{k}$ as the set of sequences $(y_t)_{t\in\mathbb{Z}}$ such that
\begin{itemize}
\item[$(1)$] $y_t\in\mathbb{Z}$ for all $t\in\mathbb{Z}$,  
\item[$(2)$] $y_{t}=k$ for $t\gg0$ and $y_t=k+t$ for $t\ll0$,
\item[$(3)$] for $t\in\mathbb{Z}$ such that $k+t\not\equiv 0, n$ (mod $2n$), it holds either $y_{t+1}=y_t$ or $y_{t+1}=y_t+1$,
\item[$(4)$] for $t\in\mathbb{Z}_{>0}$ such that $k+t\equiv 0$ or $n$ (mod $2n$), we have $y_{t+1}\geq y_t$,
\item[$(5)$] for $t\in\mathbb{Z}_{<0}$ such that $k+t\equiv 0$ or $n$ (mod $2n$), we have $y_{t+1}\leq y_t+1$.
\end{itemize}
\end{defn}
Just as in ordinary extended Young diagrams,
each element in ${\rm REYD}^{{\rm A}^{(2)}}_{k}$ and
${\rm REYD}^{{\rm D}^{(2)}}_{k}$
is described as an infinite diagram
drawn on $\mathbb{R}\times \mathbb{R}_{\leq k}$.
For example, 
let $n=3$, $k=2$ and $T=(y_t)_{t\in\mathbb{Z}}$ be
the element in ${\rm REYD}^{{\rm A}^{(2)}}_{2}$
defined as
\[
y_{\ell}=\ell+2\ (\ell\leq -3),\ y_{-2}=0,\ y_{-1}=y_0=y_1=-2,\ y_2=y_3=-1,\ y_4=y_5=y_6=y_7=1,\ y_t=2 (t\geq 8).
\]
Then $T$ is described as
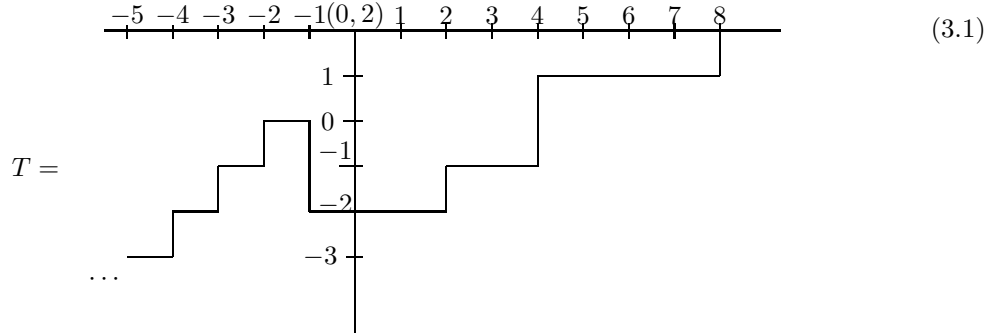
\begin{equation}\label{reydA-ex1}
\begin{xy}
(-42,-18) *{T=}="YY",
(-33,0) *{}="-6",
(-6,2) *{-1}="-1",
(-12,2) *{-2}="-2",
(-18,2) *{-3}="-3",
(-24,2) *{-4}="-4",
(-30,2) *{-5}="-5",
(-6,-1) *{}="-1a",
(-12,-1) *{}="-2a",
(-18,-1) *{}="-3a",
(-24,-1) *{}="-4a",
(-30,-1) *{}="-5a",
(0,0) *{}="1",
(56,0)*{}="2",
(0,-40)*{}="3",
(0,2)*{(0,2)}="4",
(6,2) *{1}="10",
(6,-1) *{}="1010",
(12,2) *{2}="11",
(12,-1) *{}="1111",
(18,2) *{3}="12",
(18,-1) *{}="1212",
(24,2) *{4}="13",
(24,-1) *{}="1313",
(30,2) *{5}="14",
(36,2) *{6}="6newnew",
(42,2) *{7}="7newnew",
(36,-1) *{}="6newnew-1",
(42,-1) *{}="7newnew-1",
(30,-1) *{}="1414",
(-3,-6)*{1\ }="5",
(1,-6)*{}="55",
(-3,-12)*{0\ }="6",
(-6,-12)*{}="6a",
(-12,-12)*{}="6aa",
(-12,-18)*{}="6aaa",
(-18,-18)*{}="6aaaa",
(-18,-24)*{}="st1",
(-24,-24)*{}="st2",
(-24,-30)*{}="st3",
(-30,-30)*{}="st4",
(-33,-33)*{\cdots}="stdot",
(1,-12)*{}="66",
(-6,-18)*{}="7",
(-6,-24)*{}="7new",
(-2,-16)*{-1\ }="7a",
(-2,-18)*{}="m1l",
(1,-18)*{}="m1r",
(-2,-23)*{-2\ }="8",
(1,-24)*{}="88",
(6,-18)*{}="8107",
(6,-24)*{}="8107new",
(12,-18)*{}="810711",
(12,-24)*{}="810711new",
(12,-12)*{}="8107116",
(24,-12)*{}="810711613",
(24,-18)*{}="810711613new",
(24,-6)*{}="8107116130",
(48,-6)*{}="8107116130new",
(48,2)*{8}="8107116130new2",
(30,-6)*{}="81071161300",
(30,0)*{}="810711613000",
(-4,-30)*{-3\ }="9",
(1,-30)*{}="99",
\ar@{-} "1";"-6"^{}
\ar@{-} "1";"2"^{}
\ar@{-} "1";"3"^{}
\ar@{-} "-1";"-1a"^{}
\ar@{-} "-2";"-2a"^{}
\ar@{-} "-3";"-3a"^{}
\ar@{-} "-4";"-4a"^{}
\ar@{-} "-5";"-5a"^{}
\ar@{-} "5";"55"^{}
\ar@{-} "6";"66"^{}
\ar@{-} "6newnew";"6newnew-1"^{}
\ar@{-} "7newnew";"7newnew-1"^{}
\ar@{-} "m1r";"m1l"^{}
\ar@{-} "810711";"810711613new"^{}
\ar@{-} "810711613";"810711613new"^{}
\ar@{-} "810711613";"8107116130"^{}
\ar@{-} "8107116130";"8107116130new"^{}
\ar@{-} "8107116130new2";"8107116130new"^{}
\ar@{-} "7new";"6a"^{}
\ar@{-} "7new";"810711new"^{}
\ar@{-} "810711";"810711new"^{}
\ar@{-} "6aa";"6a"^{}
\ar@{-} "6aaa";"6aa"^{}
\ar@{-} "6aaaa";"6aaa"^{}
\ar@{-} "st1";"6aaaa"^{}
\ar@{-} "st1";"st2"^{}
\ar@{-} "st2";"st3"^{}
\ar@{-} "st3";"st4"^{}
\ar@{-} "8";"88"^{}
\ar@{-} "9";"99"^{}
\ar@{-} "10";"1010"^{}
\ar@{-} "11";"1111"^{}
\ar@{-} "12";"1212"^{}
\ar@{-} "13";"1313"^{}
\ar@{-} "14";"1414"^{}
\end{xy}
\end{equation}

As an analogy of concave corners and convex corners of extended Young diagrams,
one can define admissible points and removable points of REYD:
\begin{defn}\label{ad-rem-pt}\cite{Ka}
We set ${\rm REYD}_{k}={\rm REYD}^{{\rm A}^{(2)}}_{k}$ or ${\rm REYD}_{k}={\rm REYD}^{{\rm D}^{(2)}}_{k}$ in Definition \ref{AEYD}, \ref{DEYD} .
Let $T=(y_t)_{t\in\mathbb{Z}}\in{\rm REYD}_{k}$ and $i\in\mathbb{Z}$.
\begin{enumerate}
\item Let $T'=(y_t')_{t\in\mathbb{Z}}$ be the sequence of integers such that
$y_i'=y_i-1$ and $y_t'=y_t$ $(t\neq i)$. 
The point $(i,y_i)$ is said to be an admissible point of $T$ if $T'\in{\rm REYD}_{k}$.
\item Let $T''=(y_t'')_{t\in\mathbb{Z}}$ be the sequence of integers such that
$y_{i-1}''=y_{i-1}+1$ and $y_t''=y_t$ $(t\neq i-1)$. 
The point $(i,y_{i-1})$ is said to be 
a removable point of $T$
if $T''\in{\rm REYD}_{k}$.
\end{enumerate}
\end{defn}

\begin{defn}\label{ad-rem-pt2}\cite{Ka}
Let $T=(y_t)_{t\in\mathbb{Z}}\in{\rm REYD}^{{\rm A}^{(2)}}_{k}$ and $i\in\mathbb{Z}$. 
\begin{enumerate}
\item 
The point $(i,y_i)$ is said to be a double $1$-admissible point of $T$
if $y_{i-1}<y_i=y_{i+1}$ and it holds either 
\begin{itemize}
\item $i+k\equiv 1$ (mod $2n-1$) and $i<0$ or
\item $i+k\equiv 0$ (mod $2n-1$) and $i>0$.
\end{itemize}
\item 
The point $(i,y_{i-1})$ is said to be a double $1$-removable point of $T$
if $y_{i-2}=y_{i-1}<y_{i}$ and it holds
either
\begin{itemize}
\item
$i+k-1\equiv 1$ (mod $2n-1$) and $i>1$ or
\item
$i+k-1\equiv 0$ (mod $2n-1$) and $i<1$.
\end{itemize} 
\item Other admissible (resp. removable) points $(i,y_i)$ (resp. $(i,y_{i-1})$) than (i) (resp. (ii)) are said to be single $\pi_{{\rm A}^{(2)}}(i+k)$-admissible (resp. $\pi_{{\rm A}^{(2)}}(i+k-1)$-removable) points.
\end{enumerate}
\end{defn}

\begin{defn}\label{ad-rem-pt3}\cite{Ka}
Let $T=(y_t)_{t\in\mathbb{Z}}\in{\rm REYD}^{{\rm D}^{(2)}}_{k}$, $i\in\mathbb{Z}$
and $\ell\in\{0,n\}$.
\begin{enumerate}
\item The point $(i,y_i)$ is said to be a double $\pi_{{\rm D}^{(2)}}(\ell)$-admissible point of $T$ 
if $y_{i-1}<y_i=y_{i+1}$ and it holds either 
\begin{itemize}
\item
$i+k\equiv \ell+1$ (mod $2n$) and $i<0$ or 
\item
$i+k\equiv \ell$ (mod $2n$) and $i>0$.
\end{itemize}
\item The point $(i,y_{i-1})$ is said to be a double $\pi_{{\rm D}^{(2)}}(\ell)$-removable point of $T$
if $y_{i-2}=y_{i-1}<y_{i}$ and it holds
either
\begin{itemize}
\item
$i+k-1\equiv \ell+1$ (mod $2n$) and $i>1$ or 
\item 
$i+k-1\equiv \ell$ (mod $2n$) and $i<1$. 
\end{itemize}
\item Other admissible (resp. removable) points $(i,y_i)$ (resp. $(i,y_{i-1})$) than (i) (resp. (ii)) are said to be single $\pi_{{\rm D}^{(2)}}(i+k)$-admissible (resp. $\pi_{{\rm D}^{(2)}}(i+k-1)$-removable) points.
\end{enumerate}
\end{defn}

Note that in Definition \ref{ad-rem-pt2} and \ref{ad-rem-pt3}, a double admissible (resp. removable) point 
$(i,y_i)$ (resp. $(i,y_{i-1})$) is an admissible (resp. a removable) point in Definition \ref{ad-rem-pt} by $y_{i-1}<y_i=y_{i+1}$
(resp. $y_{i-2}=y_{i-1}<y_{i}$).

\begin{ex}
Let us consider the revised extended Young diagram $T$ in (\ref{reydA-ex1}).
The points $(-2,0)$, $(-1,-2)$ and $(4,1)$ are single $1$-admissible points in $T$.
The point $(2,-1)$ is a single $2$-admissible point and
the point $(8,2)$ is a double $1$-admissible point in $T$.
We also see that the point $(2,-2)$ is a single $3$-removable point in $T$.
The point $(4,-1)$ is a single $1$-removable point and
the point $(8,1)$ is a single $2$-removable point in $T$.

\end{ex}

\subsection{Young walls}

Following \cite{Kang}, let us review the notion of Young walls of type ${\rm A}^{(2)}_{2n-2}$ and ${\rm D}^{(2)}_n$. 
We consider
$I$-colored blocks of two different
shapes:
\begin{enumerate}
\item[(1)] block with unit width, unit height and unit thickness:
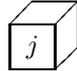
\[
\begin{xy}
(3,3) *{j}="0",
(0,0) *{}="1",
(6,0)*{}="2",
(6,6)*{}="3",
(0,6)*{}="4",
(3,9)*{}="5",
(9,9)*{}="6",
(9,3)*{}="7",
\ar@{-} "1";"2"^{}
\ar@{-} "1";"4"^{}
\ar@{-} "2";"3"^{}
\ar@{-} "3";"4"^{}
\ar@{-} "5";"4"^{}
\ar@{-} "5";"6"^{}
\ar@{-} "3";"6"^{}
\ar@{-} "2";"7"^{}
\ar@{-} "6";"7"^{}
\end{xy}
\]
\item[(2)] block with unit width, half-unit height and unit thickness:
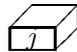
\[
\begin{xy}
(3,1.5) *{j}="0",
(0,0) *{}="1",
(6,0)*{}="2",
(6,3)*{}="3",
(0,3)*{}="4",
(3,6)*{}="5",
(9,6)*{}="6",
(9,3)*{}="7",
\ar@{-} "1";"2"^{}
\ar@{-} "1";"4"^{}
\ar@{-} "2";"3"^{}
\ar@{-} "3";"4"^{}
\ar@{-} "5";"4"^{}
\ar@{-} "5";"6"^{}
\ar@{-} "3";"6"^{}
\ar@{-} "2";"7"^{}
\ar@{-} "6";"7"^{}
\end{xy}
\]
\end{enumerate}
The block (1) with color $j\in I$ is simply expressed as
\begin{equation}\label{smpl1}
\begin{xy}
(3,3) *{j}="0",
(0,0) *{}="1",
(6,0)*{}="2",
(6,6)*{}="3",
(0,6)*{}="4",
\ar@{-} "1";"2"^{}
\ar@{-} "1";"4"^{}
\ar@{-} "2";"3"^{}
\ar@{-} "3";"4"^{}
\end{xy}
\end{equation}
and (2) with color $j\in I$ is expressed as
\begin{equation}\label{smpl2}
\begin{xy}
(1.5,1.5) *{\ \ j}="0",
(0,0) *{}="1",
(6,0)*{}="2",
(6,3.5)*{}="3",
(0,3.5)*{}="4",
\ar@{-} "1";"2"^{}
\ar@{-} "1";"4"^{}
\ar@{-} "2";"3"^{}
\ar@{-} "3";"4"^{}
\end{xy}
\end{equation}
We call
the blocks (\ref{smpl1}) and (\ref{smpl2}) $j$-blocks.
For example, if colored blocks are stacked as follows
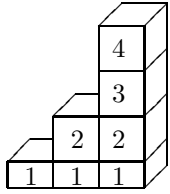
\[
\begin{xy}
(-9.5,1.5) *{\ 1}="000",
(-3.5,1.5) *{\ 1}="00",
(1.5,1.5) *{\ \ 1}="0",
(1.5,6.5) *{\ \ 2}="02",
(-4,6.5) *{\ \ 2}="002",
(1.5,12.5) *{\ \ 3}="03",
(1.5,18.5) *{\ \ 4}="04",
(9,3)*{}="2-a",
(9,6.5)*{}="3-a",
(9,12.5)*{}="3-1-a",
(9,18.5)*{}="3-2-a",
(9,24.5)*{}="3-3-a",
(3,24.5)*{}="4-3-a",
(-3,12.5)*{}="5-1-a",
(0,12.5)*{}="5-1-ab",
(-9,6.5)*{}="7-a",
(-6,6.5)*{}="7-ab",
(0,0) *{}="1",
(6,0)*{}="2",
(6,3.5)*{}="3",
(0,3.5)*{}="4",
(-6,3.5)*{}="5",
(-6,0)*{}="6",
(-12,3.5)*{}="7",
(-12,0)*{}="8",
(6,9.5)*{}="3-1",
(6,15.5)*{}="3-2",
(6,21.5)*{}="3-3",
(0,9.5)*{}="4-1",
(0,15.5)*{}="4-2",
(0,21.5)*{}="4-3",
(-6,9.5)*{}="5-1",
\ar@{-} "1";"2"^{}
\ar@{-} "1";"4"^{}
\ar@{-} "2";"3"^{}
\ar@{-} "3";"4"^{}
\ar@{-} "5";"6"^{}
\ar@{-} "5";"4"^{}
\ar@{-} "1";"6"^{}
\ar@{-} "7";"8"^{}
\ar@{-} "7";"5"^{}
\ar@{-} "6";"8"^{}
\ar@{-} "3";"3-1"^{}
\ar@{-} "4-1";"3-1"^{}
\ar@{-} "4-1";"4"^{}
\ar@{-} "5-1";"5"^{}
\ar@{-} "4-1";"5-1"^{}
\ar@{-} "4-1";"4-2"^{}
\ar@{-} "3-2";"3-1"^{}
\ar@{-} "3-2";"4-2"^{}
\ar@{-} "3-2";"3-3"^{}
\ar@{-} "4-3";"4-2"^{}
\ar@{-} "4-3";"3-3"^{}
\ar@{-} "2";"2-a"^{}
\ar@{-} "3";"3-a"^{}
\ar@{-} "3-1";"3-1-a"^{}
\ar@{-} "3-2";"3-2-a"^{}
\ar@{-} "3-3";"3-3-a"^{}
\ar@{-} "2-a";"3-3-a"^{}
\ar@{-} "4-3";"4-3-a"^{}
\ar@{-} "3-3-a";"4-3-a"^{}
\ar@{-} "5-1";"5-1-a"^{}
\ar@{-} "5-1-ab";"5-1-a"^{}
\ar@{-} "7";"7-a"^{}
\ar@{-} "7-ab";"7-a"^{}
\end{xy}
\]
then
it is simply expressed as
\begin{equation}\label{smpl3}
\begin{xy}
(-9.5,1.5) *{\ 1}="000",
(-3.5,1.5) *{\ 1}="00",
(1.5,1.5) *{\ \ 1}="0",
(1.5,6.5) *{\ \ 2}="02",
(-4,6.5) *{\ \ 2}="002",
(1.5,12.5) *{\ \ 3}="03",
(1.5,18.5) *{\ \ 4}="04",
(0,0) *{}="1",
(6,0)*{}="2",
(6,3.5)*{}="3",
(0,3.5)*{}="4",
(-6,3.5)*{}="5",
(-6,0)*{}="6",
(-12,3.5)*{}="7",
(-12,0)*{}="8",
(6,9.5)*{}="3-1",
(6,15.5)*{}="3-2",
(6,21.5)*{}="3-3",
(0,9.5)*{}="4-1",
(0,15.5)*{}="4-2",
(0,21.5)*{}="4-3",
(-6,9.5)*{}="5-1",
\ar@{-} "1";"2"^{}
\ar@{-} "1";"4"^{}
\ar@{-} "2";"3"^{}
\ar@{-} "3";"4"^{}
\ar@{-} "5";"6"^{}
\ar@{-} "5";"4"^{}
\ar@{-} "1";"6"^{}
\ar@{-} "7";"8"^{}
\ar@{-} "7";"5"^{}
\ar@{-} "6";"8"^{}
\ar@{-} "3";"3-1"^{}
\ar@{-} "4-1";"3-1"^{}
\ar@{-} "4-1";"4"^{}
\ar@{-} "5-1";"5"^{}
\ar@{-} "4-1";"5-1"^{}
\ar@{-} "4-1";"4-2"^{}
\ar@{-} "3-2";"3-1"^{}
\ar@{-} "3-2";"4-2"^{}
\ar@{-} "3-2";"3-3"^{}
\ar@{-} "4-3";"4-2"^{}
\ar@{-} "4-3";"3-3"^{}
\end{xy}
\end{equation}
A diagram formed by stacking blocks is called a wall.

For $X={\rm A}^{(2)}_{2n-2}$, ${\rm D}^{(2)}_{n}$
and $\lambda\in P^+$ of level $1$,
we define a wall $Y_{\lambda}$ called a {\it ground state wall}:

\begin{itemize}
\item For $X={\rm A}^{(2)}_{2n-2}$ and $\lambda=\Lambda_1$, one defines
\[
Y_{\Lambda_1}=
\begin{xy}
(-15.5,1.5) *{\ \cdots}="0000",
(-9.5,1.5) *{\ 1}="000",
(-3.5,1.5) *{\ 1}="00",
(1.5,1.5) *{\ \ 1}="0",
(0,0) *{}="1",
(6,0)*{}="2",
(6,3.5)*{}="3",
(0,3.5)*{}="4",
(-6,3.5)*{}="5",
(-6,0)*{}="6",
(-12,3.5)*{}="7",
(-12,0)*{}="8",
\ar@{-} "1";"2"^{}
\ar@{-} "1";"4"^{}
\ar@{-} "2";"3"^{}
\ar@{-} "3";"4"^{}
\ar@{-} "5";"6"^{}
\ar@{-} "5";"4"^{}
\ar@{-} "1";"6"^{}
\ar@{-} "7";"8"^{}
\ar@{-} "7";"5"^{}
\ar@{-} "6";"8"^{}
\end{xy}
\]
Here, the wall $Y_{\Lambda_1}$ has infinitely many $1$-blocks of half-unit height
and
extends infinitely to the left.
\item
For $X={\rm D}^{(2)}_{n}$ and $\lambda=\Lambda_1$, $\Lambda_n$, one defines
\[
Y_{\Lambda_1}=
\begin{xy}
(-15.5,1.5) *{\ \cdots}="0000",
(-9.5,1.5) *{\ 1}="000",
(-3.5,1.5) *{\ 1}="00",
(1.5,1.5) *{\ \ 1}="0",
(0,0) *{}="1",
(6,0)*{}="2",
(6,3.5)*{}="3",
(0,3.5)*{}="4",
(-6,3.5)*{}="5",
(-6,0)*{}="6",
(-12,3.5)*{}="7",
(-12,0)*{}="8",
\ar@{-} "1";"2"^{}
\ar@{-} "1";"4"^{}
\ar@{-} "2";"3"^{}
\ar@{-} "3";"4"^{}
\ar@{-} "5";"6"^{}
\ar@{-} "5";"4"^{}
\ar@{-} "1";"6"^{}
\ar@{-} "7";"8"^{}
\ar@{-} "7";"5"^{}
\ar@{-} "6";"8"^{}
\end{xy}
\]
and
\[
Y_{\Lambda_n}=
\begin{xy}
(-15.5,1.5) *{\ \cdots}="0000",
(-9.5,1.5) *{\ n}="000",
(-3.5,1.5) *{\ n}="00",
(1.5,1.5) *{\ \ n}="0",
(0,0) *{}="1",
(6,0)*{}="2",
(6,3.5)*{}="3",
(0,3.5)*{}="4",
(-6,3.5)*{}="5",
(-6,0)*{}="6",
(-12,3.5)*{}="7",
(-12,0)*{}="8",
\ar@{-} "1";"2"^{}
\ar@{-} "1";"4"^{}
\ar@{-} "2";"3"^{}
\ar@{-} "3";"4"^{}
\ar@{-} "5";"6"^{}
\ar@{-} "5";"4"^{}
\ar@{-} "1";"6"^{}
\ar@{-} "7";"8"^{}
\ar@{-} "7";"5"^{}
\ar@{-} "6";"8"^{}
\end{xy}
\]
\end{itemize}

\begin{defn}\cite{Kang}\label{def-YW}
For $X={\rm A}^{(2)}_{2n-2}$, ${\rm D}^{(2)}_{n}$
and $\lambda\in P^+$ of level $1$,
a wall $Y$ is called a {\it Young wall} of ground state $\lambda$
of type $X$ if it satisfies the following:
\begin{enumerate}
\item The wall $Y$ is obtained
by stacking finitely many colored blocks
on the ground state wall $Y_{\lambda}$.
\item The colored blocks are stacked following the patterns we give below for each type and $\lambda$.
\item Let $h_j$ be the height of $j$-th column of the wall $Y$ from the right. Then we get $h_j\geq h_{j+1}$
for all $j\in\mathbb{Z}_{\geq1}$. 
\end{enumerate}
The patterns mentioned in (ii) are as follows: 
In the case ${\rm A}^{(2)}_{2n-2}$ and $\lambda=\Lambda_1$:
\[
\begin{xy}
(-15.5,-2) *{\ 1}="000-1",
(-9.5,-2) *{\ 1}="00-1",
(-3.5,-2) *{\ 1}="0-1",
(1.5,-2) *{\ \ 1}="0-1",
(-15.5,1.5) *{\ 1}="0000",
(-9.5,1.5) *{\ 1}="000",
(-3.5,1.5) *{\ 1}="00",
(1.5,1.5) *{\ \ 1}="0",
(1.5,6.5) *{\ \ 2}="02",
(-4,6.5) *{\ \ 2}="002",
(-10,6.5) *{\ \ 2}="0002",
(-16,6.5) *{\ \ 2}="00002",
(1.5,12.5) *{\ \ 3}="03",
(-4,12.5) *{\ \ 3}="003",
(-10,12.5) *{\ \ 3}="0003",
(-16,12.5) *{\ \ 3}="00003",
(1.5,20.5) *{\ \ \vdots}="04",
(-4,20.5) *{\ \ \vdots}="004",
(-10,20.5) *{\ \ \vdots}="0004",
(-16,20.5) *{\ \ \vdots}="00004",
(1.5,27.5) *{\ \ _{n-1}}="05",
(-4,27.5) *{\ \ _{n-1}}="005",
(-10,27.5) *{\ \ _{n-1}}="0005",
(-16,27.5) *{\ \ _{n-1}}="00005",
(1.5,33.5) *{\ \ n}="06",
(-4,33.5) *{\ \ n}="006",
(-25,34) *{\cdots}="dots",
(-10,33.5) *{\ \ n}="0006",
(-16,33.5) *{\ \ n}="00006",
(1.5,39.5) *{\ \ _{n-1}}="07",
(-4,39.5) *{\ \ _{n-1}}="007",
(-10,39.5) *{\ \ _{n-1}}="0007",
(-16,39.5) *{\ \ _{n-1}}="00007",
(1.5,48) *{\ \ \vdots}="08",
(-4,48) *{\ \ \vdots}="008",
(-10,48) *{\ \ \vdots}="0008",
(-16,48) *{\ \ \vdots}="00008",
(1.5,55) *{\ \ 2}="09",
(-4,55) *{\ \ 2}="009",
(-10,55) *{\ \ 2}="0009",
(-16,55) *{\ \ 2}="00009",
(1.5,60) *{\ \ 1}="010",
(-4,60) *{\ \ 1}="0010",
(-10,60) *{\ \ 1}="00010",
(-16,60) *{\ \ 1}="000010",
(1.5,64) *{\ \ 1}="011",
(-4,64) *{\ \ 1}="0011",
(-10,64) *{\ \ 1}="00011",
(-16,64) *{\ \ 1}="000011",
(1.5,69) *{\ \ 2}="012",
(-4,69) *{\ \ 2}="0012",
(-10,69) *{\ \ 2}="00012",
(-16,69) *{\ \ 2}="000012",
(-18,-3.5) *{}="0.5-l",
(-20,-3.5) *{}="0-l",
(-20,0) *{}="1-l",
(-20,3.5) *{}="2-l",
(-20,9.5) *{}="3-l",
(-20,15.5) *{}="4-l",
(-20,24.5) *{}="5-l",
(-20,30.5) *{}="6-l",
(-20,36.5) *{}="7-l",
(-20,42.5) *{}="8-l",
(-20,52) *{}="9-l",
(-20,58) *{}="10-l",
(-20,62) *{}="11-l",
(-20,66) *{}="12-l",
(-20,72) *{}="13-l",
(-18,74) *{}="13.5-l",
(6,-3.5)*{}="r--1",
(6,0)*{}="r",
(6,3.5)*{}="r-0",
(0,-3.5) *{}="b1",
(-6,-3.5)*{}="b2",
(-12,-3.5)*{}="b3",
(0,74) *{}="t1",
(-6,74)*{}="t2",
(-12,74)*{}="t3",
(6,9.5)*{}="r-1",
(6,15.5)*{}="r-2",
(6,24.5)*{}="r-3",
(6,30.5)*{}="r-4",
(6,36.5)*{}="r-5",
(6,42.5)*{}="r-6",
(6,52)*{}="r-7",
(6,58)*{}="r-8",
(6,62)*{}="r-9",
(6,66)*{}="r-10",
(6,72)*{}="r-11",
(6,74)*{}="r-11.5",
\ar@{-} "b1";"t1"^{}
\ar@{-} "b2";"t2"^{}
\ar@{-} "b3";"t3"^{}
\ar@{-} "0.5-l";"13.5-l"^{}
\ar@{-} "r-11.5";"r--1"^{}
\ar@{-} "0-l";"r--1"^{}
\ar@{-} "1-l";"r"^{}
\ar@{-} "2-l";"r-0"^{}
\ar@{-} "3-l";"r-1"^{}
\ar@{-} "4-l";"r-2"^{}
\ar@{-} "5-l";"r-3"^{}
\ar@{-} "6-l";"r-4"^{}
\ar@{-} "7-l";"r-5"^{}
\ar@{-} "8-l";"r-6"^{}
\ar@{-} "9-l";"r-7"^{}
\ar@{-} "10-l";"r-8"^{}
\ar@{-} "11-l";"r-9"^{}
\ar@{-} "12-l";"r-10"^{}
\ar@{-} "13-l";"r-11"^{}
\end{xy}
\]
In the case $X={\rm D}^{(2)}_{n}$,
\[
\begin{xy}
(-40,35) *{\lambda=\Lambda_1:}="weight",
(-15.5,-2) *{\ 1}="000-1",
(-9.5,-2) *{\ 1}="00-1",
(-3.5,-2) *{\ 1}="0-1",
(1.5,-2) *{\ \ 1}="0-1",
(-15.5,1.5) *{\ 1}="0000",
(-9.5,1.5) *{\ 1}="000",
(-3.5,1.5) *{\ 1}="00",
(1.5,1.5) *{\ \ 1}="0",
(1.5,6.5) *{\ \ 2}="02",
(-4,6.5) *{\ \ 2}="002",
(-10,6.5) *{\ \ 2}="0002",
(-16,6.5) *{\ \ 2}="00002",
(1.5,12.5) *{\ \ 3}="03",
(-4,12.5) *{\ \ 3}="003",
(-10,12.5) *{\ \ 3}="0003",
(-16,12.5) *{\ \ 3}="00003",
(1.5,20.5) *{\ \ \vdots}="04",
(-4,20.5) *{\ \ \vdots}="004",
(-10,20.5) *{\ \ \vdots}="0004",
(-16,20.5) *{\ \ \vdots}="00004",
(1.5,27.5) *{\ \ _{n-1}}="05",
(-4,27.5) *{\ \ _{n-1}}="005",
(-10,27.5) *{\ \ _{n-1}}="0005",
(-16,27.5) *{\ \ _{n-1}}="00005",
(1.5,32) *{\ \ _n}="06",
(-4,32) *{\ \ _n}="006",
(-25,34) *{\cdots}="dots",
(-10,32) *{\ \ _n}="0006",
(-16,32) *{\ \ _n}="00006",
(1.5,35) *{\ \ _n}="06.5",
(-4,35) *{\ \ _n}="006.5",
(-10,35) *{\ \ _n}="0006.5",
(-16,35) *{\ \ _n}="00006.5",
(1.5,39.5) *{\ \ _{n-1}}="07",
(-4,39.5) *{\ \ _{n-1}}="007",
(-10,39.5) *{\ \ _{n-1}}="0007",
(-16,39.5) *{\ \ _{n-1}}="00007",
(1.5,48) *{\ \ \vdots}="08",
(-4,48) *{\ \ \vdots}="008",
(-10,48) *{\ \ \vdots}="0008",
(-16,48) *{\ \ \vdots}="00008",
(1.5,55) *{\ \ 2}="09",
(-4,55) *{\ \ 2}="009",
(-10,55) *{\ \ 2}="0009",
(-16,55) *{\ \ 2}="00009",
(1.5,60) *{\ \ 1}="010",
(-4,60) *{\ \ 1}="0010",
(-10,60) *{\ \ 1}="00010",
(-16,60) *{\ \ 1}="000010",
(1.5,64) *{\ \ 1}="011",
(-4,64) *{\ \ 1}="0011",
(-10,64) *{\ \ 1}="00011",
(-16,64) *{\ \ 1}="000011",
(1.5,69) *{\ \ 2}="012",
(-4,69) *{\ \ 2}="0012",
(-10,69) *{\ \ 2}="00012",
(-16,69) *{\ \ 2}="000012",
(-18,-3.5) *{}="0.5-l",
(-20,-3.5) *{}="0-l",
(-20,0) *{}="1-l",
(-20,3.5) *{}="2-l",
(-20,9.5) *{}="3-l",
(-20,15.5) *{}="4-l",
(-20,24.5) *{}="5-l",
(-20,30.5) *{}="6-l",
(-20,33.5) *{}="65-l",
(-20,36.5) *{}="7-l",
(-20,42.5) *{}="8-l",
(-20,52) *{}="9-l",
(-20,58) *{}="10-l",
(-20,62) *{}="11-l",
(-20,66) *{}="12-l",
(-20,72) *{}="13-l",
(-18,74) *{}="13.5-l",
(6,-3.5)*{}="r--1",
(6,0)*{}="r",
(6,3.5)*{}="r-0",
(0,-3.5) *{}="b1",
(-6,-3.5)*{}="b2",
(-12,-3.5)*{}="b3",
(0,74) *{}="t1",
(-6,74)*{}="t2",
(-12,74)*{}="t3",
(6,9.5)*{}="r-1",
(6,15.5)*{}="r-2",
(6,24.5)*{}="r-3",
(6,30.5)*{}="r-4",
(6,33.5)*{}="r-45",
(6,36.5)*{}="r-5",
(6,42.5)*{}="r-6",
(6,52)*{}="r-7",
(6,58)*{}="r-8",
(6,62)*{}="r-9",
(6,66)*{}="r-10",
(6,72)*{}="r-11",
(6,74)*{}="r-11.5",
\ar@{-} "b1";"t1"^{}
\ar@{-} "b2";"t2"^{}
\ar@{-} "b3";"t3"^{}
\ar@{-} "0.5-l";"13.5-l"^{}
\ar@{-} "r-11.5";"r--1"^{}
\ar@{-} "0-l";"r--1"^{}
\ar@{-} "1-l";"r"^{}
\ar@{-} "2-l";"r-0"^{}
\ar@{-} "3-l";"r-1"^{}
\ar@{-} "4-l";"r-2"^{}
\ar@{-} "5-l";"r-3"^{}
\ar@{-} "6-l";"r-4"^{}
\ar@{-} "65-l";"r-45"^{}
\ar@{-} "7-l";"r-5"^{}
\ar@{-} "8-l";"r-6"^{}
\ar@{-} "9-l";"r-7"^{}
\ar@{-} "10-l";"r-8"^{}
\ar@{-} "11-l";"r-9"^{}
\ar@{-} "12-l";"r-10"^{}
\ar@{-} "13-l";"r-11"^{}
\end{xy}\qquad \qquad
\begin{xy}
(-40,35) *{\lambda=\Lambda_n:}="weight",
(-15.5,-2) *{\ _n}="000-1",
(-9.5,-2) *{\ _n}="00-1",
(-3.5,-2) *{\ _n}="0-1",
(1.5,-2) *{\ \ _n}="0-1",
(-15.5,1.5) *{\ _n}="0000",
(-9.5,1.5) *{\ _n}="000",
(-3.5,1.5) *{\ _n}="00",
(1.5,1.5) *{\ \ _n}="0",
(1.5,6.5) *{\ \ _{n-1}}="02",
(-4,6.5) *{\ \ _{n-1}}="002",
(-10,6.5) *{\ \ _{n-1}}="0002",
(-16,6.5) *{\ \ _{n-1}}="00002",
(1.5,12.5) *{\ \ _{n-2}}="03",
(-4,12.5) *{\ \ _{n-2}}="003",
(-10,12.5) *{\ \ _{n-2}}="0003",
(-16,12.5) *{\ \ _{n-2}}="00003",
(1.5,20.5) *{\ \ \vdots}="04",
(-4,20.5) *{\ \ \vdots}="004",
(-10,20.5) *{\ \ \vdots}="0004",
(-16,20.5) *{\ \ \vdots}="00004",
(1.5,27.5) *{\ \ 2}="05",
(-4,27.5) *{\ \ 2}="005",
(-10,27.5) *{\ \ 2}="0005",
(-16,27.5) *{\ \ 2}="00005",
(1.5,32) *{\ \ 1}="06",
(-4,32) *{\ \ 1}="006",
(-25,34) *{\cdots}="dots",
(-10,32) *{\ \ 1}="0006",
(-16,32) *{\ \ 1}="00006",
(1.5,35) *{\ \ 1}="06.5",
(-4,35) *{\ \ 1}="006.5",
(-10,35) *{\ \ 1}="0006.5",
(-16,35) *{\ \ 1}="00006.5",
(1.5,39.5) *{\ \ 2}="07",
(-4,39.5) *{\ \ 2}="007",
(-10,39.5) *{\ \ 2}="0007",
(-16,39.5) *{\ \ 2}="00007",
(1.5,48) *{\ \ \vdots}="08",
(-4,48) *{\ \ \vdots}="008",
(-10,48) *{\ \ \vdots}="0008",
(-16,48) *{\ \ \vdots}="00008",
(1.5,55) *{\ \ _{n-1}}="09",
(-4,55) *{\ \ _{n-1}}="009",
(-10,55) *{\ \ _{n-1}}="0009",
(-16,55) *{\ \ _{n-1}}="00009",
(1.5,60) *{\ \ _{n}}="010",
(-4,60) *{\ \ _{n}}="0010",
(-10,60) *{\ \ _{n}}="00010",
(-16,60) *{\ \ _{n}}="000010",
(1.5,64) *{\ \ _{n}}="011",
(-4,64) *{\ \ _{n}}="0011",
(-10,64) *{\ \ _{n}}="00011",
(-16,64) *{\ \ _{n}}="000011",
(1.5,69) *{\ \ _{n-1}}="012",
(-4,69) *{\ \ _{n-1}}="0012",
(-10,69) *{\ \ _{n-1}}="00012",
(-16,69) *{\ \ _{n-1}}="000012",
(-18,-3.5) *{}="0.5-l",
(-20,-3.5) *{}="0-l",
(-20,0) *{}="1-l",
(-20,3.5) *{}="2-l",
(-20,9.5) *{}="3-l",
(-20,15.5) *{}="4-l",
(-20,24.5) *{}="5-l",
(-20,30.5) *{}="6-l",
(-20,33.5) *{}="65-l",
(-20,36.5) *{}="7-l",
(-20,42.5) *{}="8-l",
(-20,52) *{}="9-l",
(-20,58) *{}="10-l",
(-20,62) *{}="11-l",
(-20,66) *{}="12-l",
(-20,72) *{}="13-l",
(-18,74) *{}="13.5-l",
(6,-3.5)*{}="r--1",
(6,0)*{}="r",
(6,3.5)*{}="r-0",
(0,-3.5) *{}="b1",
(-6,-3.5)*{}="b2",
(-12,-3.5)*{}="b3",
(0,74) *{}="t1",
(-6,74)*{}="t2",
(-12,74)*{}="t3",
(6,9.5)*{}="r-1",
(6,15.5)*{}="r-2",
(6,24.5)*{}="r-3",
(6,30.5)*{}="r-4",
(6,33.5)*{}="r-45",
(6,36.5)*{}="r-5",
(6,42.5)*{}="r-6",
(6,52)*{}="r-7",
(6,58)*{}="r-8",
(6,62)*{}="r-9",
(6,66)*{}="r-10",
(6,72)*{}="r-11",
(6,74)*{}="r-11.5",
\ar@{-} "b1";"t1"^{}
\ar@{-} "b2";"t2"^{}
\ar@{-} "b3";"t3"^{}
\ar@{-} "0.5-l";"13.5-l"^{}
\ar@{-} "r-11.5";"r--1"^{}
\ar@{-} "0-l";"r--1"^{}
\ar@{-} "1-l";"r"^{}
\ar@{-} "2-l";"r-0"^{}
\ar@{-} "3-l";"r-1"^{}
\ar@{-} "4-l";"r-2"^{}
\ar@{-} "5-l";"r-3"^{}
\ar@{-} "6-l";"r-4"^{}
\ar@{-} "65-l";"r-45"^{}
\ar@{-} "7-l";"r-5"^{}
\ar@{-} "8-l";"r-6"^{}
\ar@{-} "9-l";"r-7"^{}
\ar@{-} "10-l";"r-8"^{}
\ar@{-} "11-l";"r-9"^{}
\ar@{-} "12-l";"r-10"^{}
\ar@{-} "13-l";"r-11"^{}
\end{xy}
\]
Note that the first row of each pattern from the bottom is
the ground state wall.
\end{defn}

\begin{ex}\label{ex-YW}
The following is a Young wall of ground state $\Lambda_1$ of type ${\rm A}^{(2)}_{4}$:
\begin{equation*}
\begin{xy}
(-15.5,-2) *{\ 1}="000-1",
(-9.5,-2) *{\ 1}="00-1",
(-3.5,-2) *{\ 1}="0-1",
(1.5,-2) *{\ \ 1}="0-1",
(-21.5,-2) *{\dots}="00000",
(-9.5,1.5) *{\ 1}="000",
(-3.5,1.5) *{\ 1}="00",
(1.5,1.5) *{\ \ 1}="0",
(1.5,6.5) *{\ \ 2}="02",
(-4,6.5) *{\ \ 2}="002",
(1.5,12.5) *{\ \ 3}="03",
(1.5,18.5) *{\ \ 2}="04",
(0,0) *{}="1",
(0,-3.5) *{}="1-u",
(6,0)*{}="2",
(6,-3.5)*{}="2-u",
(6,3.5)*{}="3",
(0,3.5)*{}="4",
(-6,3.5)*{}="5",
(-6,0)*{}="6",
(-6,-3.5)*{}="6-u",
(-12,3.5)*{}="7",
(-12,0)*{}="8",
(-12,-3.6)*{}="8-u",
(-18,3.5)*{}="9",
(-18,0)*{}="10-a",
(-18,0)*{}="10",
(-18,-3.5)*{}="10-u",
(6,9.5)*{}="3-1",
(6,15.5)*{}="3-2",
(6,21.5)*{}="3-3",
(0,9.5)*{}="4-1",
(0,15.5)*{}="4-2",
(0,21.5)*{}="4-3",
(-6,9.5)*{}="5-1",
\ar@{-} "8";"10-a"^{}
\ar@{-} "10-u";"10-a"^{}
\ar@{-} "10-u";"2-u"^{}
\ar@{-} "8";"8-u"^{}
\ar@{-} "6";"6-u"^{}
\ar@{-} "2";"2-u"^{}
\ar@{-} "1";"1-u"^{}
\ar@{-} "1";"2"^{}
\ar@{-} "1";"4"^{}
\ar@{-} "2";"3"^{}
\ar@{-} "3";"4"^{}
\ar@{-} "5";"6"^{}
\ar@{-} "5";"4"^{}
\ar@{-} "1";"6"^{}
\ar@{-} "7";"8"^{}
\ar@{-} "7";"5"^{}
\ar@{-} "6";"8"^{}
\ar@{-} "3";"3-1"^{}
\ar@{-} "4-1";"3-1"^{}
\ar@{-} "4-1";"4"^{}
\ar@{-} "5-1";"5"^{}
\ar@{-} "4-1";"5-1"^{}
\ar@{-} "4-1";"4-2"^{}
\ar@{-} "3-2";"3-1"^{}
\ar@{-} "3-2";"4-2"^{}
\ar@{-} "3-2";"3-3"^{}
\ar@{-} "4-3";"4-2"^{}
\ar@{-} "4-3";"3-3"^{}
\end{xy}
\end{equation*}
\end{ex}

\begin{defn}\cite{Kang}\label{def-YW2}
Let $X={\rm A}^{(2)}_{2n-2}$ or ${\rm D}^{(2)}_{n}$ and
$Y$ be a Young wall of ground state $\lambda$.
\begin{enumerate}
\item A column in $Y$ is called a {\it full column} if its height is a multiple of the unit length.
\item The wall $Y$ is said to be {\it proper} if none of two full columns in $Y$ have the same height.
\end{enumerate}
\end{defn}

\begin{defn}\cite{Kang}\label{def-YW2a}
Let $X={\rm A}^{(2)}_{2n-2}$ or ${\rm D}^{(2)}_{n}$ and
$Y$ be a proper Young wall of ground state $\lambda$.
\begin{enumerate}
\item 
Let $i\in I$ and $Y'$ be a wall obtained from $Y$ by 
removing an $i$-block. This $i$-block is said to be a {\it removable} $i$-{\it block} of $Y$
if $Y'$ is a proper Young wall. 
\item Let $i\in I$. If we obtain a proper Young wall by adding
an $i$-block to a place of $Y$ then the place is said to be an
$i$-{\it admissible} {\it slot} of $Y$.
\end{enumerate}
\end{defn}

\begin{defn}\cite{Ka}
Let $X={\rm A}^{(2)}_{2n-2}$ or ${\rm D}^{(2)}_{n}$ and
$Y$ be a proper Young wall of ground state $\lambda$. We take $t\in\{1,n\}$. 
\begin{enumerate}
\item[(i)]
Let $Y'$ be a wall
obtained by adding two $t$-blocks of shape (\ref{smpl2}) to the top of a column in $Y$:
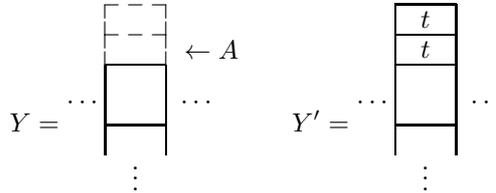
\[
Y=
\begin{xy}
(14,10) *{\leftarrow A}="A",
(-3,3) *{\cdots}="dot1",
(12,3) *{\cdots}="dot2",
(0,-4) *{}="-1",
(8,-4)*{}="-2",
(4,-6)*{\vdots}="dot3",
(0,0) *{}="1",
(8,0)*{}="2",
(8,8)*{}="3",
(0,8)*{}="4",
(8,12)*{}="5",
(0,12)*{}="6",
(8,16)*{}="7",
(0,16)*{}="8",
\ar@{-} "1";"-1"^{}
\ar@{-} "-2";"2"^{}
\ar@{-} "1";"2"^{}
\ar@{-} "1";"4"^{}
\ar@{-} "2";"3"^{}
\ar@{-} "3";"4"^{}
\ar@{--} "3";"5"^{}
\ar@{--} "4";"6"^{}
\ar@{--} "6";"5"^{}
\ar@{--} "7";"5"^{}
\ar@{--} "6";"8"^{}
\ar@{--} "7";"8"^{}
\end{xy}\qquad
Y'=
\begin{xy}
(14,14) *{}="B",
(4,10)*{t}="t1",
(4,14)*{t}="t2",
(-3,3) *{\cdots}="dot1",
(12,3) *{\cdots}="dot2",
(0,-4) *{}="-1",
(8,-4)*{}="-2",
(4,-6)*{\vdots}="dot3",
(0,0) *{}="1",
(8,0)*{}="2",
(8,8)*{}="3",
(0,8)*{}="4",
(8,12)*{}="5",
(0,12)*{}="6",
(8,16)*{}="7",
(0,16)*{}="8",
\ar@{-} "1";"-1"^{}
\ar@{-} "-2";"2"^{}
\ar@{-} "1";"2"^{}
\ar@{-} "1";"4"^{}
\ar@{-} "2";"3"^{}
\ar@{-} "3";"4"^{}
\ar@{-} "3";"5"^{}
\ar@{-} "4";"6"^{}
\ar@{-} "6";"5"^{}
\ar@{-} "7";"5"^{}
\ar@{-} "6";"8"^{}
\ar@{-} "7";"8"^{}
\end{xy}
\]
In $Y$, a slot is named $A$ as above. 
The slot $A$ in $Y$ is said to be {\it double} $t$-{\it admissible} if $Y'$ is also a proper Young wall.
\item[(ii)]
Let $Y''$ be a wall
obtained by removing two $t$-blocks of shape (\ref{smpl2}) from the top of a column in $Y$:
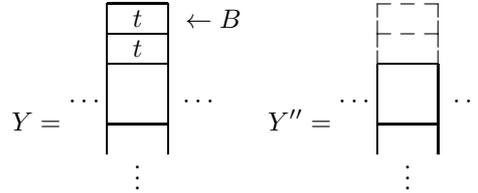
\[
Y=
\begin{xy}
(14,14) *{\leftarrow B}="B",
(4,10)*{t}="t1",
(4,14)*{t}="t2",
(-3,3) *{\cdots}="dot1",
(12,3) *{\cdots}="dot2",
(0,-4) *{}="-1",
(8,-4)*{}="-2",
(4,-6)*{\vdots}="dot3",
(0,0) *{}="1",
(8,0)*{}="2",
(8,8)*{}="3",
(0,8)*{}="4",
(8,12)*{}="5",
(0,12)*{}="6",
(8,16)*{}="7",
(0,16)*{}="8",
\ar@{-} "1";"-1"^{}
\ar@{-} "-2";"2"^{}
\ar@{-} "1";"2"^{}
\ar@{-} "1";"4"^{}
\ar@{-} "2";"3"^{}
\ar@{-} "3";"4"^{}
\ar@{-} "3";"5"^{}
\ar@{-} "4";"6"^{}
\ar@{-} "6";"5"^{}
\ar@{-} "7";"5"^{}
\ar@{-} "6";"8"^{}
\ar@{-} "7";"8"^{}
\end{xy}\quad
Y''=
\begin{xy}
(14,10) *{}="A",
(-3,3) *{\cdots}="dot1",
(12,3) *{\cdots}="dot2",
(0,-4) *{}="-1",
(8,-4)*{}="-2",
(4,-6)*{\vdots}="dot3",
(0,0) *{}="1",
(8,0)*{}="2",
(8,8)*{}="3",
(0,8)*{}="4",
(8,12)*{}="5",
(0,12)*{}="6",
(8,16)*{}="7",
(0,16)*{}="8",
\ar@{-} "1";"-1"^{}
\ar@{-} "-2";"2"^{}
\ar@{-} "1";"2"^{}
\ar@{-} "1";"4"^{}
\ar@{-} "2";"3"^{}
\ar@{-} "3";"4"^{}
\ar@{--} "3";"5"^{}
\ar@{--} "4";"6"^{}
\ar@{--} "6";"5"^{}
\ar@{--} "7";"5"^{}
\ar@{--} "6";"8"^{}
\ar@{--} "7";"8"^{}
\end{xy}
\]
A block in $Y$ is named $B$ as above.
The block $B$ in $Y$ is said to be {\it double} $t$-{\it removable}
if $Y''$ is also a proper Young wall.
\item[(iii)] Other $i$-admissible slot (resp. removable $i$-block)
is said to be single admissible (resp. single removable) for $i\in I$.
\end{enumerate}
\end{defn}

\begin{ex}\label{ex-A2-1}
We consider the following Young wall of type ${\rm A}^{(2)}_{4}$ in Example \ref{ex-YW}.
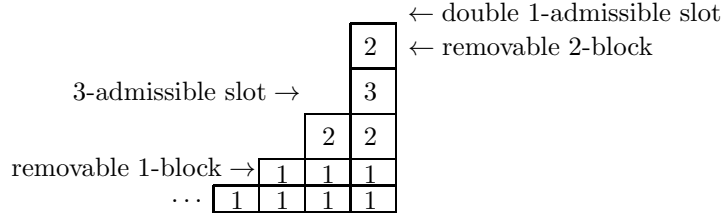
\begin{equation*}
\begin{xy}
(-15.5,-2) *{\ 1}="000-1",
(-9.5,-2) *{\ 1}="00-1",
(-3.5,-2) *{\ 1}="0-1",
(1.5,-2) *{\ \ 1}="0-1",
(-21.5,-2) *{\dots}="00000",
(-9.5,1.5) *{\ 1}="000",
(-3.5,1.5) *{\ 1}="00",
(1.5,1.5) *{\ \ 1}="0",
(1.5,6.5) *{\ \ 2}="02",
(-4,6.5) *{\ \ 2}="002",
(1.5,12.5) *{\ \ 3}="03",
(1.5,18.5) *{\ \ 2}="04",
(22.5,18.5) *{\ \ \leftarrow\text{removable }2\text{-block}}="arrow1",
(27,23) *{\ \ \leftarrow\text{double }1\text{-admissible slot}}="arrow2",
(-21.5,12.5) *{3\text{-admissible slot}\rightarrow}="arrow3",
(-28.5,2.5) *{\text{removable }1\text{-block}\rightarrow}="arrow4",
(0,0) *{}="1",
(0,-3.5) *{}="1-u",
(6,0)*{}="2",
(6,-3.5)*{}="2-u",
(6,3.5)*{}="3",
(0,3.5)*{}="4",
(-6,3.5)*{}="5",
(-6,0)*{}="6",
(-6,-3.5)*{}="6-u",
(-12,3.5)*{}="7",
(-12,0)*{}="8",
(-12,-3.6)*{}="8-u",
(-18,3.5)*{}="9",
(-18,0)*{}="10-a",
(-18,0)*{}="10",
(-18,-3.5)*{}="10-u",
(6,9.5)*{}="3-1",
(6,15.5)*{}="3-2",
(6,21.5)*{}="3-3",
(0,9.5)*{}="4-1",
(0,15.5)*{}="4-2",
(0,21.5)*{}="4-3",
(-6,9.5)*{}="5-1",
\ar@{-} "8";"10-a"^{}
\ar@{-} "10-u";"10-a"^{}
\ar@{-} "10-u";"2-u"^{}
\ar@{-} "8";"8-u"^{}
\ar@{-} "6";"6-u"^{}
\ar@{-} "2";"2-u"^{}
\ar@{-} "1";"1-u"^{}
\ar@{-} "1";"2"^{}
\ar@{-} "1";"4"^{}
\ar@{-} "2";"3"^{}
\ar@{-} "3";"4"^{}
\ar@{-} "5";"6"^{}
\ar@{-} "5";"4"^{}
\ar@{-} "1";"6"^{}
\ar@{-} "7";"8"^{}
\ar@{-} "7";"5"^{}
\ar@{-} "6";"8"^{}
\ar@{-} "3";"3-1"^{}
\ar@{-} "4-1";"3-1"^{}
\ar@{-} "4-1";"4"^{}
\ar@{-} "5-1";"5"^{}
\ar@{-} "4-1";"5-1"^{}
\ar@{-} "4-1";"4-2"^{}
\ar@{-} "3-2";"3-1"^{}
\ar@{-} "3-2";"4-2"^{}
\ar@{-} "3-2";"3-3"^{}
\ar@{-} "4-3";"4-2"^{}
\ar@{-} "4-3";"3-3"^{}
\end{xy}
\end{equation*}
Since it holds $h_1=4$, $h_2=2$, $h_3=1$ and $h_{\ell}=\frac{1}{2}$ ($\ell\geq4$), 
the first, second and third columns are full columns, but the fourth column is not since its height is half-unit.
We see that the Young wall is proper.
It has only one double $1$-admissible slot and other admissible slots and removable blocks are single.
\end{ex}

\section{Combinatorial descriptions of polyhedral realizations for $B(\lambda)$}

In this section, we present our main results.

\subsection{Adaptedness and several notation}\label{seno}

\begin{defn}\label{adapt}\cite{KaN}
Let $A=(a_{i,j})_{i,j\in I}$ be the symmetrizable generalized Cartan matrix of $\mathfrak{g}$.
We assume that a sequence $\io=(\cdots,i_3,i_2,i_1)$
satisfies $(\ref{seq-con})$.
We say $\iota$ is {\it adapted} to $A$ if the following condition holds : 
For $i,j\in I$ such that $a_{i,j}<0$, the subsequence of $\iota$ consisting of all $i$, $j$ is
\[
(\cdots,i,j,i,j,i,j,i,j)\quad {\rm or}\quad (\cdots,j,i,j,i,j,i,j,i).
\]
In case of the Cartan matrix is fixed, the sequence $\iota$ is shortly said to be {\it adapted}.
\end{defn}

\begin{ex}
In case that
$\mathfrak{g}$ is of type ${\rm A}^{(1)}_2$, $\iota=(\cdots,2,1,3,2,1,3,2,1,3)$,
it holds $a_{1,2}=a_{2,3}=a_{1,3}=-1$. 
\begin{itemize}
\item
The subsequence consisting of $1$ and $2$ is $(\cdots,2,1,2,1,2,1)$.
\item
The subsequence consisting of $2$ and $3$ is $(\cdots,2,3,2,3,2,3)$.
\item
The subsequence consisting of $1$ and $3$ is $(\cdots,1,3,1,3,1,3)$.
\end{itemize}
Thus, $\iota$ is an adapted sequence.
\end{ex}

In the rest of article, let us fix a sequence $\iota=(\cdots,i_3,i_2,i_1)$ adapted to $A$.
One defines a set of integers $(p_{i,j})_{i,j\in I;a_{i,j}<0}$ by
\begin{equation}\label{pij}
p_{i,j}=\begin{cases}
1 & {\rm if}\ {\rm the\ subsequence\ of\ }\iota{\rm\ consisting\ of}\ i,j\ {\rm is}\ (\cdots,j,i,j,i,j,i),\\
0 & {\rm if}\ {\rm the\ subsequence\ of\ }\iota{\rm\ consisting\ of}\ i,j\ {\rm is}\ (\cdots,i,j,i,j,i,j).
\end{cases}
\end{equation}
Note that if $a_{i,j}<0$ then we have
\begin{equation}\label{pij2}
p_{i,j}+p_{j,i}=1.
\end{equation}
We identify each single index $r\in\mathbb{Z}_{\geq1}$ with a pair $(s,k)\in \mathbb{Z}_{\geq1}\times I$
if $i_r=k$ holds and $k$ appears $s$ times in $i_r$, $i_{r-1}$, $\cdots,i_1$.
For example, when $\iota=(\cdots,2,1,3,2,1,3,2,1,3)$, single indices
$\cdots,6,5,4,3,2,1$ are identified with 
\[
\cdots,(2,2),(2,1),(2,3),(1,2),(1,1),(1,3).
\]
The notation $x_{r}$, $\beta_{r}$, $\beta^{\pm}_{r}$ and $\what{S}'_{r}$ in the subsections \ref{poly-uqm}, \ref{poly-uql} are also rewritten as
\[
x_r=x_{s,k},\quad \beta_r=\beta_{s,k},\quad \beta^{\pm}_r=\beta^{\pm}_{s,k},\quad \what{S}'_r=\what{S}'_{s,k}.
\]
If $(s,k)\notin \mathbb{Z}_{\geq1}\times I$ then we understand $x_{s,k}=0$.
By this identification and the ordinary order on $\mathbb{Z}_{\geq1}$, that is, $1<2<3<4<5<6<\cdots$,
an order on $\mathbb{Z}_{\geq1}\times I$ is defined. In case of $\iota=(\cdots,2,1,3,2,1,3,2,1,3)$, the order is
\[
\cdots>(2,2)>(2,1)>(2,3)>(1,2)>(1,1)>(1,3).
\]
Using the notation in (\ref{pij}), by the adaptedness of $\iota$, one can write $\beta_{s,k}$ as
\begin{equation}\label{pij3}
\beta_{s,k}=x_{s,k}+x_{s+1,k}+\sum_{j\in I; a_{k,j}<0} a_{k,j}x_{s+p_{j,k},j}.
\end{equation}
The above
$\beta_{s,k}$ can be regarded as an analog of
simple root $\alpha_{k}$ of $^L\mathfrak{g}$
since $\alpha_{k}$ is expressed by 
\[
\alpha_{k}=\Lambda_{k}+\Lambda_{k}+\sum_{j\in I; a_{k,j}<0} a_{k,j}\Lambda_{j}
\]
on $\bigoplus_{i\in I}\mathbb{Z}h_i$.
This is a reason why the
inequalities defining ${\rm Im}(\Psi^{(\lambda)}_{\iota})$
is expressed by combinatorial objects related to $^L\mathfrak{g}$ rather than $\mathfrak{g}$
in Theorem \ref{thmA1}, \ref{thmA2}.

\subsection{Non-negative integers $P^X_k$}

Let $X={\rm A}^{(1)},\ {\rm C}^{(1)},\ {\rm A}^{(2)}$ or ${\rm D}^{(2)}$ and
\[
\begin{cases}
k\in I & \text{if }X={\rm A}^{(1)},\ \text{or }{\rm C}^{(1)},\\
k\in I\setminus\{1\} & \text{if }X={\rm A}^{(2)},\\
k\in I\setminus\{1,n\} & \text{if }X={\rm D}^{(2)}.
\end{cases}
\]
For $t\in\mathbb{Z}$, let $P^X_k(t)\in\mathbb{Z}_{\geq0}$ denote the non-negative integer defined as follows:
Setting $P^X_k(k):=0$, one recursively defines as
\begin{equation}\label{pk1}
P^X_k(t):=P^X_k(t-1)+p_{\pi_X(t),\pi_X(t-1)}\ (\text{for } t>k),
\end{equation}
\begin{equation}\label{pk2}
P^X_k(t):=P^X_k(t+1)+p_{\pi_X(t),\pi_X(t+1)}\ (\text{for } t<k),
\end{equation}
where the notation $\pi_X$ is defined in Definition \ref{pi1-def} and we understand $p_{j,j}:=0$ for all $j\in I$.

Next,
considering the map
$\{1,2,\cdots,2n-2\}\rightarrow \{1,2,\cdots,n\}$
defined as
\[
\ell\mapsto \ell,\ 2n-\ell\mapsto \ell \qquad (2\leq \ell\leq n-1),
\]
\[
1\mapsto1,\ n\mapsto n,
\]
we extend it to a map 
\begin{equation}\label{piprime}
\pi':\mathbb{Z}_{\geq1}\rightarrow\{1,2,\cdots,n\}
\end{equation}
by periodicity $2n-2$.
For 
$X={\rm A}^{(2)}$ or ${\rm D}^{(2)}$
and
\[
k=
\begin{cases}
1 & \text{if }X={\rm A}^{(2)},\\
1 \text{ or }n & \text{if }X={\rm D}^{(2)},
\end{cases}
\]
we recursively define integers $P^X_k(\ell)$ ($\ell\in\mathbb{Z}_{\geq k}$) as
\[
P^X_k(k):=0,\ \ P^X_k(\ell)=P^X_k(\ell-1)+p_{\pi'(\ell),\pi'(\ell-1)}\ (\ell>k).
\]

\subsection{Type ${\rm A}_{n-1}^{(1)}$ and ${\rm D}^{(2)}_n$-case}

For any pair $(i,j)\in\mathbb{Z}\times \mathbb{Z}$, $s\in\mathbb{Z}$, $k\in I$ and $X={\rm A}^{(1)}$ or ${\rm C}^{(1)}$, we set
\begin{equation}\label{ovlL1}
L^X_{s,k,\iota}(i,j):=x_{s+P^X_k(i+j)+{\rm min}\{k-j,i\},\pi_X(i+j)}\in (\mathbb{Q}^{\infty})^*.
\end{equation}
For an extended Young diagram $T$ with $y_{\infty}=k$, we define
\begin{equation}\label{ovlL2}
L^X_{s,k,\iota}(T)
:=\sum_{P:{\rm concave\ corner\ of\ }T} L^X_{s,k,\iota}(P) - \sum_{P:{\rm convex\ corner\ of\ }T} L^X_{s,k,\iota}(P)
\in (\mathbb{Q}^{\infty})^*.
\end{equation}
For $k\in I$,
let ${\rm EYD}_{k}$ denote the set of extended Young diagrams with $y_{\infty}=k$.
One defines
\begin{equation}\label{r-tilde-A0}
\fbox{$r$}^{X}_k:=x_{P^X_k(r),\pi_X(r)}
-x_{1+P^X_k(r-1),\pi_X(r-1)} \qquad (r\in\mathbb{Z}_{\geq k+1}),
\end{equation}
\begin{equation}\label{r-tilde-A1}
\fbox{$\tilde{r}$}_k^{X}:=x_{P^X_k(r-1),\pi_X(r-1)}
-x_{1+P^X_k(r),\pi_X(r)} \qquad (r\in\mathbb{Z}_{\leq k}).
\end{equation}
We also set
\begin{equation}\label{comb-lam-1}
{\rm Comb}_{k,\iota}^{X}[\lambda]:=
\begin{cases}
\{-x_{1,k}+\langle h_k,\lambda \rangle\} & \text{if }(1,k)<(1,\pi_X(k+1)),\ (1,k)<(1,\pi_X(k-1)),\\
\{\fbox{$\tilde{r}$}_k^{X}+\langle h_k,\lambda \rangle\ | r\in\mathbb{Z}_{\leq k}\} & \text{if }(1,k)<(1,\pi_X(k+1)),\ (1,k)>(1,\pi_X(k-1)),\\
\{\fbox{$r$}^{X}_k+\langle h_k,\lambda \rangle\ | r\in\mathbb{Z}_{\geq k+1}\} & \text{if }(1,k)>(1,\pi_X(k+1)),\ (1,k)<(1,\pi_X(k-1)),\\
\{
L^X_{0,k,\iota}(T)+\langle h_k,\lambda \rangle | T\in{\rm EYD}_{k}\setminus\{\phi^k\}
\} & \text{if }(1,k)>(1,\pi_X(k+1)),\ (1,k)>(1,\pi_X(k-1))\\ 
\end{cases}
\end{equation}
and
\begin{equation}\label{comb-inf-1}
{\rm Comb}_{\iota}^{X}[\infty]:=
\{ L^X_{s,k,\iota}(T) | s\in\mathbb{Z}_{\geq1},\ k\in I,\ T\in {\rm EYD}_k \}.
\end{equation}
Here, $\phi^k\in{\rm EYD}_k$ is the extended Young diagram without boxes, that is, $\phi^k=(y_{\ell})_{\ell\in\mathbb{Z}_{\geq0}}$, 
$y_{\ell}=k$ for any $\ell\in\mathbb{Z}_{\geq0}$.

\begin{thm}\label{thmA1}
We suppose that $\mathfrak{g}$ is of type $X={\rm A}_{n-1}^{(1)}$ $(n\geq2)$ or $X={\rm D}_{n}^{(2)}$ $(n\geq3)$ and $\iota$ is adapted.
Then for any $\lambda\in P^+$, the pair $(\iota,\lambda)$ is $\Xi'$-ample and
\[
{\rm Im}(\Psi_{\iota}^{(\lambda)})
=\left\{ \mathbf{a}\in\mathbb{Z}^{\infty}
\left|
\varphi(\mathbf{a})\geq0\\
\text{ for any }
\varphi\in
{\rm Comb}^{X^L}_{\iota}[\infty]\cup
\bigcup_{k\in I}
{\rm Comb}^{X^L}_{k,\iota}[\lambda]
\right.
\right\}.
\]
Here, $({\rm A}^{(1)}_{n-1})^L={\rm A}^{(1)}_{n-1}$ and $({\rm D}^{(2)}_{n})^L={\rm C}^{(1)}_{n-1}$.
\end{thm}

\begin{rem}\label{rem-1}
Let us put
\[
[r]^X_k:= 
\frac{X_{P^X_k(r),\pi_X(r)}}
{X_{1+P^X_k(r-1),\pi_X(r-1)}},\qquad
[\tilde{r}]^X_k:=
\frac{X_{P^X_k(r-1),\pi_X(r-1)}}
{X_{1+P^X_k(r),\pi_X(r)}},
\]
where tropicalizing these Laurent monomials, we get linear functions
in (\ref{r-tilde-A0}), (\ref{r-tilde-A1}).
Then
the set
$\{[r]^X_k | r\in\mathbb{Z}\}$ (resp. $\{[\tilde{r}]^X_k | r\in\mathbb{Z}\}$) coincides with
the connected component of $\mathcal{M}_c$ in \cite{K}
with $c_{m,j}=p_{m,j}$ for $m,j\in I$ such that $a_{m,j}<0$ 
containing
$\frac{X_{P^X_k(k+1),\pi_X(k+1)}}
{X_{1,k}}$
(resp. $\frac{X_{P^X_k(k-1),\pi_X(k-1)}}
{X_{1,k}}$),
which is isomorphic to the crystal base $B(\Lambda_{\pi_X(k+1)}-\Lambda_{k})$ (resp. $B(\Lambda_{\pi_X(k-1)}-\Lambda_{k})$)
of the extremal weight module $V(\Lambda_{\pi_X(k+1)}-\Lambda_{k})$ (resp. $V(\Lambda_{\pi_X(k-1)}-\Lambda_{k})$) of $U_q(\mathfrak{g})$.
In this way, functions in ${\rm Comb}^{X^L}_{\iota}[\infty]$ and
${\rm Comb}^{X^L}_{k,\iota}[\lambda]$ are deeply related to representation
theory of $U_q(^L\mathfrak{g})$.
In Theorem 3.2 of \cite{HN2}, monomial realizations of level 0 extremal
fundamental weight crystals are given. 
\end{rem}

\begin{ex}\label{ex-2}
Let us consider the case
$\mathfrak{g}$ is of type ${\rm A}_2^{(1)}$ and $\iota=(\cdots,3,1,2,3,1,2)$, which is the same setting
as in Example 4.4 of \cite{Ka}.
In the paper \cite{Ka},
we computed several inequalities defining
${\rm Im}(\Psi_{\iota})=\{\mathbf{a}\in \mathbb{Z}^{\infty} |
 \varphi(\mathbf{a})\geq0\ \text{ for any }\varphi\in {\rm Comb}^{{\rm A}^{(1)}}_{\iota}[\infty]\}$:
\[
{\rm Im}(\Psi_{\iota})
=\left\{ 
\textbf{a}\in\mathbb{Z}^{\infty} \left| 
\begin{array}{l}
s\in\mathbb{Z}_{\geq1},\ a_{s,1}\geq0,
a_{s+1,2}+a_{s,3}-a_{s+1,1}\geq0,
a_{s+1,3}+a_{s,3}-a_{s+2,2}\geq0,\\
2a_{s+1,2}-a_{s+1,3}\geq0,
a_{s+1,2}+a_{s+1,1}-a_{s+2,2}\geq0,
a_{s+1,2}+a_{s+1,3}-a_{s+2,1}\geq0,\cdots\\
a_{s,2}\geq0,\ 
a_{s,1}+a_{s,3}-a_{s+1,2}\geq0,
a_{s,1}+a_{s+1,1}-a_{s+1,3}\geq0,\\
2a_{s,3}-a_{s+1,1}\geq0,
a_{s,3}+a_{s+1,2}-a_{s+1,3}\geq0,
a_{s,3}+a_{s+1,1}-a_{s+2,2}\geq0,\cdots\\
a_{s,3}\geq0,
a_{s+1,1}+a_{s+1,2}-a_{s+1,3}\geq0,
a_{s+2,2}+a_{s+1,2}-a_{s+2,1}\geq0,\\
2a_{s+1,1}-a_{s+2,2}\geq0,
a_{s+1,1}+a_{s+1,3}-a_{s+2,1}\geq0,
a_{s+1,1}+a_{s+2,2}-a_{s+2,3}\geq0,\cdots
\end{array}
\right.
\right\}.
\]
Let us compute several inequalities coming from ${\rm Comb}_{k,\iota}^{{\rm A}^{(1)}}[\lambda]$ for $k=1,2,3$.
One obtains
\[
\cdots,\ P_1^{{\rm A}^{(1)}}(-1)=1,\ P_1^{{\rm A}^{(1)}}(0)=0,\ P_1^{{\rm A}^{(1)}}(1)=0,\ P_1^{{\rm A}^{(1)}}(2)=1,\ P_1^{{\rm A}^{(1)}}(3)=1,\cdots,
\]
\[
\cdots,\ P_2^{{\rm A}^{(1)}}(-1)=1,\ P_2^{{\rm A}^{(1)}}(0)=0,\ P_2^{{\rm A}^{(1)}}(1)=0,\ P_2^{{\rm A}^{(1)}}(2)=0,\ P_2^{{\rm A}^{(1)}}(3)=0,\ P_2^{{\rm A}^{(1)}}(4)=1,\cdots,
\]
\[
\cdots,\ P_3^{{\rm A}^{(1)}}(1)=1,\ P_3^{{\rm A}^{(1)}}(2)=1,\ P_3^{{\rm A}^{(1)}}(3)=0,\ P_3^{{\rm A}^{(1)}}(4)=1,\ P_3^{{\rm A}^{(1)}}(5)=2,\cdots .
\]
By the definition (\ref{comb-lam-1}), we have ${\rm Comb}_{2,\iota}^{{\rm A}^{(1)}}[\lambda]=\{-x_{1,2}+\langle h_2,\lambda \rangle\}$
and
\[
{\rm Comb}_{1,\iota}^{{\rm A}^{(1)}}[\lambda]=
\{\fbox{$r$}^{{\rm A}^{(1)}}_1+\langle h_1,\lambda \rangle\ | r\in\mathbb{Z}_{\geq 2}\}.
\]
One can compute
\[
\fbox{$2$}^{{\rm A}^{(1)}}_1=x_{1,2}-x_{1,1}, \ 
\fbox{$3$}^{{\rm A}^{(1)}}_1=x_{1,3}-x_{2,2},\ 
\fbox{$4$}^{{\rm A}^{(1)}}_1=x_{2,1}-x_{2,3},\ 
\]
\[
\fbox{$5$}^{{\rm A}^{(1)}}_1=x_{3,2}-x_{3,1}, \ 
\fbox{$6$}^{{\rm A}^{(1)}}_1=x_{3,3}-x_{4,2}.
\]
Next, let us turn to ${\rm Comb}_{3,\iota}^{{\rm A}^{(1)}}[\lambda]$.
The following diagrams are elements in ${\rm EYD}_{3}\setminus\{\phi^3\}$:
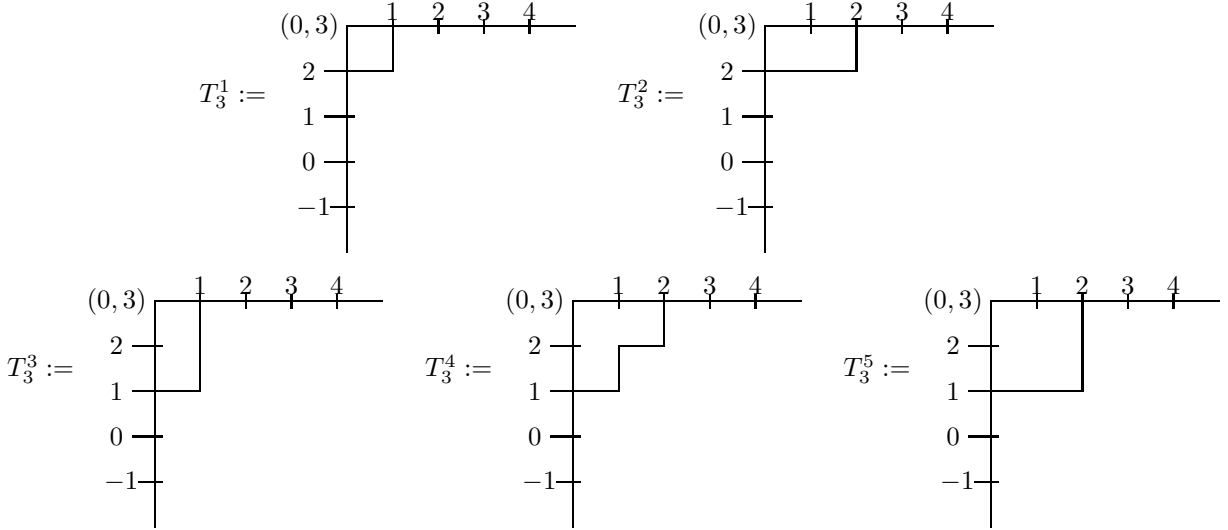
\[
\begin{xy}
(40,-9) *{T_3^1:=}="Y1",
(55,0) *{}="1a",
(85,0)*{}="2a",
(55,-30)*{}="3a",
(50,0)*{(0,3)}="4a",
(61,2) *{1}="10a",
(61,-1) *{}="1010a",
(61,-6) *{}="1-0a",
(67,2) *{2}="11a",
(67,-1) *{}="1111a",
(73,2) *{3}="12a",
(73,-1) *{}="1212a",
(79,2) *{4}="13a",
(79,-1) *{}="1313a",
(50,-6)*{\ 2\ }="5a",
(56,-6)*{}="55a",
(50,-12)*{\ 1\ }="6a",
(56,-12)*{}="66a",
(50,-18)*{\ 0\ }="7a",
(56,-18)*{}="77a",
(50,-24)*{\ -1}="8a",
(56,-24)*{}="88a",
(95,-9) *{T_3^2:=}="Y2",
(110,0) *{}="1b",
(140,0)*{}="2b",
(110,-30)*{}="3b",
(105,0)*{(0,3)}="4b",
(116,2) *{1}="10b",
(116,-1) *{}="1010b",
(122,-6) *{}="1-0b",
(122,2) *{2}="11b",
(122,-1) *{}="1111b",
(128,2) *{3}="12b",
(128,-1) *{}="1212b",
(134,2) *{4}="13b",
(134,-1) *{}="1313b",
(105,-6)*{\ 2\ }="5b",
(111,-6)*{}="55b",
(105,-12)*{\ 1\ }="6b",
(111,-12)*{}="66b",
(105,-18)*{\ 0\ }="7b",
(111,-18)*{}="77b",
(105,-24)*{\ -1}="8b",
(111,-24)*{}="88b",
\ar@{-} "1a";"2a"^{}
\ar@{-} "1a";"3a"^{}
\ar@{-} "5a";"55a"^{}
\ar@{-} "6a";"66a"^{}
\ar@{-} "7a";"77a"^{}
\ar@{-} "8a";"88a"^{}
\ar@{-} "10a";"1010a"^{}
\ar@{-} "11a";"1111a"^{}
\ar@{-} "12a";"1212a"^{}
\ar@{-} "13a";"1313a"^{}
\ar@{-} "55a";"1-0a"^{}
\ar@{-} "1010a";"1-0a"^{}
\ar@{-} "1b";"2b"^{}
\ar@{-} "1b";"3b"^{}
\ar@{-} "5b";"55b"^{}
\ar@{-} "6b";"66b"^{}
\ar@{-} "7b";"77b"^{}
\ar@{-} "8b";"88b"^{}
\ar@{-} "10b";"1010b"^{}
\ar@{-} "11b";"1111b"^{}
\ar@{-} "12b";"1212b"^{}
\ar@{-} "13b";"1313b"^{}
\ar@{-} "55b";"1-0b"^{}
\ar@{-} "11b";"1-0b"^{}
\end{xy}
\]
\[
\begin{xy}
(-15,-9) *{T_3^3:=}="Y3",
(0,0) *{}="1",
(30,0)*{}="2",
(0,-30)*{}="3",
(-5,0)*{(0,3)}="4",
(6,2) *{1}="10",
(6,-1) *{}="1010",
(12,2) *{2}="11",
(12,-1) *{}="1111",
(18,2) *{3}="12",
(18,-1) *{}="1212",
(24,2) *{4}="13",
(24,-1) *{}="1313",
(-5,-6)*{\ 2\ }="5",
(1,-6)*{}="55",
(6,-12)*{}="6-12",
(-5,-12)*{\ 1\ }="6",
(1,-12)*{}="66",
(-5,-18)*{\ 0\ }="7",
(1,-18)*{}="77",
(-5,-24)*{\ -1}="8",
(1,-24)*{}="88",
(40,-9) *{T_3^4:=}="Y4",
(55,0) *{}="1a",
(85,0)*{}="2a",
(55,-30)*{}="3a",
(50,0)*{(0,3)}="4a",
(61,2) *{1}="10a",
(61,-1) *{}="1010a",
(61,-6) *{}="6-6a",
(61,-12)*{}="6-12a",
(67,-6)*{}="12-6a",
(67,2) *{2}="11a",
(67,-1) *{}="1111a",
(73,2) *{3}="12a",
(73,-1) *{}="1212a",
(79,2) *{4}="13a",
(79,-1) *{}="1313a",
(50,-6)*{\ 2\ }="5a",
(56,-6)*{}="55a",
(50,-12)*{\ 1\ }="6a",
(56,-12)*{}="66a",
(50,-18)*{\ 0\ }="7a",
(56,-18)*{}="77a",
(50,-24)*{\ -1}="8a",
(56,-24)*{}="88a",
(95,-9) *{T_3^5:=}="Y5",
(110,0) *{}="1b",
(140,0)*{}="2b",
(110,-30)*{}="3b",
(105,0)*{(0,3)}="4b",
(116,2) *{1}="10b",
(116,-1) *{}="1010b",
(122,-12) *{}="1-0b",
(122,2) *{2}="11b",
(122,-1) *{}="1111b",
(128,2) *{3}="12b",
(128,-1) *{}="1212b",
(134,2) *{4}="13b",
(134,-1) *{}="1313b",
(105,-6)*{\ 2\ }="5b",
(111,-6)*{}="55b",
(105,-12)*{\ 1\ }="6b",
(111,-12)*{}="66b",
(105,-18)*{\ 0\ }="7b",
(111,-18)*{}="77b",
(105,-24)*{\ -1}="8b",
(111,-24)*{}="88b",
\ar@{-} "1";"2"^{}
\ar@{-} "1";"3"^{}
\ar@{-} "5";"55"^{}
\ar@{-} "6";"66"^{}
\ar@{-} "6";"6-12"^{}
\ar@{-} "1010";"6-12"^{}
\ar@{-} "7";"77"^{}
\ar@{-} "8";"88"^{}
\ar@{-} "10";"1010"^{}
\ar@{-} "11";"1111"^{}
\ar@{-} "12";"1212"^{}
\ar@{-} "13";"1313"^{}
\ar@{-} "1a";"2a"^{}
\ar@{-} "1a";"3a"^{}
\ar@{-} "5a";"55a"^{}
\ar@{-} "6a";"66a"^{}
\ar@{-} "7a";"77a"^{}
\ar@{-} "8a";"88a"^{}
\ar@{-} "10a";"1010a"^{}
\ar@{-} "11a";"1111a"^{}
\ar@{-} "12a";"1212a"^{}
\ar@{-} "13a";"1313a"^{}
\ar@{-} "6-12a";"66a"^{}
\ar@{-} "6-12a";"6-6a"^{}
\ar@{-} "12-6a";"6-6a"^{}
\ar@{-} "12-6a";"1111a"^{}
\ar@{-} "1b";"2b"^{}
\ar@{-} "1b";"3b"^{}
\ar@{-} "5b";"55b"^{}
\ar@{-} "6b";"66b"^{}
\ar@{-} "7b";"77b"^{}
\ar@{-} "8b";"88b"^{}
\ar@{-} "10b";"1010b"^{}
\ar@{-} "11b";"1111b"^{}
\ar@{-} "12b";"1212b"^{}
\ar@{-} "13b";"1313b"^{}
\ar@{-} "6b";"1-0b"^{}
\ar@{-} "11b";"1-0b"^{}
\end{xy}
\]
The diagram $T_3^1$ has concave corners $(1,3)$ and $(0,2)$ and a convex corner $(1,2)$. Thus,
\[
L^{{\rm A}^{(1)}}_{0,3,\iota}(T^1_3)=x_{P_3^{{\rm A}^{(1)}}(4),\pi_{{\rm A}^{(1)}}(4)}+x_{P_3^{{\rm A}^{(1)}}(2),\pi_{{\rm A}^{(1)}}(2)}
-x_{P_3^{{\rm A}^{(1)}}(2),\pi_{{\rm A}^{(1)}}(3)}=x_{1,1}+x_{1,2}-x_{1,3}.
\]
Similarly, it follows
\[
L^{{\rm A}^{(1)}}_{0,3,\iota}(T^2_3)=x_{P_3^{{\rm A}^{(1)}}(5),\pi_{{\rm A}^{(1)}}(5)}+x_{P_3^{{\rm A}^{(1)}}(2),\pi_{{\rm A}^{(1)}}(2)}-x_{P_3^{{\rm A}^{(1)}}(4)+1,\pi_{{\rm A}^{(1)}}(4)}
=x_{2,2}+x_{1,2}-x_{2,1},
\]
\[
L^{{\rm A}^{(1)}}_{0,3,\iota}(T^3_3)=x_{P_3^{{\rm A}^{(1)}}(4),\pi_{{\rm A}^{(1)}}(4)}+x_{P_3^{{\rm A}^{(1)}}(1),\pi_{{\rm A}^{(1)}}(1)}-x_{P_3^{{\rm A}^{(1)}}(2)+1,\pi_{{\rm A}^{(1)}}(2)}
=2x_{1,1}-x_{2,2},
\]
\begin{eqnarray*}
L^{{\rm A}^{(1)}}_{0,3,\iota}(T^4_3)&=&x_{P_3^{{\rm A}^{(1)}}(1),\pi_{{\rm A}^{(1)}}(1)}+
x_{P_3^{{\rm A}^{(1)}}(3)+1,\pi_{{\rm A}^{(1)}}(3)}
+x_{P_3^{{\rm A}^{(1)}}(5),\pi_{{\rm A}^{(1)}}(5)}\\
& &-x_{P_3^{{\rm A}^{(1)}}(2)+1,\pi_{{\rm A}^{(1)}}(2)}-x_{P_3^{{\rm A}^{(1)}}(4)+1,\pi_{{\rm A}^{(1)}}(4)}
=
x_{1,1}+x_{1,3}-x_{2,1},
\end{eqnarray*}
\[
L^{{\rm A}^{(1)}}_{0,3,\iota}(T^5_3)=x_{P_3^{{\rm A}^{(1)}}(1),\pi_{{\rm A}^{(1)}}(1)}+x_{P_3^{{\rm A}^{(1)}}(5),\pi_{{\rm A}^{(1)}}(5)}
-x_{P_3^{{\rm A}^{(1)}}(3)+2,\pi_{{\rm A}^{(1)}}(3)}
=x_{1,1}+x_{2,2}-x_{2,3}.
\]
In this way, we obtain several inequalities defining ${\rm Im}(\Psi_{\iota}^{(\lambda)})$:
\[
{\rm Im}(\Psi_{\iota}^{(\lambda)})
=\left\{ 
\textbf{a}\in\mathbb{Z}^{\infty} \left| 
\begin{array}{l}
s\in\mathbb{Z}_{\geq1},\ a_{s,1}\geq0,
a_{s+1,2}+a_{s,3}-a_{s+1,1}\geq0,
a_{s+1,3}+a_{s,3}-a_{s+2,2}\geq0,\\
2a_{s+1,2}-a_{s+1,3}\geq0,
a_{s+1,2}+a_{s+1,1}-a_{s+2,2}\geq0,
a_{s+1,2}+a_{s+1,3}-a_{s+2,1}\geq0,\cdots\\
a_{s,2}\geq0,\ 
a_{s,1}+a_{s,3}-a_{s+1,2}\geq0,
a_{s,1}+a_{s+1,1}-a_{s+1,3}\geq0,\\
2a_{s,3}-a_{s+1,1}\geq0,
a_{s,3}+a_{s+1,2}-a_{s+1,3}\geq0,
a_{s,3}+a_{s+1,1}-a_{s+2,2}\geq0,\cdots\\
a_{s,3}\geq0,
a_{s+1,1}+a_{s+1,2}-a_{s+1,3}\geq0,
a_{s+2,2}+a_{s+1,2}-a_{s+2,1}\geq0,\\
2a_{s+1,1}-a_{s+2,2}\geq0,
a_{s+1,1}+a_{s+1,3}-a_{s+2,1}\geq0,
a_{s+1,1}+a_{s+2,2}-a_{s+2,3}\geq0,\cdots\\
a_{1,2}\leq \langle h_2,\lambda\rangle,
a_{1,2}-a_{1,1}+\langle h_1,\lambda\rangle\geq0,
a_{1,3}-a_{2,2}+\langle h_1,\lambda\rangle\geq0,\\
a_{2,1}-a_{2,3}+\langle h_1,\lambda\rangle\geq0,
a_{3,2}-a_{3,1}+\langle h_1,\lambda\rangle\geq0,
a_{3,3}-a_{4,2}+\langle h_1,\lambda\rangle\geq0,\cdots\\
a_{1,1}+a_{1,2}-a_{1,3}+\langle h_3,\lambda\rangle\geq0,
a_{2,2}+a_{1,2}-a_{2,1}+\langle h_3,\lambda\rangle\geq0,\\
2a_{1,1}-a_{2,2}+\langle h_3,\lambda\rangle\geq0,
a_{1,1}+a_{1,3}-a_{2,1}+\langle h_3,\lambda\rangle\geq0,\\
a_{1,1}+a_{2,2}-a_{2,3}+\langle h_3,\lambda\rangle\geq0, \cdots
\end{array}
\right.
\right\}.
\]
\end{ex}

\begin{ex}
Let us consider another adapted sequence.
We take
$\mathfrak{g}$ as of type ${\rm A}_2^{(1)}$ just as in the previous example
and the sequence $\iota$ as $\iota=(\cdots,3,2,1,3,2,1)$. 
Then one can
compute inequalities defining ${\rm Im}(\Psi_{\iota}^{(\lambda)})$
by a similar way to Example \ref{ex-2}:
\[
{\rm Im}(\Psi_{\iota}^{(\lambda)})
=\left\{ 
\textbf{a}\in\mathbb{Z}^{\infty} \left| 
\begin{array}{l}
s\in\mathbb{Z}_{\geq1},\ 
a_{s,1}\geq0,
a_{s,2}+a_{s,3}-a_{s+1,1}\geq0,
2a_{s,3}-a_{s+1,2}\geq0,\\
a_{s,2}+a_{s+1,2}-a_{s+1,3}\geq0,
a_{s,3}+a_{s+1,1}-a_{s+1,3}\geq0,
a_{s+1,2}+a_{s,3}-a_{s+2,1}\geq0,\cdots\\
a_{s,2}\geq0,
a_{s+1,1}+a_{s,3}-a_{s+1,2}\geq0,
2a_{s+1,1}-a_{s+1,3}\geq0,\\
a_{s+1,3}+a_{s,3}-a_{s+2,1}\geq0,
a_{s+1,1}+a_{s+1,2}-a_{s+2,1}
\geq0,
a_{s+1,3}+a_{s+1,1}-a_{s+2,2}\geq0,\cdots\\
a_{s,3}\geq0,
a_{s+1,1}+a_{s+1,2}-a_{s+1,3}\geq0,
2a_{s+1,2}-a_{s+2,1}\geq0,\\
a_{s+2,1}+a_{s+1,1}-a_{s+2,2}\geq0,
a_{s+1,2}+a_{s+1,3}-a_{s+2,2}
\geq0,
a_{s+2,1}+a_{s+1,2}-a_{s+2,3}\geq0,\cdots\\
a_{1,1}\leq \langle h_1,\lambda\rangle,
a_{1,1}-a_{1,2}+\langle h_2,\lambda\rangle\geq0,
a_{1,3}-a_{2,1}+\langle h_2,\lambda\rangle\geq0,\\
a_{2,2}-a_{2,3}+\langle h_2,\lambda\rangle\geq0,
a_{3,1}-a_{3,2}+\langle h_2,\lambda\rangle\geq0,
a_{3,3}-a_{4,1}+\langle h_2,\lambda\rangle\geq0,\cdots\\
a_{1,2}+a_{1,1}-a_{1,3}+\langle h_3,\lambda\rangle\geq0,
2a_{1,2}-a_{2,1}+\langle h_3,\lambda\rangle\geq0,\\
a_{1,1}+a_{2,1}-a_{2,2}+\langle h_3,\lambda\rangle\geq0,
a_{1,3}+a_{1,2}-a_{2,2}+\langle h_3,\lambda\rangle\geq0,\\
a_{2,1}+a_{1,2}-a_{2,3}+\langle h_3,\lambda\rangle\geq0, \cdots
\end{array}
\right.
\right\}.
\]
In this case, ${\rm Comb}_{2,\iota}^{{\rm A}^{(1)}}[\lambda]$ is a set of linear functions corresponding to boxes and
${\rm Comb}_{1,\iota}^{{\rm A}^{(1)}}[\lambda]$ has a single element. 
Just as in Example \ref{ex-2}, elements in
${\rm Comb}_{3,\iota}^{{\rm A}^{(1)}}[\lambda]$ and ${\rm Comb}^{{\rm A}^{(1)}}_{\iota}[\infty]$ are described in terms of extended Young diagrams. However, 
the assigned linear functions are different from Example \ref{ex-2} since the sequence $\iota$ is changed and so is the values of $P_k^{{\rm A}^{(1)}}(t)$.
\end{ex}

\subsection{Type ${\rm A}_{2n-2}^{(2)}$ and ${\rm C}_{n-1}^{(1)}$-case}

In this subsection, let
$X={\rm A}^{(2)}$ or $X={\rm D}^{(2)}$.
We set
\[
I_X:=
\begin{cases}
I\setminus \{1\} & \text{if }X={\rm A}_{2n-2}^{(2)},\\
I\setminus\{1,n\} & \text{if }X={\rm D}_{n}^{(2)}.
\end{cases}
\]
First, we take $k$ as $k\in I_X$.
For $(i,j)\in \mathbb{Z}\times\mathbb{Z}$, $s\in\mathbb{Z}$, 
one
defines
\[
L^X_{s,k,{\rm ad}}(i,j)=x_{s+P^X_k(i+k)+[i]_-+k-j,\pi_X(i+k)},\quad
L^X_{s,k,{\rm re}}(i,j)=x_{s+P^X_k(i+k-1)+[i-1]_-+k-j,\pi_X(i+k-1)},
\]
where $[i]_-={\rm min}\{i,0\}$.
For each $T\in{\rm REYD}^{X}_{k}$, we set
\begin{eqnarray}
L^X_{s,k,\iota}(T)&:=&
\sum_{t\in I} \left(
\sum_{P:\text{single }\text{-admissible point of }T} L^X_{s,k,{\rm ad}}(P)
-\sum_{P:\text{single }\text{-removable point of }T} L^X_{s,k,{\rm re}}(P)\right)\label{L1kdef}\\
& & + 
\sum_{P:\text{double }\text{-admissible point of }T} 2L^X_{s,k,{\rm ad}}(P)
-\sum_{P:\text{double }\text{-removable point of }T} 2L^X_{s,k,{\rm re}}(P)\in (\mathbb{Q}^{\infty})^*.\nonumber
\end{eqnarray}
One defines for $r\in\mathbb{Z}_{\geq k+1}$,
\begin{equation}\label{r-tilde-D0}
\fbox{$r$}^{X}_k:=
\begin{cases}
2x_{P^X_k(r),\pi_X(r)}
-x_{1+P^X_k(r-1),\pi_X(r-1)}, & \text{if }\pi_X(r)\in I\setminus I_X\text{ and }\pi_X(r)\neq \pi_X(r-1),\\
x_{P^X_k(r),\pi_X(r)}
-2x_{1+P^X_k(r-1),\pi_X(r-1)}, & \text{if }\pi_X(r-1)\in I\setminus I_X\text{ and }\pi_X(r)\neq \pi_X(r-1),\\
x_{P^X_k(r),\pi_X(r)}
-x_{1+P^X_k(r-1),\pi_X(r-1)}, & \text{otherwise},
\end{cases}
\end{equation}
for $r\in\mathbb{Z}_{\leq k}$,
\begin{equation}\label{r-tilde-D1}
\fbox{$\tilde{r}$}^{X}_k:=
\begin{cases}
2x_{P^X_k(r-1),\pi_X(r-1)}
-x_{1+P^X_k(r),\pi_X(r)} & \text{if }\pi_X(r-1)\in I\setminus I_X\text{ and }\pi_X(r)\neq \pi_X(r-1),\\
x_{P^X_k(r-1),\pi_X(r-1)}
-2x_{1+P^X_k(r),\pi_X(r)} & \text{if }\pi_X(r)\in I\setminus I_X\text{ and }\pi_X(r)\neq \pi_X(r-1),\\
x_{P^X_k(r-1),\pi_X(r-1)}
-x_{1+P^X_k(r),\pi_X(r)} & \text{otherwise}.
\end{cases}
\end{equation}
We also define
\begin{multline}
{\rm Comb}_{k,\iota}^{X}[\lambda]:=\label{comb-lam-2}\\
\begin{cases}
\{-x_{1,k}+\langle h_k,\lambda \rangle\} & \text{if }(1,k)<(1,\pi_X(k+1)),\ (1,k)<(1,\pi_X(k-1)),\\
\{\fbox{$\tilde{r}$}^{X}_k+\langle h_k,\lambda \rangle\ | r\in\mathbb{Z}_{\leq k}\}
 & \text{if }(1,k)<(1,\pi_X(k+1)),\ (1,k)>(1,\pi_X(k-1)),\\
\left\{\fbox{$r$}^{X}_k+\langle h_k,\lambda \rangle\ | r\in\mathbb{Z}_{\geq k+1}\right\} 
& \text{if }(1,k)>(1,\pi_X(k+1)),\ (1,k)<(1,\pi_X(k-1)),\\
\{
L^X_{0,k,\iota}(T)+\langle h_k,\lambda \rangle | T\in{\rm REYD}^{X}_{k}\setminus\{\phi^k\}
\} & \text{if }(1,k)>(1,\pi_X(k+1)),\ (1,k)>(1,\pi_X(k-1)).\\ 
\end{cases}
\end{multline}
Here, $\phi^k\in{\rm REYD}^{X}_{k}$ is the revised extended Young diagram without boxes, that is, $\phi^k=(y_{\ell})_{\ell\in\mathbb{Z}}$, 
$y_{\ell}=k$ for any $\ell\in\mathbb{Z}$.

Next, we take $k$ as $k\in I\setminus I_X$.
We draw Young walls on $\mathbb{R}_{\leq0}\times \mathbb{R}_{\geq k}$.
For example, the Young wall in Example \ref{ex-A2-1} is drawn as follow:
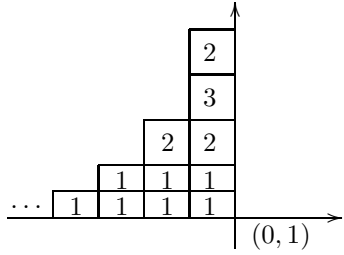
\begin{equation*}
\begin{xy}
(-15.5,-2) *{\ 1}="000-1",
(-9.5,-2) *{\ 1}="00-1",
(-3.5,-2) *{\ 1}="0-1",
(1.5,-2) *{\ \ 1}="0-1",
(-21.5,-2) *{\dots}="00000",
(-9.5,1.5) *{\ 1}="000",
(-3.5,1.5) *{\ 1}="00",
(1.5,1.5) *{\ \ 1}="0",
(1.5,6.5) *{\ \ 2}="02",
(-4,6.5) *{\ \ 2}="002",
(1.5,12.5) *{\ \ 3}="03",
(1.5,18.5) *{\ \ 2}="04",
(12,-6) *{(0,1)}="origin",
(0,0) *{}="1",
(0,-3.5) *{}="1-u",
(6,0)*{}="2",
(6,-3.5)*{}="2-u",
(6,-7.5)*{}="-1-u",
(6,3.5)*{}="3",
(0,3.5)*{}="4",
(-6,3.5)*{}="5",
(-6,0)*{}="6",
(-6,-3.5)*{}="6-u",
(-12,3.5)*{}="7",
(-12,0)*{}="8",
(-12,-3.6)*{}="8-u",
(-18,3.5)*{}="9",
(-18,0)*{}="10-a",
(-18,0)*{}="10",
(-18,-3.5)*{}="10-u",
(-24,-3.5)*{}="11-u",
(20,-3.5)*{}="0-u",
(6,9.5)*{}="3-1",
(6,15.5)*{}="3-2",
(6,21.5)*{}="3-3",
(6,25)*{}="3-4",
(0,9.5)*{}="4-1",
(0,15.5)*{}="4-2",
(0,21.5)*{}="4-3",
(-6,9.5)*{}="5-1",
\ar@{-} "-1-u";"2-u"^{}
\ar@{->} "2-u";"0-u"^{}
\ar@{-} "10-u";"11-u"^{}
\ar@{-} "8";"10-a"^{}
\ar@{-} "10-u";"10-a"^{}
\ar@{-} "10-u";"2-u"^{}
\ar@{-} "8";"8-u"^{}
\ar@{-} "6";"6-u"^{}
\ar@{-} "2";"2-u"^{}
\ar@{-} "1";"1-u"^{}
\ar@{-} "1";"2"^{}
\ar@{-} "1";"4"^{}
\ar@{-} "2";"3"^{}
\ar@{-} "3";"4"^{}
\ar@{-} "5";"6"^{}
\ar@{-} "5";"4"^{}
\ar@{-} "1";"6"^{}
\ar@{-} "7";"8"^{}
\ar@{-} "7";"5"^{}
\ar@{-} "6";"8"^{}
\ar@{-} "3";"3-1"^{}
\ar@{-} "4-1";"3-1"^{}
\ar@{-} "4-1";"4"^{}
\ar@{-} "5-1";"5"^{}
\ar@{-} "4-1";"5-1"^{}
\ar@{-} "4-1";"4-2"^{}
\ar@{-} "3-2";"3-1"^{}
\ar@{-} "3-2";"4-2"^{}
\ar@{-} "3-2";"3-3"^{}
\ar@{-} "4-3";"4-2"^{}
\ar@{-} "4-3";"3-3"^{}
\ar@{->} "3-3";"3-4"^{}
\end{xy}
\end{equation*}
Let $i\in\mathbb{Z}_{\geq0}$, $\ell\in\mathbb{Z}_{\geq k}$ and $S$ be a slot or block
\[
S=
\begin{xy}
(-8,15) *{(-i-1,\ell+1)}="0000",
(17,15) *{(-i,\ell+1)}="000",
(15,-3) *{(-i,\ell)}="00",
(-7,-3) *{(-i-1,\ell)}="0",
(0,0) *{}="1",
(12,0)*{}="2",
(12,12)*{}="3",
(0,12)*{}="4",
\ar@{-} "1";"2"^{}
\ar@{-} "1";"4"^{}
\ar@{-} "2";"3"^{}
\ar@{-} "3";"4"^{}
\end{xy}
\]
in $\mathbb{R}_{\leq0}\times \mathbb{R}_{\geq k}$.
If $S$ is colored by $t\in I$ in the pattern of Definition \ref{def-YW} (ii)
then we define
\begin{equation}\label{l1def1}
L^X_{s,k,{\rm ad}}(S):=x_{s+P^X_k(\ell)+i,t},\quad L^X_{s,k,{\rm re}}(S):=x_{s+P^X_k(\ell)+i+1,t}.
\end{equation}
Let $i\in\mathbb{Z}_{\geq0}$, $\ell\in\mathbb{Z}_{\geq k}$ and $S'$ be
a slot or block colored by $t=1$ or $n$ in $\mathbb{R}_{\leq 0}\times \mathbb{R}_{\geq k}$
such that the place is one of the following two: 
\[
S'=
\begin{xy}
(-8,10) *{(-i-1,\ell+\frac{1}{2})}="0000",
(17,10) *{(-i,\ell+\frac{1}{2})}="000",
(15,-3) *{(-i,\ell)}="00",
(-7,-3) *{(-i-1,\ell)}="0",
(0,0) *{}="1",
(12,0)*{}="2",
(12,6)*{}="3",
(0,6)*{}="4",
(32,3) *{{\rm or}}="or",
(52,10) *{(-i-1,\ell+1)}="a0000",
(77,10) *{(-i,\ell+1)}="a000",
(75,-3) *{(-i,\ell+\frac{1}{2})}="a00",
(54,-3) *{(-i-1,\ell+\frac{1}{2})}="a0",
(60,0) *{}="a1",
(72,0)*{}="a2",
(72,6)*{}="a3",
(60,6)*{}="a4",
\ar@{-} "1";"2"^{}
\ar@{-} "1";"4"^{}
\ar@{-} "2";"3"^{}
\ar@{-} "3";"4"^{}
\ar@{-} "a1";"a2"^{}
\ar@{-} "a1";"a4"^{}
\ar@{-} "a2";"a3"^{}
\ar@{-} "a3";"a4"^{}
\end{xy}
\]
Then we set
\begin{equation}\label{l1def2}
L^X_{s,k,{\rm ad}}(S'):=x_{s+P^X_k(\ell)+i,t},\quad L^X_{s,k,{\rm re}}(S'):=x_{s+P^X_k(\ell)+i+1,t}.
\end{equation}
If there exists $S$ or $S'$ colored by $t\in I$ as above,
it holds $\pi'(\ell)=t$, where $\pi'$ is defined as (\ref{piprime}).
For a proper Young wall $Y$ of type $X$ of ground state $\Lambda_k$, one defines
\begin{eqnarray}
L^X_{s,k,\iota}(Y)&:=&
\sum_{t\in I} \left(
\sum_{P:\text{single admissible slot}} L^X_{s,k,{\rm ad}}(P)
-\sum_{P:\text{single removable block}} L^X_{s,k,{\rm re}}(P)\right)\nonumber\\
& & + 
\sum_{P:\text{double admissible slot}} 2L^X_{s,k,{\rm ad}}(P)
-\sum_{P:\text{double removable block}} 2L^X_{s,k,{\rm re}}(P). \label{L11-def}
\end{eqnarray}
Let ${\rm YW}^{X}_{k}$ denote the set of all proper Young walls of type $X$ of ground state $\Lambda_k$.
We define
\begin{equation}\label{comb-lam-3}
{\rm Comb}_{k,\iota}^{X}[\lambda]:=
\begin{cases}
\{-x_{1,k}+\langle h_k,\lambda \rangle\} & \text{if }(1,k)<(1,\pi'(k+1)),\\
\{
L^X_{0,k,\iota}(Y)+\langle h_k,\lambda \rangle | Y\in{\rm YW}^{X}_{k}\setminus\{\phi^k\}
\} & \text{if }(1,k)>(1,\pi'(k+1)).
\end{cases}
\end{equation}
Here, $\phi^k=Y_{\Lambda_k}\in {\rm YW}^{X}_{k}$ is the ground state wall.
One also defines
\begin{equation}\label{comb-inf-2}
{\rm Comb}^{X}_{\iota}[\infty]
:=\{
L^{X}_{s,k,\iota}(T)| 
s\in\mathbb{Z}_{\geq1},\ k\in I_{X}\ 
\text{and}\ T\in{\rm REYD}^{X}_{k}
\}\cup
\{
L^{X}_{s,k,\iota}(Y)|
s\in\mathbb{Z}_{\geq1},\ k\in I\setminus I_{X}\text{ and }Y\in{\rm YW}^{X}_{k}
\}.
\end{equation}
\begin{thm}\label{thmA2}
We suppose that $\mathfrak{g}$ is of type $X={\rm A}_{2n-2}^{(2)}$ or $X={\rm C}_{n-1}^{(1)}$ $(n\geq3)$ and $\iota$ is adapted. 
Then 
for any $\lambda\in P^+$, the pair
$(\iota,\lambda)$ is $\Xi'$-ample and
\[
{\rm Im}(\Psi^{(\lambda)}_{\iota})
=
\left\{ 
\mathbf{a}\in\mathbb{Z}^{\infty} \left|
\varphi(\mathbf{a})\geq0\text{ for any }
\varphi\in
{\rm Comb}^{X^L}_{\iota}[\infty]\cup
\bigcup_{k\in I}
{\rm Comb}^{X^L}_{k,\iota}[\lambda]
\right.
\right\}.
\]
Here, $({\rm A}^{(2)}_{2n-2})^L={\rm A}^{(2)}_{2n-2}$, $({\rm C}^{(1)}_{n-1})^L={\rm D}^{(2)}_{n}$.
\end{thm}

Similarly to Remark \ref{rem-1}, we see that
functions in ${\rm Comb}^{X^L}_{\iota}[\infty]$ and ${\rm Comb}^{X^L}_{k,\iota}[\lambda]$
are related to representation theory of $U_q(^L\mathfrak{g})$.

\begin{ex}
We consider the case $\mathfrak{g}$ is of type ${\rm A}_4^{(2)}$ and $\iota=(\cdots,3,1,2,3,1,2)$, which
is the same setting as in Example 4.9 of \cite{Ka}.
Let us
compute several inequalities defining ${\rm Im}(\Psi_{\iota}^{(\lambda)})$.
In the same example of \cite{Ka},
we already found several inequalities defining
${\rm Im}(\Psi_{\iota})=\{\mathbf{a}\in \mathbb{Z}^{\infty} |
 \varphi(\mathbf{a})\geq0\ \text{ for any }\varphi\in {\rm Comb}^{{\rm A}^{(2)}}_{\iota}[\infty]\}$:
\[
{\rm Im}(\Psi_{\iota})
=\left\{ 
\textbf{a}\in\mathbb{Z}^{\infty} \left|
\begin{array}{l}
s\in\mathbb{Z}_{\geq1},\ a_{s,2}\geq0,
2a_{s,1}+a_{s,3}-a_{s+1,2}\geq0,
a_{s,1}+a_{s,3}-a_{s+1,1}\geq0,\\
a_{s,1}+2a_{s+1,2}-a_{s+1,3}-a_{s+1,1}\geq0,\\
a_{s,1}+a_{s+1,2}+a_{s+1,1}-a_{s+2,2}\geq0,
a_{s,1}+a_{s+1,2}-a_{s+2,1}\geq0,\cdots\\
a_{s,3}\geq0,\ 
2a_{s+1,2}-a_{s+1,3}\geq0,
2a_{s+1,1}+a_{s+1,2}-a_{s+2,2}\geq0,\cdots\\
a_{s,1}\geq0,
a_{s+1,2}-a_{s+1,1}\geq0,
a_{s+1,3}+a_{s+1,1}-a_{s+2,2}\geq0,\\
a_{s+1,3}-a_{s+2,1}\geq0,
a_{s+2,2}+a_{s+1,1}-a_{s+2,3}\geq0,\cdots
\end{array} \right.
\right\}.
\]
Let us compute several inequalities coming from ${\rm Comb}_{k,\iota}^{{\rm A}^{(2)}}[\lambda]$ for $k=1,2,3$.
We obtain
\[
P_1^{{\rm A}^{(2)}}(1)=0,\ P_1^{{\rm A}^{(2)}}(2)=1,\ P_1^{{\rm A}^{(2)}}(3)=1,\ P_1^{{\rm A}^{(2)}}(4)=2,\cdots,
\]
\[
\cdots,\ P_2^{{\rm A}^{(2)}}(-1)=1,\ P_2^{{\rm A}^{(2)}}(0)=0,\ P_2^{{\rm A}^{(2)}}(1)=0,\ P_2^{{\rm A}^{(2)}}(2)=0,\ P_2^{{\rm A}^{(2)}}(3)=0,\ P_2^{{\rm A}^{(2)}}(4)=1,\cdots,
\]
\[
\cdots,\ P_3^{{\rm A}^{(2)}}(1)=1,\ P_3^{{\rm A}^{(2)}}(2)=1,\ P_3^{{\rm A}^{(2)}}(3)=0,\ P_3^{{\rm A}^{(2)}}(4)=1,\ P_3^{{\rm A}^{(2)}}(5)=1,\cdots.
\]
It follows from definitions that
${\rm Comb}_{2,\iota}^{{\rm A}^{(2)}}[\lambda]=\{-x_{1,2}+\langle h_2,\lambda\rangle\}$.
Next, we consider ${\rm Comb}_{1,\iota}^{{\rm A}^{(2)}}[\lambda]$.
The following Young walls are elements in ${\rm YW}^{{\rm A}^{(2)}}_{1}\setminus\{\phi^1\}$:
\[
\begin{xy}
(30.5,6) *{Y_1:=}="0000-1a",
(52.5,1.5) *{\ 1}="000-1-1a",
(34.5,-2) *{\ 1}="000-1a",
(40.5,-2) *{\ 1}="00-1a",
(46.5,-2) *{\ 1}="0-1a",
(52.5,-2) *{\ 1}="0-1sa",
(29.5,-2) *{\dots}="00000a",
(62,-6) *{(0,1)}="origina",
(50,0) *{}="1a",
(56,18) *{}="ya",
(56,3.5)*{}="6+35a",
(50,3.5)*{}="0+35a",
(50,0)*{}="0+0a",
(50,-3.5) *{}="1-ua",
(56,0)*{}="2a",
(56,-3.5)*{}="2-ua",
(56,-7.5)*{}="-1-ua",
(44,0)*{}="6a",
(44,-3.5)*{}="6-ua",
(38,0)*{}="8a",
(38,-3.6)*{}="8-ua",
(32,3.5)*{}="9a",
(32,0)*{}="10-ta",
(32,0)*{}="10ta",
(32,-3.5)*{}="10-uta",
(27,-3.5)*{}="11-uta",
(70,-3.5)*{}="0-ua",
(80.5,6) *{Y_2:=}="0000-1b",
(102.5,1.5) *{\ 1}="000-1-1b",
(102.5,6.5) *{\ 2}="000-1-1-2b",
(84.5,-2) *{\ 1}="000-1b",
(90.5,-2) *{\ 1}="00-1b",
(96.5,-2) *{\ 1}="0-1b",
(102.5,-2) *{\ 1}="0-1sb",
(79.5,-2) *{\dots}="00000b",
(112,-6) *{(0,1)}="originb",
(100,0) *{}="1b",
(106,18) *{}="yb",
(106,9.5)*{}="6+95b",
(100,9.5)*{}="0+95b",
(106,3.5)*{}="6+35b",
(100,3.5)*{}="0+35b",
(100,0)*{}="0+0b",
(100,-3.5) *{}="1-ub",
(106,0)*{}="2b",
(106,-3.5)*{}="2-ub",
(106,-7.5)*{}="-1-ub",
(94,0)*{}="6b",
(94,-3.5)*{}="6-ub",
(88,0)*{}="8b",
(88,-3.6)*{}="8-ub",
(82,3.5)*{}="9b",
(82,0)*{}="10-tb",
(82,0)*{}="10tb",
(82,-3.5)*{}="10-utb",
(77,-3.5)*{}="11-utb",
(120,-3.5)*{}="0-ub",
\ar@{-} "-1-ua";"2-ua"^{}
\ar@{->} "2-ua";"0-ua"^{}
\ar@{->} "2a";"ya"^{}
\ar@{-} "10-uta";"11-uta"^{}
\ar@{-} "8a";"10-ta"^{}
\ar@{-} "10-uta";"10-ta"^{}
\ar@{-} "10-uta";"2-ua"^{}
\ar@{-} "8a";"8-ua"^{}
\ar@{-} "6a";"6-ua"^{}
\ar@{-} "2a";"2-ua"^{}
\ar@{-} "1a";"1-ua"^{}
\ar@{-} "1a";"2a"^{}
\ar@{-} "1a";"6a"^{}
\ar@{-} "6a";"8a"^{}
\ar@{-} "6+35a";"0+35a"^{}
\ar@{-} "0+0a";"0+35a"^{}
\ar@{-} "-1-ub";"2-ub"^{}
\ar@{->} "2-ub";"0-ub"^{}
\ar@{->} "2b";"yb"^{}
\ar@{-} "10-utb";"11-utb"^{}
\ar@{-} "8b";"10-tb"^{}
\ar@{-} "10-utb";"10-tb"^{}
\ar@{-} "10-utb";"2-ub"^{}
\ar@{-} "8b";"8-ub"^{}
\ar@{-} "6b";"6-ub"^{}
\ar@{-} "2b";"2-ub"^{}
\ar@{-} "1b";"1-ub"^{}
\ar@{-} "1b";"2b"^{}
\ar@{-} "1b";"6b"^{}
\ar@{-} "6b";"8b"^{}
\ar@{-} "6+35b";"0+35b"^{}
\ar@{-} "0+0b";"0+35b"^{}
\ar@{-} "0+95b";"6+95b"^{}
\ar@{-} "0+35b";"0+95b"^{}
\end{xy}
\]
\[
\begin{xy}
(30.5,6) *{Y_3:=}="0000-1a",
(52.5,1.5) *{\ 1}="000-1-1a",
(52.5,6.5) *{\ 2}="000-1-1-2a",
(34.5,-2) *{\ 1}="000-1a",
(40.5,-2) *{\ 1}="00-1a",
(46.5,-2) *{\ 1}="0-1a",
(46.5,1.5) *{\ 1}="0-1-1a",
(52.5,-2) *{\ 1}="0-1sa",
(29.5,-2) *{\dots}="00000a",
(62,-6) *{(0,1)}="origina",
(50,0) *{}="1a",
(56,18) *{}="ya",
(56,9.5)*{}="6+95a",
(50,9.5)*{}="0+95a",
(56,3.5)*{}="6+35a",
(50,3.5)*{}="0+35a",
(50,0)*{}="0+0a",
(50,-3.5) *{}="1-ua",
(56,0)*{}="2a",
(56,-3.5)*{}="2-ua",
(56,-7.5)*{}="-1-ua",
(44,3.5)*{}="-6+35a",
(44,0)*{}="6a",
(44,-3.5)*{}="6-ua",
(38,0)*{}="8a",
(38,-3.6)*{}="8-ua",
(32,3.5)*{}="9a",
(32,0)*{}="10-ta",
(32,0)*{}="10ta",
(32,-3.5)*{}="10-uta",
(27,-3.5)*{}="11-uta",
(70,-3.5)*{}="0-ua",
(80.5,6) *{Y_4:=}="0000-1b",
(102.5,1.5) *{\ 1}="000-1-1b",
(102.5,6.5) *{\ 2}="000-1-1-2b",
(102.5,12.5) *{\ 3}="000-1-1-2-3b",
(84.5,-2) *{\ 1}="000-1b",
(90.5,-2) *{\ 1}="00-1b",
(96.5,-2) *{\ 1}="0-1b",
(102.5,-2) *{\ 1}="0-1sb",
(79.5,-2) *{\dots}="00000b",
(112,-6) *{(0,1)}="originb",
(100,0) *{}="1b",
(106,18) *{}="yb",
(106,15.5)*{}="6+155b",
(100,15.5)*{}="0+155b",
(106,9.5)*{}="6+95b",
(100,9.5)*{}="0+95b",
(106,3.5)*{}="6+35b",
(100,3.5)*{}="0+35b",
(100,0)*{}="0+0b",
(100,-3.5) *{}="1-ub",
(106,0)*{}="2b",
(106,-3.5)*{}="2-ub",
(106,-7.5)*{}="-1-ub",
(94,0)*{}="6b",
(94,-3.5)*{}="6-ub",
(88,0)*{}="8b",
(88,-3.6)*{}="8-ub",
(82,3.5)*{}="9b",
(82,0)*{}="10-tb",
(82,0)*{}="10tb",
(82,-3.5)*{}="10-utb",
(77,-3.5)*{}="11-utb",
(120,-3.5)*{}="0-ub",
\ar@{-} "-1-ua";"2-ua"^{}
\ar@{->} "2-ua";"0-ua"^{}
\ar@{->} "2a";"ya"^{}
\ar@{-} "10-uta";"11-uta"^{}
\ar@{-} "8a";"10-ta"^{}
\ar@{-} "10-uta";"10-ta"^{}
\ar@{-} "10-uta";"2-ua"^{}
\ar@{-} "8a";"8-ua"^{}
\ar@{-} "6a";"6-ua"^{}
\ar@{-} "2a";"2-ua"^{}
\ar@{-} "1a";"1-ua"^{}
\ar@{-} "1a";"2a"^{}
\ar@{-} "1a";"6a"^{}
\ar@{-} "6a";"8a"^{}
\ar@{-} "-6+35a";"0+35a"^{}
\ar@{-} "-6+35a";"6a"^{}
\ar@{-} "6+35a";"0+35a"^{}
\ar@{-} "0+0a";"0+35a"^{}
\ar@{-} "0+95a";"6+95a"^{}
\ar@{-} "0+35a";"0+95a"^{}
\ar@{-} "-1-ub";"2-ub"^{}
\ar@{->} "2-ub";"0-ub"^{}
\ar@{->} "2b";"yb"^{}
\ar@{-} "10-utb";"11-utb"^{}
\ar@{-} "8b";"10-tb"^{}
\ar@{-} "10-utb";"10-tb"^{}
\ar@{-} "10-utb";"2-ub"^{}
\ar@{-} "8b";"8-ub"^{}
\ar@{-} "6b";"6-ub"^{}
\ar@{-} "2b";"2-ub"^{}
\ar@{-} "1b";"1-ub"^{}
\ar@{-} "1b";"2b"^{}
\ar@{-} "1b";"6b"^{}
\ar@{-} "6b";"8b"^{}
\ar@{-} "6+35b";"0+35b"^{}
\ar@{-} "0+0b";"0+35b"^{}
\ar@{-} "0+95b";"6+95b"^{}
\ar@{-} "0+35b";"0+95b"^{}
\ar@{-} "0+155b";"6+155b"^{}
\ar@{-} "0+95b";"0+155b"^{}
\end{xy}
\]
We see that $Y_1$ has a single $2$-admissible slot
and a single removable $1$-block so that $L^{{\rm A}^{(2)}}_{0,1,\iota}(Y_1)=x_{P_1^{{\rm A}^{(2)}}(2),2}-
x_{P_1^{{\rm A}^{(2)}}(1)+1,1}
=x_{1,2}-x_{1,1}$. Similarly,
\[
L^{{\rm A}^{(2)}}_{0,1,\iota}(Y_2)=x_{P_1^{{\rm A}^{(2)}}(3),3}+x_{P_1^{{\rm A}^{(2)}}(1)+1,1}-x_{P_1^{{\rm A}^{(2)}}(2)+1,2}=x_{1,3}+x_{1,1}-x_{2,2}, \]
\[
L^{{\rm A}^{(2)}}_{0,1,\iota}(Y_3)=x_{P_1^{{\rm A}^{(2)}}(3),3}-x_{P_1^{{\rm A}^{(2)}}(1)+2,1}=x_{1,3}-x_{2,1},
\]
\[
L^{{\rm A}^{(2)}}_{0,1,\iota}(Y_4)=x_{P_1^{{\rm A}^{(2)}}(4),2}+x_{P_1^{{\rm A}^{(2)}}(1)+1,1}-x_{P_1^{{\rm A}^{(2)}}(3)+1,3}=x_{2,2}+x_{1,1}-x_{2,3}.
\]
Next, we consider ${\rm Comb}_{3,\iota}^{{\rm A}^{(2)}}[\lambda]$.
The following diagrams are elements in ${\rm REYD}^{{\rm A}^{(2)}}_{3}\setminus\{\phi^3\}$:
\[
\begin{xy}
(40,-9) *{T_3^1:=}="Y1",
(55,0) *{}="1a",
(85,0)*{}="2a",
(50,0)*{(0,3)}="4a",
(61,2) *{1}="10a",
(61,-1) *{}="1010a",
(61,-6) *{}="1010-1a",
(67,2) *{2}="11a",
(67,-1) *{}="1111a",
(73,2) *{3}="12a",
(73,-1) *{}="1212a",
(79,2) *{4}="13a",
(79,-1) *{}="1313a",
(55,-6)*{}="5a",
(49,-6)*{}="55a",
(49,-12)*{}="6a",
(43,-12)*{}="66a",
(43,-18)*{}="7a",
(37,-18)*{}="77a",
(37,-24)*{}="8a",
(31,-24)*{}="88a",
(95,-9) *{T_3^2:=}="Y2",
(110,0) *{}="1b",
(110,-12) *{}="1b-12",
(140,0)*{}="2b",
(105,0)*{(0,3)}="4b",
(116,2) *{1}="10b",
(116,-1) *{}="1010b",
(116,-6) *{}="1010-1b",
(122,2) *{2}="11b",
(122,-1) *{}="1111b",
(128,2) *{3}="12b",
(128,-1) *{}="1212b",
(134,2) *{4}="13b",
(134,-1) *{}="1313b",
(110,-6)*{}="5b",
(104,-6)*{}="55b",
(104,-12)*{}="6b",
(98,-12)*{}="66b",
(98,-18)*{}="7b",
(92,-18)*{}="77b",
(92,-24)*{}="8b",
(86,-24)*{}="88b",
\ar@{-} "1a";"2a"^{}
\ar@{-} "1a";"5a"^{}
\ar@{-} "5a";"55a"^{}
\ar@{-} "55a";"6a"^{}
\ar@{-} "6a";"66a"^{}
\ar@{-} "66a";"7a"^{}
\ar@{-} "7a";"77a"^{}
\ar@{-} "77a";"8a"^{}
\ar@{-} "8a";"88a"^{}
\ar@{-} "10a";"1010a"^{}
\ar@{-} "1010a";"1010-1a"^{}
\ar@{-} "5a";"1010-1a"^{}
\ar@{-} "11a";"1111a"^{}
\ar@{-} "12a";"1212a"^{}
\ar@{-} "13a";"1313a"^{}
\ar@{-} "1b";"2b"^{}
\ar@{-} "1b";"1b-12"^{}
\ar@{-} "66b";"1b-12"^{}
\ar@{-} "1b";"5b"^{}
\ar@{-} "5b";"55b"^{}
\ar@{-} "55b";"6b"^{}
\ar@{-} "6b";"66b"^{}
\ar@{-} "66b";"7b"^{}
\ar@{-} "7b";"77b"^{}
\ar@{-} "77b";"8b"^{}
\ar@{-} "8b";"88b"^{}
\ar@{-} "10b";"1010b"^{}
\ar@{-} "1010b";"1010-1b"^{}
\ar@{-} "5b";"1010-1b"^{}
\ar@{-} "11b";"1111b"^{}
\ar@{-} "12b";"1212b"^{}
\ar@{-} "13b";"1313b"^{}
\end{xy}
\]
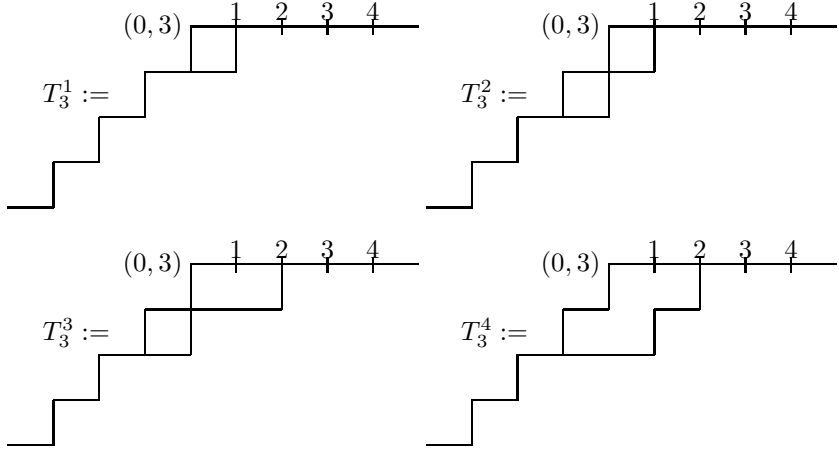
\[
\begin{xy}
(-15,-9) *{T_3^3:=}="Y3",
(0,0) *{}="1",
(30,0)*{}="2",
(-5,0)*{(0,3)}="4",
(6,2) *{1}="10",
(6,-1) *{}="1010",
(6,-6) *{}="1010-1",
(12,2) *{2}="11",
(12,-6) *{}="12-6",
(12,-1) *{}="1111",
(18,2) *{3}="12",
(18,-1) *{}="1212",
(24,2) *{4}="13",
(24,-1) *{}="1313",
(0,-6)*{}="5",
(0,-12) *{}="0-12",
(-6,-6)*{}="55",
(-6,-12)*{}="6",
(-6,-18)*{}="6-18",
(0,-18)*{}="0-18",
(-12,-12)*{}="66",
(-12,-18)*{}="7",
(-18,-18)*{}="77",
(-18,-24)*{}="8",
(-24,-24)*{}="88",
(40,-9) *{T_3^4:=}="Y4",
(55,0) *{}="1a",
(85,0)*{}="2a",
(50,0)*{(0,3)}="4a",
(61,2) *{1}="10a",
(61,-1) *{}="1010a",
(61,-6) *{}="6-6a",
(61,-12) *{}="6-12a",
(67,-6) *{}="1010-1a",
(67,2) *{2}="11a",
(67,-1) *{}="1111a",
(73,2) *{3}="12a",
(73,-1) *{}="1212a",
(79,2) *{4}="13a",
(79,-1) *{}="1313a",
(55,-6)*{}="5a",
(49,-6)*{}="55a",
(49,-12)*{}="6a",
(43,-12)*{}="66a",
(43,-18)*{}="7a",
(37,-18)*{}="77a",
(37,-24)*{}="8a",
(31,-24)*{}="88a",
\ar@{-} "1";"2"^{}
\ar@{-} "1";"5"^{}
\ar@{-} "5";"55"^{}
\ar@{-} "55";"6"^{}
\ar@{-} "6";"66"^{}
\ar@{-} "5";"0-12"^{}
\ar@{-} "6";"0-12"^{}
\ar@{-} "66";"7"^{}
\ar@{-} "7";"77"^{}
\ar@{-} "77";"8"^{}
\ar@{-} "8";"88"^{}
\ar@{-} "10";"1010"^{}
\ar@{-} "11";"1111"^{}
\ar@{-} "11";"12-6"^{}
\ar@{-} "12";"1212"^{}
\ar@{-} "13";"1313"^{}
\ar@{-} "12-6";"1010-1"^{}
\ar@{-} "5";"1010-1"^{}
\ar@{-} "1a";"2a"^{}
\ar@{-} "1a";"5a"^{}
\ar@{-} "5a";"55a"^{}
\ar@{-} "55a";"6a"^{}
\ar@{-} "6a";"66a"^{}
\ar@{-} "66a";"7a"^{}
\ar@{-} "7a";"77a"^{}
\ar@{-} "77a";"8a"^{}
\ar@{-} "8a";"88a"^{}
\ar@{-} "10a";"1010a"^{}
\ar@{-} "1111a";"1010-1a"^{}
\ar@{-} "6-12a";"6-6a"^{}
\ar@{-} "6-12a";"6a"^{}
\ar@{-} "1010-1a";"6-6a"^{}
\ar@{-} "11a";"1111a"^{}
\ar@{-} "12a";"1212a"^{}
\ar@{-} "13a";"1313a"^{}
\end{xy}
\]
One can compute
\[
L^{{\rm A}^{(2)}}_{0,3,\iota}(T_3^1)=x_{P_3^{{\rm A}^{(2)}}(4),2}+x_{P_3^{{\rm A}^{(2)}}(2),2}-
x_{P_3^{{\rm A}^{(2)}}(3)+1,3}
=2x_{1,2}-x_{1,3},\]
\[
L^{{\rm A}^{(2)}}_{0,3,\iota}(T_3^2)=2x_{P_3^{{\rm A}^{(2)}}(1),1}+x_{P_3^{{\rm A}^{(2)}}(4),2}-x_{P_3^{{\rm A}^{(2)}}(2)+1,2}
=2x_{1,1}+x_{1,2}-x_{2,2},
\]
\[
L^{{\rm A}^{(2)}}_{0,3,\iota}(T_3^3)=2x_{P_3^{{\rm A}^{(2)}}(1),1}+
x_{P_3^{{\rm A}^{(2)}}(3)+1,3}+
2x_{P_3^{{\rm A}^{(2)}}(5),1}
-x_{P_3^{{\rm A}^{(2)}}(4)+1,2}
-x_{P_3^{{\rm A}^{(2)}}(2)+1,2}
=4x_{1,1}+x_{1,3}
-2x_{2,2},
\]
\[
L^{{\rm A}^{(2)}}_{0,3,\iota}(T_3^4)=2x_{P_3^{{\rm A}^{(2)}}(1),1}+
2x_{P_3^{{\rm A}^{(2)}}(5),1}
-x_{P_3^{{\rm A}^{(2)}}(3)+2,3}
=4x_{1,1}-x_{2,3}.
\]
Note that in $T_3^2$, $T_3^3$ and $T_3^4$, the point $(-2,1)$ is double $1$-admissible
and
in $T_3^3$ and $T_3^4$, the point $(2,3)$ is double $1$-admissible. 
Hence, we get several inequalities defining ${\rm Im}(\Psi_{\iota}^{(\lambda)})$:
\[
{\rm Im}(\Psi_{\iota}^{(\lambda)})
=\left\{ 
\textbf{a}\in\mathbb{Z}^{\infty} \left|
\begin{array}{l}
s\in\mathbb{Z}_{\geq1},\ a_{s,2}\geq0,
2a_{s,1}+a_{s,3}-a_{s+1,2}\geq0,
a_{s,1}+a_{s,3}-a_{s+1,1}\geq0,\\
a_{s,1}+2a_{s+1,2}-a_{s+1,3}-a_{s+1,1}\geq0,\\
a_{s,1}+a_{s+1,2}+a_{s+1,1}-a_{s+2,2}\geq0,
a_{s,1}+a_{s+1,2}-a_{s+2,1}\geq0,\cdots\\
a_{s,3}\geq0,\ 
2a_{s+1,2}-a_{s+1,3}\geq0,
2a_{s+1,1}+a_{s+1,2}-a_{s+2,2}\geq0,\cdots\\
a_{s,1}\geq0,
a_{s+1,2}-a_{s+1,1}\geq0,
a_{s+1,3}+a_{s+1,1}-a_{s+2,2}\geq0,\\
a_{s+1,3}-a_{s+2,1}\geq0,
a_{s+2,2}+a_{s+1,1}-a_{s+2,3}\geq0,\cdots\\
a_{1,2}\leq \langle h_2,\lambda \rangle,
a_{1,2}-a_{1,1}+\langle h_1,\lambda \rangle\geq0,
a_{1,3}+a_{1,1}-a_{2,2}+\langle h_1,\lambda \rangle\geq0,\\
a_{1,3}-a_{2,1}+\langle h_1,\lambda \rangle\geq0,
a_{2,2}+a_{1,1}-a_{2,3}+\langle h_1,\lambda \rangle\geq0,\\
2a_{1,2}-a_{1,3}+\langle h_3,\lambda \rangle\geq0,
2a_{1,1}+a_{1,2}-a_{2,2}+\langle h_3,\lambda \rangle\geq0,\\
4a_{1,1}+a_{1,3}-2a_{2,2}+\langle h_3,\lambda \rangle\geq0,
4a_{1,1}-a_{2,3}+\langle h_3,\lambda \rangle\geq0,\cdots
\end{array} \right.
\right\}.
\]
\end{ex}

\subsection{Combinatorial description of $\varepsilon_k^*$}

\begin{thm}\label{thm3}
We suppose that $\mathfrak{g}$ is of type $X={\rm A}_{n-1}^{(1)}$ $(n\geq2)$, ${\rm D}_{n}^{(2)}$,
${\rm A}_{2n-2}^{(2)}$ or ${\rm C}_{n-1}^{(1)}$ $(n\geq3)$ and $\iota$ is adapted. 
Then for $k\in I$ and $x\in{\rm Im}(\Psi_{\iota})$, we have
\[
\varepsilon_k^*(x)
={\rm max}\{-\varphi(x) | \varphi\in{\rm Comb}_{k,\iota}^{X^L}[0]\}.
\]
\end{thm}

\begin{ex}
Let
$\mathfrak{g}$ be of type ${\rm A}_2^{(1)}$ and $\iota=(\cdots,3,1,2,3,1,2)$, which is the same setting
as in Example \ref{ex-2}. In Example 4.4 of \cite{Ka},
we shown that
for any $a_{2,3},a_{2,1},a_{2,2},a_{1,3},a_{1,1},a_{1,2}\in\mathbb{Z}_{\geq0}$ such that
$a_{1,1}-a_{2,2}\geq0$, $a_{1,3}-a_{2,1}\geq0$, $a_{2,2}-a_{2,3}\geq0$ and $a_{2,1}-a_{2,3}\geq0$,
it follows
$(\cdots,0,0,0,a_{2,3},a_{2,1},a_{2,2},a_{1,3},a_{1,1},a_{1,2})\in{\rm Im}(\Psi_{\iota})$.
Taking $x\in {\rm Im}(\Psi_{\iota})$ as
\begin{equation}\label{x-def}
x=(\cdots,0,0,1,2,3,2,3,3),
\end{equation}
let us confirm Theorem \ref{thm3}. By Section 2 of \cite{NZ},
we have
\[
x^*=
\tilde{f}_2^3 \tilde{f}_1^3 \tilde{f}_3^2 \tilde{f}_2^3 \tilde{f}_1^2 \tilde{f}_3   
(\cdots,0,0,0,0,0,0)
=(\cdots,0,0,1,0,1,1,1,1,1,1,2,2,1,1,1)
\]
and $\varepsilon_1^*(x)=1$, $\varepsilon_2^*(x)=3$
and $\varepsilon_3^*(x)=0$.

As seen in Example \ref{ex-2}, it holds
${\rm Comb}_{2,\iota}^{{\rm A}^{(1)}}[0]=\{-x_{1,2}\}$. Thus, 
${\rm max}\{-\varphi(x) | \varphi\in{\rm Comb}_{2,\iota}^{{\rm A}^{(1)}}[0]\}={\rm max}\{3,0\}=3$.
We also know that
${\rm Comb}_{1,\iota}^{{\rm A}^{(1)}}[0]=
\{\fbox{$r$}^{{\rm A}^{(1)}}_1| r\in\mathbb{Z}_{\geq 2}\}$.
By $P_1^{{\rm A}^{(1)}}(4)=2$, $P_1^{{\rm A}^{(1)}}(m)\geq3$ for $m\geq5$ and
(\ref{r-tilde-A0}), it follows
$\fbox{$r$}^{{\rm A}^{(1)}}_1(x)=0$ for $r\in \mathbb{Z}_{\geq5}$.
Hence
\begin{eqnarray*}
{\rm max}\{-\varphi(x) | \varphi\in{\rm Comb}_{1,\iota}^{{\rm A}^{(1)}}[0]\}&=&
{\rm max}
\{
-\fbox{$2$}^{{\rm A}^{(1)}}_1(x),\ -\fbox{$3$}^{{\rm A}^{(1)}}_1(x),\ -\fbox{$4$}^{{\rm A}^{(1)}}_1(x)
\}\\
&=&
{\rm max}\{0,1,-1\}=1.
\end{eqnarray*}
In this way, we can confirm 
$\varepsilon_k^*(x)
={\rm max}\{-\varphi(x) | \varphi\in{\rm Comb}_{k,\iota}^{{\rm A}^{(1)}}[0]\}$ for $k=1,2$.

Next, we consider the case $k=3$. The definition (\ref{comb-lam-1}) says 
${\rm Comb}_{3,\iota}^{{\rm A}^{(1)}}[0]=\{L^{{\rm A}^{(1)}}_{0,3,\iota}(T) | T\in {\rm EYD}_3\setminus\{\phi^3\} \}$.
Note that
for $\varphi = \sum_{m\in \mathbb{Z}_{\geq1},j\in I} c_{m,j}x_{m,j}\in {\rm Comb}_{3,\iota}^{{\rm A}^{(1)}}[0]$,
if $c_{m,j}\geq0$ for all $m\in\{1,2\}$, $j\in I$ then
$-\varphi(x)\leq 0$ by (\ref{x-def}).
By $P^{{\rm A}^{(1)}}_3(3)=0$, $P^{{\rm A}^{(1)}}_3(4)=1$, $P^{{\rm A}^{(1)}}_3(t)\geq2$ for $t\in\mathbb{Z}_{\geq5}$
and $P^{{\rm A}^{(1)}}_3(2)=P^{{\rm A}^{(1)}}_3(1)=P^{{\rm A}^{(1)}}_3(0)=1$, $P^{{\rm A}^{(1)}}_3(t)\geq2$
for $t\in\mathbb{Z}_{\leq-1}$, if $(i,j)$ is a convex corner of $T\in {\rm EYD}_3$ and
$P^{{\rm A}^{(1)}}_3(i+j)+{\rm min}\{3-j,i\}\leq2$ then
\[
(i,j)=(1,2),\ (1,1),\ (1,0),\ (1,-1),\ (2,2)\ \text{or } (2,1).
\]
If $T$ has a convex corner $(i,j)=(1,2)$ then $T=T_{3}^1$ in Example \ref{ex-2} and
$L^{{\rm A}^{(1)}}_{0,3,\iota}(T^1_3)=x_{1,1}+x_{1,2}-x_{1,3}$ and
$-L^{{\rm A}^{(1)}}_{0,3,\iota}(T^1_3)(x)= 2-3-3 \leq0$.

If $T$ has a convex corner $(i,j)=(1,1)$ (resp. $(1,0)$, $(1,-1)$, $(2,2)$) then 
the point $(0,1)$ (resp. $(0,0)$, $(0,-1)$, $(2,3)$) is a concave corner in $T$.
Note that $L^{{\rm A}^{(1)}}_{0,3,\iota}(1,1)=x_{2,2}$,
$L^{{\rm A}^{(1)}}_{0,3,\iota}(1,0)=x_{2,1}$,
$L^{{\rm A}^{(1)}}_{0,3,\iota}(1,-1)=x_{2,3}$,
$L^{{\rm A}^{(1)}}_{0,3,\iota}(2,2)=x_{2,1}$
and
$L^{{\rm A}^{(1)}}_{0,3,\iota}(0,1)=x_{1,1}$,
$L^{{\rm A}^{(1)}}_{0,3,\iota}(0,0)=x_{1,3}$,
$L^{{\rm A}^{(1)}}_{0,3,\iota}(0,-1)=x_{2,2}$,
$L^{{\rm A}^{(1)}}_{0,3,\iota}(2,3)=x_{2,2}$.
If $T$ has a convex corner $(i,j)=(2,1)$ then either $(2,2)$
or $(2,3)$ is a concave corner in $T$ and 
$L^{{\rm A}^{(1)}}_{0,3,\iota}(2,1)=x_{2,3}$.
By the above argument and (\ref{x-def}), we see that
${\rm max}\{-\varphi(x) | \varphi\in{\rm Comb}_{3,\iota}^{{\rm A}^{(1)}}[0]\}=0=\varepsilon_3^*(x)$.

\end{ex}

\section{Action of $\what{S}'$}

For the proof of theorems in the previous section,
we recall action of $\what{S}'$ on (revised) extended Young diagrams and Young walls.


\subsection{On extended Young diagrams}

In the next proposition,
for two extended Young diagrams $T=(y_r)_{r\in\mathbb{Z}_{\geq0}},\ T'=(y_r')_{r\in\mathbb{Z}_{\geq0}} \in {\rm EYD}_{k}$,
we consider the setting a point $(i,j)$ is a concave corner in $T$ and 
$T'$ is obtained by 
the following replacement of the concave corner
with a convex corner
and calculate how the values $L^X_{s,k,\iota}(T)$ (defined in (\ref{ovlL2})) are changed.
\begin{equation}\label{change}
\begin{xy}
(-6,13) *{(i,j)}="0000",
(4,-3) *{(i,j-1)}="000",
(10,15) *{(i+1,j)}="00",
(72,0) *{(i+1,j-1)}="a0000",
(54,-3) *{(i,j-1)}="a000",
(60,15) *{(i+1,j)}="a00",
(30,8) *{\rightarrow}="a00",
(0,0) *{\bullet}="1",
(12,12)*{\bullet}="3",
(0,12)*{\bullet}="4",
(50,0) *{\bullet}="a1",
(62,12)*{\bullet}="a3",
(62,0)*{\bullet}="a4",
\ar@{-} "a1";"a4"^{}
\ar@{-} "a3";"a4"^{}
\ar@{-} "1";"4"^{}
\ar@{-} "3";"4"^{}
\end{xy}
\end{equation}

\begin{prop}{\rm \cite{Ka}}\label{prop-closednessAD}
We suppose that $\mathfrak{g}$ is of type $X={\rm A}_{n-1}^{(1)}$ or $X={\rm D}_{n}^{(2)}$.
Let $T=(y_r)_{r\in\mathbb{Z}_{\geq0}},\ T'=(y_r')_{r\in\mathbb{Z}_{\geq0}} \in {\rm EYD}_{k}$
and put $j:=y_i$ with some $i\in\mathbb{Z}_{\geq0}$.
We suppose that
\[
y_i'=y_i-1(=j-1),\quad y_r'=y_r\ (r\neq i)
\]
and take $s\in\mathbb{Z}_{\geq0}$. If $s+P_k^{X^L}(i+j)+{\rm min}\{i,k-j\}\geq1$
then
\[
L^{X^L}_{s,k,\iota}(T')=L^{X^L}_{s,k,\iota}(T)-\beta_{s+P_k^{X^L}(i+j)+{\rm min}\{i,k-j\},\pi_{X^L}(i+j)}.
\]
\end{prop}

\subsection{On revised extended Young diagrams}

In this subsection, we assume $k\in I_X$.

\begin{prop}\label{A2closed}\cite{Ka}
We assume $\mathfrak{g}$ is of type ${\rm A}^{(2)}_{2n-2}$.
Let $T=(y_t)_{t\in\mathbb{Z}}$ be a sequence in ${\rm REYD}^{{\rm A}^{(2)}}_{k}$ (Definition \ref{AEYD}), $s\in\mathbb{Z}_{\geq0}$
and $i\in\mathbb{Z}$.
\begin{enumerate}
\item We assume that the point $(i,y_i)$ is single or double admissible and a sequence
$T'=(y_t')_{t\in\mathbb{Z}}\in{\rm REYD}^{{\rm A}^{(2)}}_{k}$
satisfies $y_{i}'=y_i-1$ and $y_t'=y_t$ $(t\neq i)$. Putting $j:=y_i$, if $s+P_k^{{\rm A}^{(2)}}(i+k)+[i]_-+k-j\geq1$ then
\[
L^{{\rm A}^{(2)}}_{s,k,\iota}(T')-L^{{\rm A}^{(2)}}_{s,k,\iota}(T)=-\beta_{s+P_k^{{\rm A}^{(2)}}(i+k)+[i]_-+k-j,\pi_{{\rm A}^{(2)}}(i+k)}.
\]
\item We assume that the point $(i,y_{i-1})$ is single or double removable and a sequence $T''=(y_t'')_{t\in\mathbb{Z}}\in {\rm REYD}^{{\rm A}^{(2)}}_{k}$ satisfies $y_{i-1}''=y_{i-1}+1$ and $y_t''=y_t$ $(t\neq i-1)$. Putting $j:=y_{i-1}''=y_{i-1}+1$,
if $s+P_k^{{\rm A}^{(2)}}(i+k-1)+[i-1]_-+k-j\geq1$ 
then
\[
L^{{\rm A}^{(2)}}_{s,k,\iota}(T'')-L^{{\rm A}^{(2)}}_{s,k,\iota}(T)=\beta_{s+P_k^{{\rm A}^{(2)}}(i+k-1)+[i-1]_-+k-j,\pi_{{\rm A}^{(2)}}(i+k-1)}.
\]
\end{enumerate}
\end{prop}

\begin{prop}\label{D2closed}\cite{Ka}
We assume $\mathfrak{g}$ is of type ${\rm C}^{(1)}_{n-1}$.
Let $T=(y_t)_{t\in\mathbb{Z}}$ be a sequence in ${\rm REYD}^{{\rm D}^{(2)}}_{k}$ (Definition \ref{DEYD}), $s\in\mathbb{Z}_{\geq0}$ and $i\in\mathbb{Z}$.
\begin{enumerate}
\item We assume that the point $(i,y_i)$ is single or double admissible and a sequence
$T'=(y_t')_{t\in\mathbb{Z}}\in{\rm REYD}^{{\rm D}^{(2)}}_{k}$ satisfies
$y_{i}'=y_i-1$ and $y_t'=y_t$ $(t\neq i)$. Putting $j:=y_i$, if $s+P_k^{{\rm D}^{(2)}}(i+k)+[i]_-+k-j\geq1$ then
\[
L^{{\rm D}^{(2)}}_{s,k,\iota}(T')-L^{{\rm D}^{(2)}}_{s,k,\iota}(T)=-\beta_{s+P_k^{{\rm D}^{(2)}}(i+k)+[i]_-+k-j,\pi_{{\rm D}^{(2)}}(i+k)}.
\]
\item We assume that the point $(i,y_{i-1})$ is single or double removable and $T''=(y_t'')_{t\in\mathbb{Z}}\in{\rm REYD}^{{\rm D}^{(2)}}_{k}$
satisfies $y_{i-1}''=y_{i-1}+1$ and $y_t''=y_t$ $(t\neq i-1)$. Putting $j:=y_{i-1}''=y_{i-1}+1$,
if $s+P_k^{{\rm D}^{(2)}}(i+k-1)+[i-1]_-+k-j\geq1$ then
\[
L^{{\rm D}^{(2)}}_{s,k,\iota}(T'')-L^{{\rm D}^{(2)}}_{s,k,\iota}(T)=\beta_{s+P_k^{{\rm D}^{(2)}}(i+k-1)+[i-1]_-+k-j,\pi_{{\rm D}^{(2)}}(i+k-1)}.
\]
\end{enumerate}
\end{prop}

\subsection{On proper Young walls}

\begin{prop}\label{prop-closednessAw-YW}\cite{Ka}
We assume $\mathfrak{g}$ is of type ${\rm A}^{(2)}_{2n-2}$.
\begin{enumerate}
\item
For $t\in I\setminus\{1\}$,
we assume that a wall $Y\in {\rm YW}^{{\rm A}^{(2)}}_{1}$ has a $t$-admissible slot
\[
\begin{xy}
(-8,15) *{(-i-1,\ell+1)}="0000",
(17,15) *{(-i,\ell+1)}="000",
(15,-3) *{(-i,\ell)}="00",
(-7,-3) *{(-i-1,\ell)}="0",
(0,0) *{}="1",
(12,0)*{}="2",
(12,12)*{}="3",
(0,12)*{}="4",
\ar@{--} "1";"2"^{}
\ar@{--} "1";"4"^{}
\ar@{--} "2";"3"^{}
\ar@{--} "3";"4"^{}
\end{xy}
\]
Let $Y'$ be a wall obtained from $Y$
by adding the $t$-block to the slot.
If $Y'\in{\rm YW}^{{\rm A}^{(2)}}_{1}$
then for $s\in\mathbb{Z}_{\geq0}$
such that $s+P_1^{{\rm A}^{(2)}}(\ell)+i\geq1$,
 it follows
\[
L^{{\rm A}^{(2)}}_{s,1,\iota}(Y')=L^{{\rm A}^{(2)}}_{s,1,\iota}(Y)-\beta_{s+P_1^{{\rm A}^{(2)}}(\ell)+i,t}.
\]
\item We assume that a wall $Y\in {\rm YW}^{{\rm A}^{(2)}}_{1}$ has a $1$-admissible slot
\[
\begin{xy}
(-8,9) *{(-i-1,\ell+\frac{1}{2})}="0000",
(17,9) *{(-i,\ell+\frac{1}{2})}="000",
(15,-3) *{(-i,\ell)}="00",
(-7,-3) *{(-i-1,\ell)}="0",
(0,0) *{}="1",
(12,0)*{}="2",
(12,6)*{}="3",
(0,6)*{}="4",
(38,3)*{{\rm or}}="or",
(60,9) *{(-i-1,\ell+1)}="a0000",
(85,9) *{(-i,\ell+1)}="a000",
(83,-3) *{(-i,\ell+\frac{1}{2})}="a00",
(61,-3) *{(-i-1,\ell+\frac{1}{2})}="a0",
(68,0) *{}="a1",
(80,0)*{}="a2",
(80,6)*{}="a3",
(68,6)*{}="a4",
\ar@{--} "1";"2"^{}
\ar@{--} "1";"4"^{}
\ar@{--} "2";"3"^{}
\ar@{--} "3";"4"^{}
\ar@{--} "a1";"a2"^{}
\ar@{--} "a1";"a4"^{}
\ar@{--} "a2";"a3"^{}
\ar@{--} "a3";"a4"^{}
\end{xy}
\]
Let $Y'$ be a wall obtained from $Y$
by adding the $1$-block to the slot.
If $Y'\in{\rm YW}^{{\rm A}^{(2)}}_{1}$
then for $s\in\mathbb{Z}_{\geq0}$
such that $s+P_1^{{\rm A}^{(2)}}(\ell)+i\geq1$,
it follows
\[
L^{{\rm A}^{(2)}}_{s,1,\iota}(Y')=L^{{\rm A}^{(2)}}_{s,1,\iota}(Y)-\beta_{s+P_1^{{\rm A}^{(2)}}(\ell)+i,1}.
\]
\end{enumerate}
\end{prop}

\begin{prop}\label{prop-closednessCw-YW}\cite{Ka}
We assume $\mathfrak{g}$ is of type ${\rm C}^{(1)}_{n-1}$ and $k\in\{1,n\}$.
\begin{enumerate}
\item
Let $t\in I\setminus\{1,n\}$ and
we assume that a wall $Y\in {\rm YW}^{{\rm D}^{(2)}}_{k}$ has a $t$-admissible slot
\[
\begin{xy}
(-8,15) *{(-i-1,\ell+1)}="0000",
(17,15) *{(-i,\ell+1)}="000",
(15,-3) *{(-i,\ell)}="00",
(-7,-3) *{(-i-1,\ell)}="0",
(0,0) *{}="1",
(12,0)*{}="2",
(12,12)*{}="3",
(0,12)*{}="4",
\ar@{--} "1";"2"^{}
\ar@{--} "1";"4"^{}
\ar@{--} "2";"3"^{}
\ar@{--} "3";"4"^{}
\end{xy}
\]
Let $Y'$ be a wall obtained from $Y$
by adding the $t$-block to the slot.
If $Y'\in{\rm YW}^{{\rm D}^{(2)}}_{k}$
then for $s\in\mathbb{Z}_{\geq0}$
such that $s+P_k^{{\rm D}^{(2)}}(\ell)+i\geq1$,
it follows
\[
L^{{\rm D}^{(2)}}_{s,k,\iota}(Y')=L^{{\rm D}^{(2)}}_{s,k,\iota}(Y)-\beta_{s+P_k^{{\rm D}^{(2)}}(\ell)+i,t}.
\]
\item Let $t$ be $t=1$ or $t=n$. We assume that a wall $Y\in {\rm YW}^{{\rm D}^{(2)}}_{k}$ has a $t$-admissible slot
\[
\begin{xy}
(-8,9) *{(-i-1,\ell+\frac{1}{2})}="0000",
(17,9) *{(-i,\ell+\frac{1}{2})}="000",
(15,-3) *{(-i,\ell)}="00",
(-7,-3) *{(-i-1,\ell)}="0",
(0,0) *{}="1",
(12,0)*{}="2",
(12,6)*{}="3",
(0,6)*{}="4",
(38,3)*{{\rm or}}="or",
(60,9) *{(-i-1,\ell+1)}="a0000",
(85,9) *{(-i,\ell+1)}="a000",
(83,-3) *{(-i,\ell+\frac{1}{2})}="a00",
(61,-3) *{(-i-1,\ell+\frac{1}{2})}="a0",
(68,0) *{}="a1",
(80,0)*{}="a2",
(80,6)*{}="a3",
(68,6)*{}="a4",
\ar@{--} "1";"2"^{}
\ar@{--} "1";"4"^{}
\ar@{--} "2";"3"^{}
\ar@{--} "3";"4"^{}
\ar@{--} "a1";"a2"^{}
\ar@{--} "a1";"a4"^{}
\ar@{--} "a2";"a3"^{}
\ar@{--} "a3";"a4"^{}
\end{xy}
\]
Let $Y'$ be a wall obtained from $Y$
by adding the $t$-block to the slot.
If $Y'\in{\rm YW}^{{\rm D}^{(2)}}_{k}$
then for $s\in\mathbb{Z}_{\geq0}$ such that $s+P_k^{{\rm D}^{(2)}}(\ell)+i\geq1$,
it follows
\[
L^{{\rm D}^{(2)}}_{s,k,\iota}(Y')=L^{{\rm D}^{(2)}}_{s,k,\iota}(Y)-\beta_{s+P_k^{{\rm D}^{(2)}}(\ell)+i,t}.
\]
\end{enumerate}
\end{prop}

\section{Proofs}

In this section, we complete proofs of theorems in Section 4.

\subsection{Type ${\rm A}_{n-1}^{(1)}$-case and ${\rm D}_{n}^{(2)}$-case}

In this subsection, let $\mathfrak{g}$ be of type $X={\rm A}_{n-1}^{(1)}$ or 
$X={\rm D}_{n}^{(2)}$ and prove Theorem \ref{thmA1}. 
First, we consider the first set
$\{\what{S}'_{j_{\ell}}\cdots \what{S}'_{j_1}x_{j_0} | \ell\in\mathbb{Z}_{\geq0}, j_0,\cdots,j_{\ell}\in \mathbb{Z}_{\geq1} \}$
of right hand side of (\ref{xilamdef}).
By Theorem 4.3, 4.5 of \cite{Ka}, the sequence $\iota$ satisfies the $\Xi'$-positivity condition
and it follows from (\ref{twoS2}) and the proofs of Theorem 4.3, 4.5 of \cite{Ka} that
\begin{eqnarray}
\{\what{S}'_{j_{\ell}}\cdots \what{S}'_{j_1}x_{j_0} | \ell\in\mathbb{Z}_{\geq0}, j_0,\cdots,j_{\ell}\in \mathbb{Z}_{\geq1} \}
&=&\{ S'_{j_{\ell}}\cdots S'_{j_1}x_{j_0} | \ell\in\mathbb{Z}_{\geq0}, j_0,\cdots,j_{\ell}\in \mathbb{Z}_{\geq1} \}\nonumber\\
&=&\{ S'_{j_{\ell}}\cdots S'_{j_1}x_{s,k} | \ell\in\mathbb{Z}_{\geq0}, j_1,\cdots,j_{\ell}\in \mathbb{Z}_{\geq1},\ k\in I,\ s\in\mathbb{Z}_{\geq1} \}
\nonumber\\
&=& \{ L^{X^L}_{s,k,\iota}(T)| k\in I,\ s\in\mathbb{Z}_{\geq1}, T\in {\rm EYD}_k \}\label{half1}\\
&=&{\rm Comb}_{\iota}^{X^L}[\infty].\nonumber
\end{eqnarray}
Next, let us consider the second set
$\{\what{S}'_{j_{\ell}}\cdots \what{S}'_{j_1} \lambda^{(k)} | \ell\in\mathbb{Z}_{\geq0},\ k\in I,\ j_1,\cdots,j_{\ell}\in \mathbb{Z}_{\geq1} \}$
of (\ref{xilamdef}). We fix an index $k\in I$. In this proof, we simply write $P_k^{X^L}$ and $\pi_{X^L}$ as
$P_k$ and $\pi$.

\vspace{2mm}

\nd
\underline{Case 1 : it holds $(1,k)<(1,\pi(k-1)), (1,\pi(k+1))$}

\vspace{2mm}

Since $\lambda^{(k)}=-x_{1,k}+\langle h_k,\lambda \rangle$, it is clear that
$\{\what{S}'_{j_{\ell}}\cdots \what{S}'_{j_1} \lambda^{(k)} | \ell\in\mathbb{Z}_{\geq0},\ j_1,\cdots,j_{\ell}\in \mathbb{Z}_{\geq1} \}
=\{0,\lambda^{(k)}\}$.

\vspace{2mm}

\nd
\underline{Case 2 : it holds $(1,k)>(1,\pi(k-1))$ and $(1,k)<(1,\pi(k+1))$}

\vspace{2mm}

In this case, the definition of $\lambda^{(k)}$ implies
$\lambda^{(k)}=-x_{1,k}+x_{1,\pi(k-1)}+\langle h_k,\lambda \rangle$.
Using the notation in (\ref{r-tilde-A1}), one can express $\lambda^{(k)}$ as $\lambda^{(k)}=\fbox{$\tilde{k}$}^{{\rm X}^{L}}_k
+\langle h_k,\lambda \rangle$.
Following the definition of $\what{S}'$,
we get the following diagram of actions of $\what{S}'$:
\[
\begin{xy}
(-8,0) *{-\langle h_k,\lambda \rangle}="00",
(0,0) *{\ }="0",
(10,0) *{\fbox{$\tilde{k}$}^{X^L}}="r-1r-1",
(40,0)*{\fbox{$\overset{\sim}{k-1}$}^{X^L}}="r-2r-2",
(70,0)*{\fbox{$\overset{\sim}{k-2}$}^{X^L}}="r-3r-3",
(100,0)*{\fbox{$\overset{\sim}{k-3}$}^{X^L}}="r-4r-4",
(130,0)*{\cdots}="dot",
\ar@{->} "r-1r-1";"00"^{\what{S}'_{1,k}}
\ar@/^/ @{->} "r-1r-1";"r-2r-2"^{\what{S}'_{P_k(k-1),\pi(k-1)}}
\ar@/^/ @{->} "r-2r-2";"r-3r-3"^{\what{S}'_{P_k(k-2),\pi(k-2)}}
\ar@/^/ @{->} "r-3r-3";"r-4r-4"^{\what{S}'_{P_k(k-3),\pi(k-3)}}
\ar@/^/ @{->} "r-4r-4";"dot"^{\what{S}'_{P_k(k-4),\pi(k-4)}}
\ar@/_/ @{<-} "r-1r-1";"r-2r-2"_{\what{S}'_{1+P_k(k-1),\pi(k-1)}}
\ar@/_/ @{<-} "r-2r-2";"r-3r-3"_{\what{S}'_{1+P_k(k-2),\pi(k-2)}}
\ar@/_/ @{<-} "r-3r-3";"r-4r-4"_{\what{S}'_{1+P_k(k-3),\pi(k-3)}}
\ar@/_/ @{<-} "r-4r-4";"dot"_{\what{S}'_{1+P_k(k-4),\pi(k-4)}}
\end{xy}
\]
Here, we omitted index $k$ of boxes.
Other actions of $\what{S}'$ are trivial. 
Note that by $p_{\pi(k-1),\pi(k)}=1$, if $r<k$ then
\begin{equation}\label{stpspr1}
1+P_k(r)\geq2,
\end{equation}
which is the left index of $-x_{1+P_k(r),\pi(r)}$ in (\ref{r-tilde-A1}).
Thus,
\[
\{\what{S}'_{j_{\ell}}\cdots \what{S}'_{j_1} \lambda^{(k)} | \ell\in\mathbb{Z}_{\geq0},\ j_1,\cdots,j_{\ell}\in \mathbb{Z}_{\geq1} \}
=\{0,\fbox{$\tilde{r}$}^{X^L}_k+\langle h_k,\lambda \rangle\ | r\in\mathbb{Z}_{\leq k}\}.
\]

\vspace{2mm}

\nd
\underline{Case 3 : it holds $(1,k)<(1,\pi(k-1))$ and $(1,k)>(1,\pi(k+1))$}

\vspace{2mm}

In this case, we have
$\lambda^{(k)}=-x_{1,k}+x_{1,\pi(k+1)}+\langle h_k,\lambda \rangle$.
Using the notation in (\ref{r-tilde-A0}), 
it holds 
$\lambda^{(k)}=\fbox{$k+1$}^{X^L}_k+\langle h_k,\lambda \rangle$.
One obtains the following diagram of actions of $\what{S}'$:
\[
\begin{xy}
(-8,0) *{-\langle h_k,\lambda \rangle}="00",
(0,0) *{\ }="0",
(10,0) *{\fbox{$k+1$}^{X^L}}="r-1r-1",
(40,0)*{\fbox{$k+2$}^{X^L}}="r-2r-2",
(70,0)*{\fbox{$k+3$}^{X^L}}="r-3r-3",
(100,0)*{\fbox{$k+4$}^{X^L}}="r-4r-4",
(130,0)*{\cdots}="dot",
\ar@{->} "r-1r-1";"00"^{\what{S}'_{1,k}}
\ar@/^/ @{->} "r-1r-1";"r-2r-2"^{\what{S}'_{P_k(k+1),\pi(k+1)}}
\ar@/^/ @{->} "r-2r-2";"r-3r-3"^{\what{S}'_{P_k(k+2),\pi(k+2)}}
\ar@/^/ @{->} "r-3r-3";"r-4r-4"^{\what{S}'_{P_k(k+3),\pi(k+3)}}
\ar@/^/ @{->} "r-4r-4";"dot"^{\what{S}'_{P_k(k+4),\pi(k+4)}}
\ar@/_/ @{<-} "r-1r-1";"r-2r-2"_{\what{S}'_{1+P_k(k+1),\pi(k+1)}}
\ar@/_/ @{<-} "r-2r-2";"r-3r-3"_{\what{S}'_{1+P_k(k+2),\pi(k+2)}}
\ar@/_/ @{<-} "r-3r-3";"r-4r-4"_{\what{S}'_{1+P_k(k+3),\pi(k+3)}}
\ar@/_/ @{<-} "r-4r-4";"dot"_{\what{S}'_{1+P_k(k+4),\pi(k+4)}}
\end{xy}
\]
Here, we omitted index $k$.
Other actions of $\what{S}'$ are trivial. Note that for $r>k+1$, it holds
\begin{equation}\label{stpspr2}
1+P_k(r-1)\geq2,
\end{equation}
which is the left index of $-x_{1+P_k(r-1),\pi(r-1)}$ in (\ref{r-tilde-A0}).
Thus,
\[
\{\what{S}'_{j_{\ell}}\cdots \what{S}'_{j_1} \lambda^{(k)} | \ell\in\mathbb{Z}_{\geq0},\ j_1,\cdots,j_{\ell}\in \mathbb{Z}_{\geq1} \}
=\{0,\fbox{$r$}^{X^L}_k+\langle h_k,\lambda \rangle\ | r\in\mathbb{Z}_{\geq k+1}\}.
\]

\vspace{2mm}

\nd
\underline{Case 4 : it holds $(1,k)>(1,\pi(k-1))$ and $(1,k)>(1,\pi(k+1))$}

\vspace{2mm}

In this case, it holds $\lambda^{(k)}=-x_{1,k}+x_{1,\pi(k-1)}+x_{1,\pi(k+1)}+
\lan h_k,\lambda \ran$. Note that if $(i,j)\in\mathbb{Z}_{\geq0}\times\mathbb{Z}_{\leq k}\setminus\{(0,k)\}$
then 
\begin{equation}\label{geq-1}
P_k(i+j)+{\rm min}\{k-j,i\}\geq1.
\end{equation}
First, let us prove
\begin{equation}\label{A-pr-1}
\{\what{S}'_{j_{\ell}}\cdots \what{S}'_{j_1} \lambda^{(k)} | \ell\in\mathbb{Z}_{\geq0},\ j_1,\cdots,j_{\ell}\in \mathbb{Z}_{\geq1} \}
\subset
\{L^{X^L}_{0,k,\iota}(T)+\langle h_k,\lambda \rangle | T\in{\rm EYD}_{k}\setminus\{\phi^k\}\}\cup\{0\}.
\end{equation}
We put the right hand side of (\ref{A-pr-1}) as $Z$: 
$Z:=\{L^{X^L}_{0,k,\iota}(T)+\langle h_k,\lambda \rangle | T\in{\rm EYD}_{k}\setminus\{\phi^k\}\}\cup\{0\}$.
We take $T^1_k\in{\rm EYD}_{k}\setminus\{\phi^k\}$ as follows:
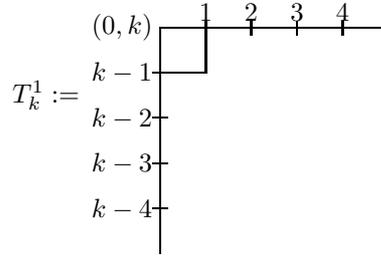
\[
\begin{xy}
(40,-9) *{T^1_k:=}="Y1",
(55,0) *{}="1a",
(85,0)*{}="2a",
(55,-30)*{}="3a",
(50,0)*{(0,k)}="4a",
(61,2) *{1}="10a",
(61,-1) *{}="1010a",
(61,-6) *{}="1-0a",
(67,2) *{2}="11a",
(67,-1) *{}="1111a",
(73,2) *{3}="12a",
(73,-1) *{}="1212a",
(79,2) *{4}="13a",
(79,-1) *{}="1313a",
(50,-6)*{k-1}="5a",
(56,-6)*{}="55a",
(50,-12)*{k-2}="6a",
(56,-12)*{}="66a",
(50,-18)*{k-3}="7a",
(56,-18)*{}="77a",
(50,-24)*{k-4}="8a",
(56,-24)*{}="88a",
\ar@{-} "1a";"2a"^{}
\ar@{-} "1a";"3a"^{}
\ar@{-} "5a";"55a"^{}
\ar@{-} "6a";"66a"^{}
\ar@{-} "7a";"77a"^{}
\ar@{-} "8a";"88a"^{}
\ar@{-} "10a";"1010a"^{}
\ar@{-} "11a";"1111a"^{}
\ar@{-} "12a";"1212a"^{}
\ar@{-} "13a";"1313a"^{}
\ar@{-} "55a";"1-0a"^{}
\ar@{-} "1010a";"1-0a"^{}
\end{xy}
\]
The points $(1,k)$ and $(0,k-1)$ are concave corners and
$(1,k-1)$ is the unique convex corners in $T^1_k$.
Considering $P_k(k+1)=p_{\pi(k+1),\pi(k)}=1$, $P_k(k-1)=p_{\pi(k-1),\pi(k)}=1$,
it follows
\[
L^{X^L}_{0,k,\iota}(T^1_k)+\lan h_k,\lambda \ran=
-x_{1,k}+x_{1,\pi(k-1)}+x_{1,\pi(k+1)}+
\lan h_k,\lambda \ran=\lambda^{(k)}.
\]
Hence, $\lambda^{(k)}\in \{L^{X^L}_{0,k,\iota}(T)+\langle h_k,\lambda \rangle | T\in{\rm EYD}_{k}\setminus\{\phi^k\}\}$.
Thus, to prove (\ref{A-pr-1}), we need to show
$Z$
is closed under the action of $\what{S}'_r$ for all $r\in\mathbb{Z}_{\geq1}$.
For any $T=(y_r)_{r\in\mathbb{Z}_{\geq1}}\in {\rm EYD}_{k}\setminus\{\phi^k\}$, if the coefficient of $x_r$ is positive in 
$L^{X^L}_{0,k,\iota}(T)+\langle h_k,\lambda \rangle$ then
$T$ has a concave corner $(i,j)$ such that $j=y_i$ and
\[
r=(P_k(i+j)+{\rm min}\{k-j,i\},\pi(i+j))
\]
by the definition of $L^{X^L}_{0,k,\iota}$ in (\ref{ovlL1}), (\ref{ovlL2}). Note that
we identify $\mathbb{Z}_{\geq1}$ with $\mathbb{Z}_{\geq1}\times I$ by the method of subsection \ref{seno}.
Following the definition of $\what{S}'_r$,
one obtains
\[
\what{S}'_r(L^{X^L}_{0,k,\iota}(T)+\langle h_k,\lambda \rangle)
=L^{X^L}_{0,k,\iota}(T)+\langle h_k,\lambda \rangle - \beta_{P_k(i+j)+{\rm min}\{k-j,i\},\pi(i+j)}.
\]
Taking $T'=(y_r')_{r\in\mathbb{Z}_{\geq1}}\in {\rm EYD}_{k}\setminus\{\phi^k\}$ such that
$y_i'=y_i-1$, $y_r'=y_r$ ($r\neq i$), it follows by Proposition \ref{prop-closednessAD}
that
\begin{eqnarray*}
\what{S}'_r(L^{X^L}_{0,k,\iota}(T)+\langle h_k,\lambda \rangle)
&=&L^{X^L}_{0,k,\iota}(T)- \beta_{P_k(i+j)+{\rm min}\{k-j,i\},\pi(i+j)}+\langle h_k,\lambda \rangle\\
&=&
L^{X^L}_{0,k,\iota}(T')+\langle h_k,\lambda \rangle\in
Z.
\end{eqnarray*}
If the coefficient of $x_r$ is $0$ then 
it is clear $\what{S}'_r(L^{X^L}_{0,k,\iota}(T)+\langle h_k,\lambda \rangle)=L^{X^L}_{0,k,\iota}(T)+\langle h_k,\lambda \rangle
\in
Z$.
We assume 
the coefficient of $x_r$ is negative in 
$L^{X^L}_{0,k,\iota}(T)+\langle h_k,\lambda \rangle$.
Then
$T$ has a convex corner $(i,j)$ such that $j=y_{i-1}$ and
\[
r=(P_k(i+j)+{\rm min}\{k-j,i\},\pi(i+j)).
\]
Note that if $(i,j)=(1,k-1)$ then $T=T^1_k$ and $\what{S}'_r(L^{X^L}_{0,k,\iota}(T)+\langle h_k,\lambda \rangle)=0$.
Thus, we may assume $(i,j)\neq(1,k-1)$.
Since $(i,j)$ is convex, we see that ${\rm min}\{k-j,i\}\geq1$.
In case of $i+j\neq k$, one can verify $P_k(i+j)\geq1$ by 
$p_{\pi(k+1),\pi(k)}=p_{\pi(k-1),\pi(k)}=1$ so that $P_k(i+j)+{\rm min}\{k-j,i\}\geq2$.
In case of $i+j=k$, 
since we assumed $(i,j)\neq(1,k-1)$,
it follows $i=k-j\geq2$ and $P_k(i+j)+{\rm min}\{k-j,i\}={\rm min}\{k-j,i\}\geq2$.
Hence, 
it holds
\begin{equation}\label{A-pr-1a}
P_k(i+j)+{\rm min}\{k-j,i\}\geq2
\end{equation}
and
taking $T''=(y_r'')_{r\in\mathbb{Z}_{\geq1}}$ as $y_{i-1}''=y_{i-1}+1$, $y_r''=y_r$ ($r\neq i-1$), 
the point $(i-1,y_{i-1}'')=(i-1,j+1)$ is a concave corner in
$T''$ and
it follows
from Proposition \ref{prop-closednessAD} that
\begin{eqnarray}
\what{S}'_r(L^{X^L}_{0,k,\iota}(T)+\langle h_k,\lambda \rangle)
&=&
L^{X^L}_{0,k,\iota}(T)+\langle h_k,\lambda \rangle+
\beta_{P_k(i+j)+{\rm min}\{k-j,i\}-1,\pi(i+j)}\nonumber\\
&=&
L^{X^L}_{0,k,\iota}(T)+
\beta_{P_k((i-1)+(j+1))+{\rm min}\{k-(j+1),i-1\},\pi((i-1)+(j+1))}+\langle h_k,\lambda \rangle \nonumber\\
&=&
L^{X^L}_{0,k,\iota}(T'')+\langle h_k,\lambda \rangle
\in
Z. \label{A-pr1-b}
\end{eqnarray}
Thus, we obtain (\ref{A-pr-1}).

Next, let us show
\begin{equation}\label{A-pr-2}
\{\what{S}'_{j_{\ell}}\cdots \what{S}'_{j_1} \lambda^{(k)} | \ell\in\mathbb{Z}_{\geq0},\ j_1,\cdots,j_{\ell}\in \mathbb{Z}_{\geq1} \}
\supset Z.
\end{equation}
Since we know $0\in \{\what{S}'_{j_{\ell}}\cdots \what{S}'_{j_1} \lambda^{(k)} | \ell\in\mathbb{Z}_{\geq0},\ j_1,\cdots,j_{\ell}\in \mathbb{Z}_{\geq1} \}$, let us prove
for any $T\in{\rm EYD}_{k}\setminus\{\phi^k\}$,
it holds 
\begin{equation}\label{A-pr-3}
L^{X^L}_{0,k,\iota}(T)+\langle h_k,\lambda \rangle\in
\{\what{S}'_{j_{\ell}}\cdots \what{S}'_{j_1} \lambda^{(k)} | \ell\in\mathbb{Z}_{\geq0},\ j_1,\cdots,j_{\ell}\in \mathbb{Z}_{\geq1} \}.
\end{equation}
We identify each $T\in {\rm EYD}_k\setminus\{\phi^k\}$ as a Young diagram consisting of several boxes, where each box is the square
whose length of sides are $1$. Let us show (\ref{A-pr-3}) by induction on the number of boxes in $T$. 
If $T$ has only one box then $T=T^1_k$ and
$L^{X^L}_{0,k,\iota}(T)+\langle h_k,\lambda \rangle=-x_{1,k}+x_{1,\pi(k-1)}+x_{1,\pi(k+1)}+
\lan h_k,\lambda \ran=\lambda^{(k)}\in
\{\what{S}'_{j_{\ell}}\cdots \what{S}'_{j_1} \lambda^{(k)} | \ell\in\mathbb{Z}_{\geq0},\ j_1,\cdots,j_{\ell}\in \mathbb{Z}_{\geq1} \}$.
Thus,
we may assume that $T$ has at least two boxes.
By Proposition \ref{prop-closednessAD}, (\ref{geq-1}) and
$L^{X^L}_{0,k,\iota}(T^k_1)+\langle h_k,\lambda \rangle=\lambda^{(k)}$, we see that
\[
L^{X^L}_{0,k,\iota}(T)+\lan h_k,\lambda \ran
=
\lambda^{(k)} - \sum_{(t,d)\in\mathbb{Z}_{\geq1}\times I} c_{t,d}\beta_{t,d}
\]
with some non-negative integers $\{c_{t,d}\}$ such that $c_{t,d}=0$ except for finitely many $(t,d)$.
Since $T\neq T^1_k$, 
there exists $(t',d')\in\mathbb{Z}_{\geq1}\times I$ such that
\[
(t',d')={\rm max}\{(t,d)\in\mathbb{Z}_{\geq1}\times I | c_{t,d}\neq0 \}.
\]
The definition (\ref{betak}) of $\beta_{t,d}$ and Definition \ref{adapt} mean
the coefficient of $x_{t'+1,d'}$ is negative in $L^{X^L}_{0,k,\iota}(T)$. We define $(t'',d'')$ as
\[
(t'',d'')={\rm min}\{(t,d)\in\mathbb{Z}_{\geq1}\times I | \text{the coefficient of }x_{t,d}\ \text{in } L^{X^L}_{0,k,\iota}(T)\text{ is negative} \}.
\]
By the definition of $L^{X^L}_{0,k,\iota}(T)$,
$T$ has a convex corner $(i+1,j-1)$ such that $i\in\mathbb{Z}_{\geq0}$, $j\in\mathbb{Z}_{\leq k}$ and $L^{X^L}_{0,k,\iota}(i+1,j-1)=x_{t'',d''}$. 
Just as in (\ref{A-pr-1a}), one can prove
$t''-1\geq1$.
Let $T'\in {\rm EYD}_k\setminus\{\phi^k\}$ be
the extended Young diagram obtained from $T$ by replacing the
convex corner $(i+1,j-1)$ by a concave corner. 
Just as in (\ref{A-pr1-b}), we obtain
\[
L^{X^L}_{0,k,\iota}(T)+\beta_{t''-1,d''}=L^{X^L}_{0,k,\iota}(T').
\]
The minimality of $(t'',d'')$ means the coefficient of $x_{t''-1,d''}$ in 
$L^{X^L}_{0,k,\iota}(T')$ is positive and
\begin{equation}\label{A-pr-4}
\what{S}'_{t''-1,d''}(L^{X^L}_{0,k,\iota}(T')+\lan h_k,\lambda \ran)=L^{X^L}_{0,k,\iota}(T')+\lan h_k,\lambda \ran-\beta_{t''-1,d''}=L^{X^L}_{0,k,\iota}(T)+\lan h_k,\lambda \ran.
\end{equation}
Note that the number of boxes in $T'$ is smaller than $T$. Hence,
the induction assumption implies
\[
L^{X^L}_{0,k,\iota}(T')+\lan h_k,\lambda \ran=
\what{S}'_{j_{\ell}}\cdots \what{S}'_{j_1} \lambda^{(k)}
\]
with some $j_1,\cdots,j_{\ell}\in\mathbb{Z}_{\geq1}$. Combining with
(\ref{A-pr-4}), one obtains 
$L^{X^L}_{0,k,\iota}(T)+\lan h_k,\lambda \ran\in \{\what{S}'_{j_{\ell}}\cdots \what{S}'_{j_1} \lambda^{(k)} | \ell\in\mathbb{Z}_{\geq0},\ j_1,\cdots,j_{\ell}\in \mathbb{Z}_{\geq1} \}$, which implies (\ref{A-pr-3}) so that $(\ref{A-pr-2})$.

By the above argument, it follows
\begin{equation}\label{stpspr3}
\{\what{S}'_{j_{\ell}}\cdots \what{S}'_{j_1} \lambda^{(k)} | \ell\in\mathbb{Z}_{\geq0},\ j_1,\cdots,j_{\ell}\in \mathbb{Z}_{\geq1} \}
={\rm Comb}^{X^L}_k[\lambda]\cup\{0\}.
\end{equation}
Combining with (\ref{half1}), one can confirm the $\Xi'$-ample condition $\textbf{0}=(\cdots,0,0,0)\in\Sigma'_{\iota}[\lambda]$
(Definition \ref{ample1}). Using Theorem \ref{Nthm1}, we get our claim. \qed

\subsection{Type ${\rm A}_{2n-2}^{(2)}$-case and ${\rm C}_{n-1}^{(1)}$-case}

In this subsection, let $\mathfrak{g}$ be of type $X={\rm A}_{2n-2}^{(2)}$ or 
$X={\rm C}_{n-1}^{(1)}$ and prove Theorem \ref{thmA2}. 
First, we consider the first set
$\{\what{S}'_{j_{\ell}}\cdots \what{S}'_{j_1}x_{j_0} | \ell\in\mathbb{Z}_{\geq0}, j_0,\cdots,j_{\ell}\in \mathbb{Z}_{\geq1} \}$
of right hand side of (\ref{xilamdef}).
By Theorem 4.8, 4.11 of \cite{Ka}, the sequence $\iota$ satisfies the $\Xi'$-positivity condition
and it follows from (\ref{twoS2}) and the proofs of Theorem 4.8, 4.11 of \cite{Ka} that
\begin{eqnarray}
\{\what{S}'_{j_{\ell}}\cdots \what{S}'_{j_1}x_{j_0} | \ell\in\mathbb{Z}_{\geq0}, j_0,\cdots,j_{\ell}\in \mathbb{Z}_{\geq1} \}
&=&\{ S'_{j_{\ell}}\cdots S'_{j_1}x_{j_0} | \ell\in\mathbb{Z}_{\geq0}, j_0,\cdots,j_{\ell}\in \mathbb{Z}_{\geq1} \}\nonumber\\
&=&\{ S'_{j_{\ell}}\cdots S'_{j_1}x_{s,k} | \ell\in\mathbb{Z}_{\geq0}, j_1,\cdots,j_{\ell}\in \mathbb{Z}_{\geq1},\ k\in I,\ s\in\mathbb{Z}_{\geq1} \}
\nonumber\\
&=& \{ L^{X^L}_{s,k,\iota}(T)| k\in I_{X^L},\ s\in\mathbb{Z}_{\geq1}, T\in {\rm REYD}^{X^L}_{k} \}\nonumber\\
& &\cup
\{ L^{X^L}_{s,k,\iota}(T)| k\in I\setminus I_{X^L},\ s\in\mathbb{Z}_{\geq1}, T\in {\rm YW}^{X^L}_{k} \}\label{half2}\\
&=&{\rm Comb}_{\iota}^{X^L}[\infty].\nonumber
\end{eqnarray}
Next, we consider the second set
$\{\what{S}'_{j_{\ell}}\cdots \what{S}'_{j_1} \lambda^{(k)} | \ell\in\mathbb{Z}_{\geq0},\ k\in I,\ j_1,\cdots,j_{\ell}\in \mathbb{Z}_{\geq1} \}$
of (\ref{xilamdef}). We fix an index $k\in I$. In this proof, we simply write $P_k^{X^L}$ and $\pi_{X^L}$ as
$P_k$ and $\pi$. We take $k$ as $k\in I_{X^L}$ and
for $j=k-1$, $k+1$,
\[
c_{j}:=
\begin{cases}
2 & \text{if }\pi(j)\in I\setminus I_X,\\ 
1 & \text{otherwise}.
\end{cases}
\]

\vspace{2mm}

\nd
\underline{Case 1 : it holds $(1,k)<(1,\pi(k-1)), (1,\pi(k+1))$}

\vspace{2mm}

By $\lambda^{(k)}=-x_{1,k}+\langle h_k,\lambda \rangle$, it is easy to see
$\{\what{S}'_{j_{\ell}}\cdots \what{S}'_{j_1} \lambda^{(k)} | \ell\in\mathbb{Z}_{\geq0},\ j_1,\cdots,j_{\ell}\in \mathbb{Z}_{\geq1} \}
=\{0,\lambda^{(k)}\}$.

\vspace{2mm}

\nd
\underline{Case 2 : it holds $(1,k)>(1,\pi(k-1))$ and $(1,k)<(1,\pi(k+1))$}

\vspace{2mm}

The definition of $\lambda^{(k)}$ implies
$\lambda^{(k)}=-x_{1,k}+c_{k-1}x_{1,\pi(k-1)}+\langle h_k,\lambda \rangle$.
Using the notation in (\ref{r-tilde-D1}), it holds $\lambda^{(k)}=\fbox{$\tilde{k}$}^{{\rm X}^{L}}_k+\langle h_k,\lambda \rangle$.
By the definition of $\what{S}'$,
the actions of $\what{S}'$ are as follows: 
\[
\begin{xy}
(-8,0) *{-\langle h_k,\lambda \rangle}="00",
(0,0) *{\ }="0",
(10,0) *{\fbox{$\tilde{k}$}^{X^L}}="r-1r-1",
(40,0)*{\fbox{$\overset{\sim}{k-1}$}^{X^L}}="r-2r-2",
(70,0)*{\fbox{$\overset{\sim}{k-2}$}^{X^L}}="r-3r-3",
(100,0)*{\fbox{$\overset{\sim}{k-3}$}^{X^L}}="r-4r-4",
(130,0)*{\cdots}="dot",
\ar@{->} "r-1r-1";"00"^{\what{S}'_{1,k}}
\ar@/^/ @{->} "r-1r-1";"r-2r-2"^{\what{S}'_{P_k(k-1),\pi(k-1)}}
\ar@/^/ @{->} "r-2r-2";"r-3r-3"^{\what{S}'_{P_k(k-2),\pi(k-2)}}
\ar@/^/ @{->} "r-3r-3";"r-4r-4"^{\what{S}'_{P_k(k-3),\pi(k-3)}}
\ar@/^/ @{->} "r-4r-4";"dot"^{\what{S}'_{P_k(k-4),\pi(k-4)}}
\ar@/_/ @{<-} "r-1r-1";"r-2r-2"_{\what{S}'_{1+P_k(k-1),\pi(k-1)}}
\ar@/_/ @{<-} "r-2r-2";"r-3r-3"_{\what{S}'_{1+P_k(k-2),\pi(k-2)}}
\ar@/_/ @{<-} "r-3r-3";"r-4r-4"_{\what{S}'_{1+P_k(k-3),\pi(k-3)}}
\ar@/_/ @{<-} "r-4r-4";"dot"_{\what{S}'_{1+P_k(k-4),\pi(k-4)}}
\end{xy}
\]
\[
\begin{xy}
(-5,0) *{}="00",
(0,0) *{\ }="0",
(10,0) *{\cdots}="r-1r-1",
(40,0)*{\fbox{$\tilde{2}$}^{X^L}}="r-2r-2",
(70,0)*{\fbox{$\tilde{1}$}^{X^L}}="r-3r-3",
(100,0)*{\fbox{$\tilde{0}$}^{X^L}}="r-4r-4",
(130,0)*{\cdots}="dot",
\ar@/^/ @{->} "r-1r-1";"r-2r-2"^{\what{S}'_{P_k(2),\pi(2)}}
\ar@/^/ @{->} "r-2r-2";"r-3r-3"^{\what{S}'_{P_k(1),\pi(1)}}
\ar@/^/ @{->} "r-3r-3";"r-4r-4"^{\what{S}'_{P_k(0),\pi(0)}}
\ar@/^/ @{->} "r-4r-4";"dot"^{\what{S}'_{P_k(-1),\pi(-1)}}
\ar@/_/ @{<-} "r-1r-1";"r-2r-2"_{\what{S}'_{1+P_k(2),\pi(2)}}
\ar@/_/ @{<-} "r-2r-2";"r-3r-3"_{\what{S}'_{1+P_k(1),\pi(1)}}
\ar@/_/ @{<-} "r-3r-3";"r-4r-4"_{\what{S}'_{1+P_k(0),\pi(0)}}
\ar@/_/ @{<-} "r-4r-4";"dot"_{\what{S}'_{1+P_k(-1),\pi(-1)}}
\end{xy}
\]
Here, we omitted index $k$.
Other actions of $\what{S}'$ are trivial. 
Note that by $p_{\pi(k-1),\pi(k)}=1$, if $r<k$ then
\begin{equation}\label{stpspr4}
1+P_k(r)\geq2,
\end{equation}
which is the left index of $-x_{1+P_k(r),\pi(r)}$ in (\ref{r-tilde-D1}).
Thus,
\[
\{\what{S}'_{j_{\ell}}\cdots \what{S}'_{j_1} \lambda^{(k)} | \ell\in\mathbb{Z}_{\geq0},\ j_1,\cdots,j_{\ell}\in \mathbb{Z}_{\geq1} \}
=\{0,\fbox{$\tilde{r}$}^{X^L}_k+\langle h_k,\lambda \rangle\ | r\in\mathbb{Z}_{\leq k}\}.
\]

\vspace{2mm}

\nd
\underline{Case 3 : it holds $(1,k)>(1,\pi(k+1))$ and $(1,k)<(1,\pi(k-1))$}

\vspace{2mm}

It is easy to see
$\lambda^{(k)}=-x_{1,k}+c_{k+1}x_{1,\pi(k+1)}+\langle h_k,\lambda \rangle$ so that
$\lambda^{(k)}=\fbox{$k+1$}^{{\rm X}^{L}}_k+\langle h_k,\lambda \rangle$ by using the notation in (\ref{r-tilde-D0}).
By the definition of $\what{S}'$,
the actions of $\what{S}'$ are as follows: 
\[
\begin{xy}
(-8,0) *{-\langle h_k,\lambda \rangle}="00",
(0,0) *{\ }="0",
(10,0) *{\fbox{$k+1$}^{X^L}}="r-1r-1",
(40,0)*{\fbox{$k+2$}^{X^L}}="r-2r-2",
(70,0)*{\fbox{$k+3$}^{X^L}}="r-3r-3",
(100,0)*{\fbox{$k+4$}^{X^L}}="r-4r-4",
(130,0)*{\cdots}="dot",
\ar@{->} "r-1r-1";"00"^{\what{S}'_{1,k}}
\ar@/^/ @{->} "r-1r-1";"r-2r-2"^{\what{S}'_{P_k(k+1),\pi(k+1)}}
\ar@/^/ @{->} "r-2r-2";"r-3r-3"^{\what{S}'_{P_k(k+2),\pi(k+2)}}
\ar@/^/ @{->} "r-3r-3";"r-4r-4"^{\what{S}'_{P_k(k+3),\pi(k+3)}}
\ar@/^/ @{->} "r-4r-4";"dot"^{\what{S}'_{P_k(k+4),\pi(k+4)}}
\ar@/_/ @{<-} "r-1r-1";"r-2r-2"_{\what{S}'_{1+P_k(k+1),\pi(k+1)}}
\ar@/_/ @{<-} "r-2r-2";"r-3r-3"_{\what{S}'_{1+P_k(k+2),\pi(k+2)}}
\ar@/_/ @{<-} "r-3r-3";"r-4r-4"_{\what{S}'_{1+P_k(k+3),\pi(k+3)}}
\ar@/_/ @{<-} "r-4r-4";"dot"_{\what{S}'_{1+P_k(k+4),\pi(k+4)}}
\end{xy}
\]
\[
\begin{xy}
(-5,0) *{}="00",
(0,0) *{\ }="0",
(10,0) *{\cdots}="r-1r-1",
(40,0)*{\fbox{$n-1$}^{X^L}}="r-2r-2",
(70,0)*{\fbox{$n$}^{X^L}}="r-3r-3",
(100,0)*{\fbox{$n+1$}^{X^L}}="r-4r-4",
(130,0)*{\cdots}="dot",
\ar@/^/ @{->} "r-1r-1";"r-2r-2"^{\what{S}'_{P_k(n-2),\pi(n-2)}}
\ar@/^/ @{->} "r-2r-2";"r-3r-3"^{\what{S}'_{P_k(n-1),\pi(n-1)}}
\ar@/^/ @{->} "r-3r-3";"r-4r-4"^{\what{S}'_{P_k(n),\pi(n)}}
\ar@/^/ @{->} "r-4r-4";"dot"^{\what{S}'_{P_k(n+1),\pi(n+1)}}
\ar@/_/ @{<-} "r-1r-1";"r-2r-2"_{\what{S}'_{1+P_k(n-2),\pi(n-2)}}
\ar@/_/ @{<-} "r-2r-2";"r-3r-3"_{\what{S}'_{1+P_k(n-1),\pi(n-1)}}
\ar@/_/ @{<-} "r-3r-3";"r-4r-4"_{\what{S}'_{1+P_k(n),\pi(n)}}
\ar@/_/ @{<-} "r-4r-4";"dot"_{\what{S}'_{1+P_k(n+1),\pi(n+1)}}
\end{xy}
\]
Here, we omitted index $k$.
Other actions of $\what{S}'$ are trivial. 
Note that for $r>k+1$, it holds
\begin{equation}\label{stpspr5}
1+P_k(r-1)\geq2,
\end{equation}
which is the left index of $-x_{1+P_k(r-1),\pi(r-1)}$ in (\ref{r-tilde-D0}).
Thus,
\[
\{\what{S}'_{j_{\ell}}\cdots \what{S}'_{j_1} \lambda^{(k)} | \ell\in\mathbb{Z}_{\geq0},\ j_1,\cdots,j_{\ell}\in \mathbb{Z}_{\geq1} \}
=\{0,\fbox{$r$}^{X^L}_k+\langle h_k,\lambda \rangle\ | r\in\mathbb{Z}_{\geq k+1}\}.
\]

\vspace{2mm}

\nd
\underline{Case 4 : it holds $(1,k)>(1,\pi(k+1))$, $(1,\pi(k-1))$}

\vspace{2mm}

In this case, one obtains $\lambda^{(k)}=-x_{1,k}+c_{k-1}x_{1,\pi(k-1)}+c_{k+1}x_{1,\pi(k+1)}+
\lan h_k,\lambda \ran$. 
First, let us prove
\begin{equation}\label{D-pr-1}
\{\what{S}'_{j_{\ell}}\cdots \what{S}'_{j_1} \lambda^{(k)} | \ell\in\mathbb{Z}_{\geq0},\ j_1,\cdots,j_{\ell}\in \mathbb{Z}_{\geq1} \}
\subset
\{L^{X^L}_{0,k,\iota}(T)+\langle h_k,\lambda \rangle | T\in{\rm REYD}^{X^L}_{k}\setminus\{\phi^k\}\}\cup\{0\}.
\end{equation}
We put the right hand side of (\ref{D-pr-1}) as $Z$: 
$Z:=\{L^{X^L}_{0,k,\iota}(T)+\langle h_k,\lambda \rangle | T\in{\rm REYD}^{X^L}_{k}\setminus\{\phi^k\}\}\cup\{0\}$.
Let
$T^1_k:=(y_{\ell})_{\ell\in\mathbb{Z}}\in{\rm REYD}^{X}_{k}\setminus\{\phi^k\}$ be the following:
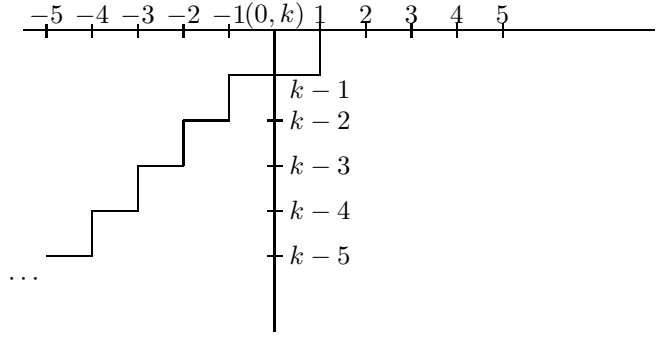
\[
\begin{xy}
(-33,0) *{}="-6",
(-6,2) *{-1}="-1",
(-12,2) *{-2}="-2",
(-18,2) *{-3}="-3",
(-24,2) *{-4}="-4",
(-30,2) *{-5}="-5",
(-6,-1) *{}="-1a",
(-12,-1) *{}="-2a",
(-18,-1) *{}="-3a",
(-24,-1) *{}="-4a",
(-30,-1) *{}="-5a",
(0,0) *{}="1",
(50,0)*{}="2",
(0,-40)*{}="3",
(0,2)*{(0,k)}="4",
(6,2) *{1}="10",
(6,-1) *{}="1010",
(12,2) *{2}="11",
(12,-1) *{}="1111",
(18,2) *{3}="12",
(18,-1) *{}="1212",
(24,2) *{4}="13",
(24,-1) *{}="1313",
(30,2) *{5}="14",
(30,-1) *{}="1414",
(6,-8)*{k-1}="k-1",
(6,-12)*{k-2}="k-2",
(6,-18)*{k-3}="k-3",
(6,-24)*{k-4}="k-4",
(6,-30)*{k-5}="k-5",
(-1,-6)*{}="5",
(1,-6)*{}="55",
(6,-6)*{}="555",
(-1,-12)*{}="6",
(-6,-12)*{}="6a",
(-12,-12)*{}="6aa",
(-12,-18)*{}="6aaa",
(-18,-18)*{}="6aaaa",
(-18,-24)*{}="st1",
(-24,-24)*{}="st2",
(-24,-30)*{}="st3",
(-30,-30)*{}="st4",
(-33,-33)*{\cdots}="stdot",
(1,-12)*{}="66",
(-1,-18)*{}="t",
(1,-18)*{}="tt",
(-6,-6)*{}="7",
(1,-18)*{}="77",
(-1,-24)*{}="8",
(1,-24)*{}="88",
(-1,-30)*{}="9",
(1,-30)*{}="99",
\ar@{-} "1";"-6"^{}
\ar@{-} "1";"2"^{}
\ar@{-} "1";"3"^{}
\ar@{-} "-1";"-1a"^{}
\ar@{-} "-2";"-2a"^{}
\ar@{-} "-3";"-3a"^{}
\ar@{-} "-4";"-4a"^{}
\ar@{-} "-5";"-5a"^{}
\ar@{-} "5";"55"^{}
\ar@{-} "55";"555"^{}
\ar@{-} "1010";"555"^{}
\ar@{-} "6";"66"^{}
\ar@{-} "t";"tt"^{}
\ar@{-} "8";"88"^{}
\ar@{-} "7";"5"^{}
\ar@{-} "7";"6a"^{}
\ar@{-} "6aa";"6a"^{}
\ar@{-} "6aaa";"6aa"^{}
\ar@{-} "6aaaa";"6aaa"^{}
\ar@{-} "st1";"6aaaa"^{}
\ar@{-} "st1";"st2"^{}
\ar@{-} "st2";"st3"^{}
\ar@{-} "st3";"st4"^{}
\ar@{-} "9";"99"^{}
\ar@{-} "10";"1010"^{}
\ar@{-} "11";"1111"^{}
\ar@{-} "12";"1212"^{}
\ar@{-} "13";"1313"^{}
\ar@{-} "14";"1414"^{}
\end{xy}
\]
The point $(1,k-1)$ is a single $k$-removable point in $T^1_k$.
We see that the point $(1,k)$ is a double (resp. single) $\pi(k+1)$-admissible point
if $c_{k+1}=2$ (resp. $c_{k+1}=1$). Similarly,
the point $(-1,k-1)$ is a double (resp. single) $\pi(k-1)$-admissible point
if $c_{k-1}=2$ (resp. $c_{k-1}=1$).
Hence, it holds $L_{0,k,\iota}^{X^L}(T^1_k)+\lan h_k,\lambda \ran = \lambda^{(k)}$ so that
$\lambda^{(k)}\in Z$.
To show (\ref{D-pr-1}), let us prove $Z$
is closed under the action of $\what{S}'_r$ for all $r\in\mathbb{Z}_{\geq1}$.
For any $T=(y_r)_{r\in\mathbb{Z}}\in {\rm REYD}^{X}_{k}\setminus\{\phi^k\}$, if the coefficient of $x_r$ is positive in 
$L^{X^L}_{0,k,\iota}(T)+\langle h_k,\lambda \rangle$ then
$T$ has an admissible point $(i,j)$ such that $j=y_i$ and
\[
r=(P_k(i+k)+[i]_-+k-j,\pi(i+k))
\]
by the definition of $L^{X^L}_{0,k,\iota}$ in (\ref{L1kdef}).
Let 
$T'=(y_t')_{t\in\mathbb{Z}}\in{\rm REYD}^{X}_{k}$ be the sequence defined as
$y_{i}'=y_i-1$ and $y_t'=y_t$ $(t\neq i)$.
It follows from Proposition \ref{A2closed}, \ref{D2closed} that
\[
\what{S}'_{r}(L^{X^L}_{0,k,\iota}(T)+\langle h_k,\lambda \rangle)
=L^{X^L}_{0,k,\iota}(T)+\langle h_k,\lambda \rangle-\beta_{r}=
L^{X^L}_{0,k,\iota}(T')+\langle h_k,\lambda \rangle\in Z.
\]
If the coefficient of $x_r$ is negative in 
$L^{X^L}_{0,k,\iota}(T)+\langle h_k,\lambda \rangle$ then
$T$ has a removable point $(i,j)$ such that $j=y_{i-1}$ and
\[
r=(P_k(i+k-1)+[i-1]_-+k-j,\pi(i+k-1)).
\]
Note that $k-j\geq1$ since $(i,j)$ is removable.
If it holds $T=T^1_k$ then
$L^{X^L}_{0,k,\iota}(T)=-x_{1,k}+c_{k-1}x_{1,\pi(k-1)}+c_{k+1}x_{1,\pi(k+1)}$ and
$r=(1,k)$. Thus,
\[
\what{S}'_r(L^{X^L}_{0,k,\iota}(T)+\langle h_k,\lambda \rangle)=0.
\]
Hence, we may assume that $T\neq T^1_k$.
In case of $i-1>0$, since $p_{\pi(k+1),\pi(k)}=1$, it holds $P_k(i+k-1)\geq1$ so that
$P_k(i+k-1)+[i-1]_-+k-j=P_k(i+k-1)+k-j\geq2$.
In case of $i-1=0$, if $k-j=1$ then $T=T^1_k$, which is absurd. If $k-j\geq2$ then
$P_k(i+k-1)+[i-1]_-+k-j=P_k(i+k-1)+k-j\geq2$.
In case of $i-1<0$, it holds $k-j\geq 2-i$
and $P_k(i+k-1)\geq p_{\pi(k-1),\pi(k)}=1$
 so that
$P_k(i+k-1)+[i-1]_-+k-j=P_k(i+k-1)+i-1+k-j\geq2$.
Hence,
\begin{equation}\label{D-pr-1a}
P_k(i+k-1)+[i-1]_-+k-j\geq2.
\end{equation}
Let
$T''=(y_t'')_{t\in\mathbb{Z}}\in{\rm REYD}^{X^L}_{k}$
be the sequence such that $y_{i-1}''=y_{i-1}+1$ and $y_t''=y_t$ $(t\neq i-1)$. 
It follows by Proposition \ref{A2closed}, \ref{D2closed} that
\begin{multline*}
\what{S}'_{r}(L^{X^L}_{0,k,\iota}(T)+\langle h_k,\lambda \rangle)
=L^{X^L}_{0,k,\iota}(T)+\langle h_k,\lambda \rangle+\beta_{r^-}\\
=
L^{X^L}_{0,k,\iota}(T)+\langle h_k,\lambda \rangle+\beta_{P_k(i+k-1)+[i-1]_-+k-j-1,\pi(i+k-1)}=
L^{X^L}_{0,k,\iota}(T'')+\langle h_k,\lambda \rangle\in Z.
\end{multline*}
Therefore, one obtains $(\ref{D-pr-1})$.

Next, we prove
\begin{equation}\label{D-pr-2}
\{\what{S}'_{j_{\ell}}\cdots \what{S}'_{j_1} \lambda^{(k)} | \ell\in\mathbb{Z}_{\geq0},\ j_1,\cdots,j_{\ell}\in \mathbb{Z}_{\geq1} \}
\supset Z.
\end{equation}
Since $0\in\{\what{S}'_{j_{\ell}}\cdots \what{S}'_{j_1} \lambda^{(k)} | \ell\in\mathbb{Z}_{\geq0},\ j_1,\cdots,j_{\ell}\in \mathbb{Z}_{\geq1} \}$
is clear, we show that for any 
$T=(y_{\ell})_{\ell\in\mathbb{Z}}\in {\rm REYD}^{X^L}_{k}\setminus\{\phi^k\}$, it holds
\begin{equation}\label{D-pr-3}
L^{X^L}_{0,k,\iota}(T)+\langle h_k,\lambda \rangle \in 
\{\what{S}'_{j_{\ell}}\cdots \what{S}'_{j_1} \lambda^{(k)} | \ell\in\mathbb{Z}_{\geq0},\ j_1,\cdots,j_{\ell}\in \mathbb{Z}_{\geq1}\}.
\end{equation}
Note that each element in ${\rm REYD}^{X^L}_{k}$ can be obtained from $\phi^k$
by putting several boxes whose length of sides are $1$.
The proof of (\ref{D-pr-3}) proceeds by induction on the number of boxes of $T$.
When the number of boxes is $1$, it holds $T=T^1_k$ and
$L^{X^L}_{0,k,\iota}(T)+\langle h_k,\lambda \rangle
=\lambda^{(k)}\in\{\what{S}'_{j_{\ell}}\cdots \what{S}'_{j_1} \lambda^{(k)} | \ell\in\mathbb{Z}_{\geq0},\ j_1,\cdots,j_{\ell}\in \mathbb{Z}_{\geq1}\}$.
Hence, we may assume
that the number of boxes of $T$ is at least $2$.
Putting $m:={\rm min}\{\ell\in\mathbb{Z}|y_{\ell}<k+\ell\}$, one obtains $m\leq 1$ by the definition of ${\rm REYD}^{X^L}_{k}$.
The equality $m=1$ means $y_{\ell}=k+\ell$ for all $\ell\in\mathbb{Z}_{\leq0}$, in particular, $y_0=k$,
which implies $y_r=k$ for any $r\in\mathbb{Z}_{\geq0}$ by $y_0\leq y_r\leq k$. Hence $T=\phi^k$, which contradicts
$T\neq \phi^k$ so that $m\leq0$ and $y_m<k$.
Putting
$m_1:={\rm min}\{y_{\ell} | m\leq \ell\}$,
we obtain $m_1<k$ and there exists $i\in\mathbb{Z}_{\geq m}$ such that $y_i=m_1$ and $y_{i+1}>y_i$. 
Note that the point $(i+1,y_i)$ is removable.
If $i=0$ and $y_0=k-1$ then $T=T^1_k$, which is absurd.
Thus, if $i=0$ then $k-1>y_0$.
In the case $i>0$, it follows from $p_{\pi(k+1),\pi(k)}=1$ that
$P_k(i+k)+[i]_-+k-y_i-1 \geq p_{\pi(k+1),\pi(k)}=1$.
In the case $i=0$, it holds $k-1>y_0$ and
$P_k(i+k)+[i]_-+k-y_i-1 =k-y_0-1\geq1$.
In the case $i<0$, using $p_{\pi(k-1),\pi(k)}=1$ and $y_i\leq k+i-1$, we obtain
$P_k(i+k)+[i]_-+k-y_i-1 \geq1$. Thus,
\begin{equation}\label{D-pr-4}
P_k(i+k)+[i]_- + k-y_i-1 \geq1.
\end{equation}
We define
$T''=(y_{\ell}'')_{\ell\in\mathbb{Z}}\in{\rm REYD}^{X^L}_{k}$
as $y_{i}''=y_{i}+1$ and $y_{\ell}''=y_{\ell}$ $(\ell\neq i)$.
Considering Proposition \ref{A2closed}, \ref{D2closed}, $L^{X^L}_{0,k,\iota}(T)$ is in the form
\[
L^{X^L}_{0,k,\iota}(T)
=L^{X^L}_{0,k,\iota}(T'')-\beta_{t_1,d_1}
\]
with $(t_1,d_1)=(P_k(i+k)+[i]_- + k-y_i-1,\pi(i+k))\in \mathbb{Z}_{\geq1}\times I$.
Repeating this argument,
\[
L^{X^L}_{0,k,\iota}(T)+\lan h_k, \lambda \ran
=L^{X^L}_{0,k,\iota}(T_k^1)-\sum_{(t,d)\in\mathbb{Z}_{\geq1}\times I} c_{t,d}\beta_{t,d}+\lan h_k, \lambda \ran
=\lambda^{(k)}-\sum_{(t,d)\in\mathbb{Z}_{\geq1}\times I} c_{t,d}\beta_{t,d}
\]
with non-negative integers $\{c_{t,d}\}$ such that $c_{t,d}=0$
except for finitely many $(t,d)\in\mathbb{Z}_{\geq1}\times I$.
Taking the assumption $T\neq T^1_k$ into account,
there exists $(t',d')\in\mathbb{Z}_{\geq1}\times I$ such that
\[
(t',d')
={\rm max} \{(t,d)\in\mathbb{Z}_{\geq1}\times I | c_{t,d}>0\}.
\]
By the definition of $\beta_{t,d}$,
the coefficient of $x_{t'+1,d'}$ is negative in 
$L^{X^L}_{0,k,\iota}(T)$. We set
\[
(t'',d''):=
{\rm min}
\{(t,d)\in\mathbb{Z}_{\geq1}\times I | \text{the coefficient of }x_{t,d}\text{ in }L^{X^L}_{0,k,\iota}(T)\text{ is negative}\}.
\]
The definition of $L^{X^L}_{0,k,\iota}$ says
there exists a removable point $(\xi,y_{\xi-1})$ in $T$
such that $L^{X^L}_{0,k,\text{re}}(\xi,y_{\xi-1})= x_{t'',d''}$
and $t''=P_k(\xi+k-1)+[\xi-1]_-+k-y_{\xi-1}$.
Just as in (\ref{D-pr-1a}), it holds
$t''\geq2$.
Let
 $T''=(y_t'')_{t\in\mathbb{Z}}\in {\rm REYD}^{X^L}_{k}$
be the sequence such that $y_{\xi-1}''=y_{\xi-1}+1$ and $y_t''=y_t$ $(t\neq \xi-1)$.
Using Proposition \ref{A2closed} and \ref{D2closed}, we have
\[
L^{X^L}_{0,k,\iota}(T)
=L^{X^L}_{0,k,\iota}(T'')-\beta_{t''-1,d''}.
\]
Considering the minimality of $(t'',d'')$ and above equality,
we see that the coefficient of $x_{t''-1,d''}$
in $L^{X^L}_{0,k,\iota}(T'')$ is positive. Thus,
\begin{equation}\label{D-pr-5}
L^{X^L}_{0,k,\iota}(T)
=\what{S}'_{t''-1,d''}L^{X^L}_{0,k,\iota}(T'').
\end{equation}
Since the number of boxes in $T''$ is less than that of $T$,
the induction assumption yields
\[
L^{X^L}_{0,k,\iota}(T'')+\langle h_k,\lambda \rangle
\in 
\{\what{S}'_{j_{\ell}}\cdots \what{S}'_{j_1} \lambda^{(k)} | \ell\in\mathbb{Z}_{\geq0},\ j_1,\cdots,j_{\ell}\in \mathbb{Z}_{\geq1}\}.
\]
In conjunction with (\ref{D-pr-5}), we get
$L^{X^L}_{0,k,\iota}(T)+\langle h_k,\lambda \rangle
\in 
\{\what{S}'_{j_{\ell}}\cdots \what{S}'_{j_1} \lambda^{(k)} | \ell\in\mathbb{Z}_{\geq0},\ j_1,\cdots,j_{\ell}\in \mathbb{Z}_{\geq1}\}$.

Finally, we take $k\in I\setminus I_{X^L}$.

\vspace{2mm}

\nd
\underline{Case 1 : it holds $(1,k)<(1,\pi'(k+1))$}

\vspace{2mm}

It holds $\lambda^{(k)}=-x_{1,k}+\langle h_k,\lambda \rangle$ and one can easily check
$\{\what{S}'_{j_{\ell}}\cdots \what{S}'_{j_1} \lambda^{(k)} | \ell\in\mathbb{Z}_{\geq0},\ j_1,\cdots,j_{\ell}\in \mathbb{Z}_{\geq1} \}
=\{0,\lambda^{(k)}\}$.

\vspace{2mm}

\nd
\underline{Case 2 : it holds $(1,k)>(1,\pi'(k+1))$}

\vspace{2mm}

We have $\lambda^{(k)}=-x_{1,k}+x_{1,\pi'(k+1)}+\langle h_k,\lambda \rangle$.
First, let us prove
\begin{equation}\label{D-pr-6}
\{\what{S}'_{j_{\ell}}\cdots \what{S}'_{j_1} \lambda^{(k)} | \ell\in\mathbb{Z}_{\geq0},\ j_1,\cdots,j_{\ell}\in \mathbb{Z}_{\geq1} \}
\subset
\{L^{X^L}_{0,k,\iota}(Y)+\langle h_k,\lambda \rangle | Y\in{\rm YW}^{X^L}_{k}\setminus\{Y_{\Lambda_k}\}\}\cup\{0\}.
\end{equation}
We put the right hand side of (\ref{D-pr-6}) as $Z$: 
$Z:=\{L^{X^L}_{0,k,\iota}(Y)+\langle h_k,\lambda \rangle | Y\in{\rm YW}^{X^L}_{k}\setminus\{Y_{\Lambda_k}\}\}\cup\{0\}$.
Let $Y^1_k\in {\rm YW}^{X^L}_{k}\setminus\{Y_{\Lambda_k}\}$ be the following:
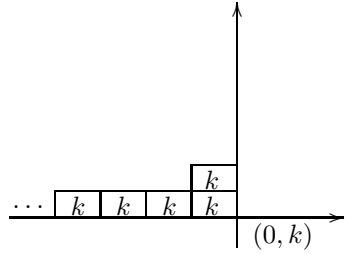
\begin{equation*}
\begin{xy}
(-15.5,-2) *{\ k}="000-1",
(-9.5,-2) *{\ k}="00-1",
(-3.5,-2) *{\ k}="0-1",
(1.5,-2) *{\ \ k}="0-1",
(-21.5,-2) *{\dots}="00000",
(1.5,1.5) *{\ \ k}="0",
(12,-6) *{(0,k)}="origin",
(0,0) *{}="1",
(0,-3.5) *{}="1-u",
(6,0)*{}="2",
(6,-3.5)*{}="2-u",
(6,-7.5)*{}="-1-u",
(6,3.5)*{}="3",
(0,3.5)*{}="4",
(-6,3.5)*{}="5",
(-6,0)*{}="6",
(-6,-3.5)*{}="6-u",
(-12,3.5)*{}="7",
(-12,0)*{}="8",
(-12,-3.6)*{}="8-u",
(-18,3.5)*{}="9",
(-18,0)*{}="10-a",
(-18,0)*{}="10",
(-18,-3.5)*{}="10-u",
(-24,-3.5)*{}="11-u",
(20,-3.5)*{}="0-u",
(6,9.5)*{}="3-1",
(6,15.5)*{}="3-2",
(6,21.5)*{}="3-3",
(6,25)*{}="3-4",
(0,9.5)*{}="4-1",
(0,15.5)*{}="4-2",
(0,21.5)*{}="4-3",
(-6,9.5)*{}="5-1",
\ar@{-} "-1-u";"2-u"^{}
\ar@{->} "2-u";"0-u"^{}
\ar@{-} "10-u";"11-u"^{}
\ar@{-} "8";"10-a"^{}
\ar@{-} "10-u";"10-a"^{}
\ar@{-} "10-u";"2-u"^{}
\ar@{-} "8";"8-u"^{}
\ar@{-} "6";"6-u"^{}
\ar@{-} "2";"2-u"^{}
\ar@{-} "1";"1-u"^{}
\ar@{-} "1";"2"^{}
\ar@{-} "1";"4"^{}
\ar@{-} "2";"3"^{}
\ar@{-} "3";"4"^{}
\ar@{-} "1";"6"^{}
\ar@{-} "6";"8"^{}
\ar@{-} "3";"3-1"^{}
\ar@{-} "3-2";"3-1"^{}
\ar@{-} "3-2";"3-3"^{}
\ar@{->} "3-3";"3-4"^{}
\end{xy}
\end{equation*}
Since $Y^1_k$ has a single removable $1$-block and 
a single $1$-admissible slot, we have
$Z\ni L^{X^L}_{0,k,\iota}(Y^1_k)+\langle h_k,\lambda \rangle=\lambda^{(k)}$.
We need to show $Z$ is closed under the action of $\what{S}'_r$ for all $r\in\mathbb{Z}_{\geq1}$.
Let $Y\in{\rm YW}^{X^L}_{k}\setminus\{Y_{\Lambda_k}\}$.
If the coefficient of $x_r$ is positive in $L^{X^L}_{0,k,\iota}(Y)$
then
the definition of $L^{X^L}_{0,k,\iota}$ in (\ref{L11-def}) says
there exists a $t$-admissible slot $S$ in $Y$ in the form
\[
\begin{xy}
(-8,15) *{(-i-1,\ell+1)}="0000",
(17,15) *{(-i,\ell+1)}="000",
(15,-3) *{(-i,\ell)}="00",
(-7,-3) *{(-i-1,\ell)}="0",
(0,0) *{}="1",
(12,0)*{}="2",
(12,12)*{}="3",
(0,12)*{}="4",
\ar@{--} "1";"2"^{}
\ar@{--} "1";"4"^{}
\ar@{--} "2";"3"^{}
\ar@{--} "3";"4"^{}
\end{xy}
\]
or
\[
\begin{xy}
(-8,9) *{(-i-1,\ell+\frac{1}{2})}="0000",
(17,9) *{(-i,\ell+\frac{1}{2})}="000",
(15,-3) *{(-i,\ell)}="00",
(-7,-3) *{(-i-1,\ell)}="0",
(0,0) *{}="1",
(12,0)*{}="2",
(12,6)*{}="3",
(0,6)*{}="4",
(38,3)*{{\rm or}}="or",
(60,9) *{(-i-1,\ell+1)}="a0000",
(85,9) *{(-i,\ell+1)}="a000",
(83,-3) *{(-i,\ell+\frac{1}{2})}="a00",
(61,-3) *{(-i-1,\ell+\frac{1}{2})}="a0",
(68,0) *{}="a1",
(80,0)*{}="a2",
(80,6)*{}="a3",
(68,6)*{}="a4",
\ar@{--} "1";"2"^{}
\ar@{--} "1";"4"^{}
\ar@{--} "2";"3"^{}
\ar@{--} "3";"4"^{}
\ar@{--} "a1";"a2"^{}
\ar@{--} "a1";"a4"^{}
\ar@{--} "a2";"a3"^{}
\ar@{--} "a3";"a4"^{}
\end{xy}
\]
such that $L^{X^L}_{0,k,{\rm ad}}(S)=x_{r}$. Hence,
$r=(P_k(\ell)+i,t)$.
Since $Y\neq Y_{\Lambda_k}$, at least one of $\ell\geq k+1$ and $i\geq1$ holds, which
implies $P_k(\ell)+i\geq1$ by $p_{\pi'(k+1),k}=1$. Let $Y'\in {\rm YW}^{X^L}_{k}\setminus\{Y_{\Lambda_k}\}$
be the proper Young wall obtained from $Y$ by adding a $t$-block to the slot $S$.
It follows from Proposition \ref{prop-closednessAw-YW} and \ref{prop-closednessCw-YW}
that
\[
\what{S}'_r(L^{X^L}_{0,k,\iota}(Y)+\langle h_k,\lambda \rangle)=
L^{X^L}_{0,k,\iota}(Y)+\langle h_k,\lambda \rangle-\beta_{r}
=
L^{X^L}_{0,k,\iota}(Y)+\langle h_k,\lambda \rangle-\beta_{P_k(\ell)+i,t}
=
L^{X^L}_{0,k,\iota}(Y')+\langle h_k,\lambda \rangle\in Z.
\]
If the coefficient of $x_r$ is negative in $L^{X^L}_{0,k,\iota}(Y)$
then
there exists a removable $t$-block $B$ in $Y$ in the form
\[
\begin{xy}
(-8,15) *{(-i-1,\ell+1)}="0000",
(17,15) *{(-i,\ell+1)}="000",
(15,-3) *{(-i,\ell)}="00",
(-7,-3) *{(-i-1,\ell)}="0",
(0,0) *{}="1",
(12,0)*{}="2",
(12,12)*{}="3",
(0,12)*{}="4",
\ar@{-} "1";"2"^{}
\ar@{-} "1";"4"^{}
\ar@{-} "2";"3"^{}
\ar@{-} "3";"4"^{}
\end{xy}
\]
or
\[
\begin{xy}
(-8,9) *{(-i-1,\ell+\frac{1}{2})}="0000",
(17,9) *{(-i,\ell+\frac{1}{2})}="000",
(15,-3) *{(-i,\ell)}="00",
(-7,-3) *{(-i-1,\ell)}="0",
(0,0) *{}="1",
(12,0)*{}="2",
(12,6)*{}="3",
(0,6)*{}="4",
(38,3)*{{\rm or}}="or",
(60,9) *{(-i-1,\ell+1)}="a0000",
(85,9) *{(-i,\ell+1)}="a000",
(83,-3) *{(-i,\ell+\frac{1}{2})}="a00",
(61,-3) *{(-i-1,\ell+\frac{1}{2})}="a0",
(68,0) *{}="a1",
(80,0)*{}="a2",
(80,6)*{}="a3",
(68,6)*{}="a4",
\ar@{-} "1";"2"^{}
\ar@{-} "1";"4"^{}
\ar@{-} "2";"3"^{}
\ar@{-} "3";"4"^{}
\ar@{-} "a1";"a2"^{}
\ar@{-} "a1";"a4"^{}
\ar@{-} "a2";"a3"^{}
\ar@{-} "a3";"a4"^{}
\end{xy}
\]
such that $L^{X^L}_{0,k,{\rm re}}(S)=x_{r}$ and
$r=(P_k(\ell)+i+1,t)$.
In the case $P_k(\ell)+i+1=1$, it holds $P_k(\ell)=i=0$ so that $\ell=k$,
which imply $t=k$ and $B$ equals
\[
\begin{xy}
(60,9) *{(-1,1)}="a0000",
(85,9) *{(0,1)}="a000",
(83,-3) *{(0,\frac{1}{2})}="a00",
(61,-3) *{(-1,\frac{1}{2})}="a0",
(68,0) *{}="a1",
(80,0)*{}="a2",
(80,6)*{}="a3",
(68,6)*{}="a4",
\ar@{-} "a1";"a2"^{}
\ar@{-} "a1";"a4"^{}
\ar@{-} "a2";"a3"^{}
\ar@{-} "a3";"a4"^{}
\end{xy}
\]
and $Y=Y^1_k$. Thus, $\what{S}'_r(L^{X^L}_{0,k,\iota}(Y)+\langle h_k,\lambda \rangle)=\what{S}'_{1,k}(\lambda^{(k)})=0$.
Hence, we may assume that
\begin{equation}\label{D-pr-a}
P_k(\ell)+i+1>1.
\end{equation}
Let $Y''\in {\rm YW}^{X^L}_{k}\setminus\{Y_{\Lambda_k}\}$
be the proper Young wall obtained from $Y$ by removing the $t$-block $B$.
It follows from Proposition \ref{prop-closednessAw-YW} and \ref{prop-closednessCw-YW}
that
\[
\what{S}'_r(L^{X^L}_{0,k,\iota}(Y)+\langle h_k,\lambda \rangle)=
L^{X^L}_{0,k,\iota}(Y)+\langle h_k,\lambda \rangle+\beta_{r^-}
=
L^{X^L}_{0,k,\iota}(Y)+\langle h_k,\lambda \rangle+\beta_{P_k(\ell)+i,t}
=
L^{X^L}_{0,k,\iota}(Y'')+\langle h_k,\lambda \rangle\in Z.
\]
If the coefficient of $x_r$ is zero then 
$\what{S}'_r(L^{X^L}_{0,k,\iota}(Y)+\langle h_k,\lambda \rangle)=
L^{X^L}_{0,k,\iota}(Y)+\langle h_k,\lambda \rangle\in Z$.
Therefore, $Z$ is closed under the action of $\what{S}'_r$ and
one gets
(\ref{D-pr-6}).

Next, let us show
\begin{equation}\label{D-pr-7}
\{\what{S}'_{j_{\ell}}\cdots \what{S}'_{j_1} \lambda^{(k)} | \ell\in\mathbb{Z}_{\geq0},\ j_1,\cdots,j_{\ell}\in \mathbb{Z}_{\geq1} \}
\supset
Z=\{L^{X^L}_{0,k,\iota}(Y)+\langle h_k,\lambda \rangle | Y\in{\rm YW}^{X^L}_{k}\setminus\{Y_{\Lambda_k}\}\}\cup\{0\}.
\end{equation}
For $Y\in{\rm YW}^{X^L}_{k}$,
we say the number of blocks in $Y$ is $m$ if $Y$ is obtained from 
$Y_{\Lambda_k}$ by adding $m$ blocks.
For $Y\in{\rm YW}^{X^L}_{k}\setminus\{Y_{\Lambda_k}\}$,
let us prove
$L^{X^L}_{0,k,\iota}(Y)+\langle h_k,\lambda \rangle\in \{\what{S}'_{j_{\ell}}\cdots \what{S}'_{j_1} \lambda^{(k)} | \ell\in\mathbb{Z}_{\geq0},\ j_1,\cdots,j_{\ell}\in \mathbb{Z}_{\geq1} \}$
by induction on the number of blocks in $Y$.
When the number is $1$, that is, $Y=Y^1_k$, it holds
$L^{X^L}_{0,k,\iota}(Y)+\langle h_k,\lambda \rangle=\lambda^{(k)}
\in \{\what{S}'_{j_{\ell}}\cdots \what{S}'_{j_1} \lambda^{(k)} | \ell\in\mathbb{Z}_{\geq0},\ j_1,\cdots,j_{\ell}\in \mathbb{Z}_{\geq1} \}$.
Hence, we may assume $Y\neq Y^1_k$.
Using Proposition \ref{prop-closednessAw-YW}, \ref{prop-closednessCw-YW},
one can describe $L^{X^L}_{0,k,\iota}(Y)$ as
\[
L^{X^L}_{0,k,\iota}(Y)+\langle h_k,\lambda \rangle=\lambda^{(k)}-\sum_{(t,d)\in\mathbb{Z}_{\geq 1}\times I} c_{t,d}\beta_{t,d}
\]
with non-negative integers $\{c_{t,d}\}$. Except for finitely many $(t,d)$, we have $c_{t,d}=0$.
As assumed  $Y\neq Y^1_k$, one can take $(t',d')\in\mathbb{Z}_{\geq 1}\times I$ as
$(t',d')
={\rm max}\{(t,d)\in\mathbb{Z}_{\geq 1}\times I |c_{t,d}>0 \}$.
Since we supposed $\iota$ is adapted, 
the coefficient of $x_{t'+1,d'}$ is negative in
$L^{X^L}_{0,k,\iota}(Y)$.
Let us take $(t'',d'')$ as
\[
(t'',d'')={\rm min}\{(t,d)\in\mathbb{Z}_{\geq 1}\times I | \text{the coefficient of }x_{t,d}\text{ in }L^{X^L}_{0,k,\iota}(Y)\text{ is negative}\}.
\]
There is a removable block $B$ in $Y$ such that
$L^{X^L}_{0,k,{\rm re}}(B)=x_{t'',d''}$.
Let $Y'\in{\rm YW}^{X^L}_{k}$ be the proper Young wall
obtained from $Y$ by removing the block $B$. From $Y\neq Y^1_k$, one can check $t''\geq2$ by a
similar argument to (\ref{D-pr-a}).
By Proposition \ref{prop-closednessAw-YW} and \ref{prop-closednessCw-YW}, we have
\[
L^{X^L}_{0,k,\iota}(Y)=
L^{X^L}_{0,k,\iota}(Y')-\beta_{t''-1.d''}.
\]
Note that
the minimality of $(t'',d'')$ means the coefficient of $x_{t''-1,d''}$ is positive in 
$L^{X^L}_{0,k,\iota}(Y')$ so that
\begin{equation}\label{D-pr-8}
\what{S}'_{t''-1,d''}(L^{X^L}_{0,k,\iota}(Y')+\lan h_k,\lambda\ran)
=L^{X^L}_{0,k,\iota}(Y')+\lan h_k,\lambda\ran-\beta_{t''-1,d''}=
L^{X^L}_{0,k,\iota}(Y)+\lan h_k,\lambda\ran.
\end{equation}
By induction assumption, one can write
$L^{X^L}_{0,k,\iota}(Y')+\lan h_k,\lambda\ran
=
\what{S}'_{j_{\xi}}\cdots \what{S}'_{j_1} \lambda^{(k)}
$
with some
$\xi\in\mathbb{Z}_{\geq0}$ and $j_1,\cdots,j_{\xi}\in \mathbb{Z}_{\geq1} $.
In conjunction with (\ref{D-pr-8}), 
one gets $L^{X^L}_{0,k,\iota}(Y)+\lan h_k,\lambda\ran \in
\{\what{S}'_{j_{\ell}}\cdots \what{S}'_{j_1} \lambda^{(k)} | \ell\in\mathbb{Z}_{\geq0},\ j_1,\cdots,j_{\ell}\in \mathbb{Z}_{\geq1} \}$
and $(\ref{D-pr-7})$. 

In this way, we get
\begin{equation}\label{stpspr6}
\{\what{S}'_{j_{\ell}}\cdots \what{S}'_{j_1} \lambda^{(k)} | \ell\in\mathbb{Z}_{\geq0},\ j_1,\cdots,j_{\ell}\in \mathbb{Z}_{\geq1} \}
={\rm Comb}^{X^L}_k[\lambda]\cup\{0\}.
\end{equation}
Combining with (\ref{half2}), we see that the $\Xi'$-ample condition $\textbf{0}=(\cdots,0,0,0)\in\Sigma'_{\iota}[\lambda]$
(Definition \ref{ample1}) holds. Using Theorem \ref{Nthm1}, we get our claim. \qed

\subsection{Proof of Theorem \ref{thm3}}

In this subsection, we prove Theorem \ref{thm3}. First, let us show $\iota$ satisfies
$\Xi'$-strict positivity condition (\ref{strict-cond}).
Since we assumed $\iota$ is adapted, it follows from Theorem 4.3, 4.5, 4.8 and 4.11 in \cite{Ka} that
$\iota$ satisfies the $\Xi'$-positivity condition.
Therefore, it holds
\[
\text{if }\ell^{(-)}=0\text{ then }\varphi_{\ell}\geq0\text{ for any }
\vp=\sum_{\ell\geq1} \varphi_{\ell}x_{\ell}\in \Xi'_{\iota}.
\]
To show (\ref{strict-cond}), we need to prove
for any $k\in I$ and $\vp=\sum_{\ell\geq1} \varphi_{\ell}x_{\ell}\in \Xi'^{(k)}_{\iota}\setminus\{\xi^{(k)}\}$,
if $\ell^{(-)}=0$ then $\varphi_{\ell}\geq0$.
One can write $\vp=S'_{j_{r}}\cdots S'_{j_1}\xi^{(k)}$
with some
$r\in\mathbb{Z}_{\geq1}$ and $j_1,\cdots,j_{r}\in\mathbb{Z}_{\geq1}$.
We may assume that for any $t=1,2,\cdots,r$, it holds 
$S'_{j_t}(S'_{j_{t-1}}\cdots S'_{j_1}\xi^{(k)})\neq S'_{j_{t-1}}\cdots S'_{j_1}\xi^{(k)}$.
By (\ref{twoS}), this assumption means
$\vp=S'_{j_{r}}\cdots S'_{j_1}\xi^{(k)}=\what{S}'_{j_r}\cdots \what{S}'_{j_1}\xi^{(k)}$.
Note that when $\lambda=0$, it holds $\lambda^{(k)}=\xi^{(k)}$.
Thus, the relation (\ref{stpspr3}) and (\ref{stpspr6}) yields $\vp\in {\rm Comb}_{k,\iota}^{X^L}[0]$.
Therefore, by the same argument as in (\ref{stpspr1}), (\ref{stpspr2}), (\ref{A-pr-1a}), (\ref{stpspr4}), (\ref{stpspr5}), (\ref{D-pr-1a}) and (\ref{D-pr-a}), one obtains if $\ell^{(-)}=0$ then $\varphi_{\ell}\geq0$. Hence,
$\Xi'$-strict positivity condition holds.

Therefore, as a consequence of (\ref{thm2-pr2}), Theorem \ref{thm2}, (\ref{stpspr3}) and (\ref{stpspr6}),
we see that
\begin{eqnarray*}
\varepsilon_k^*(x)
&=&{\rm max}\{-\varphi(x)|\varphi\in \Xi'^{(k)}_{\iota}\}\\
&=&
{\rm max}\{-\varphi(x)|\varphi\in 
\{\what{S}'_{j_{\ell}}\cdots \what{S}'_{j_1} 0^{(k)} | \ell\in\mathbb{Z}_{\geq0},\ j_1,\cdots,j_{\ell}\in \mathbb{Z}_{\geq1} \}\setminus\{0\}\}\\
&=&
{\rm max}\{-\varphi(x)|\varphi\in {\rm Comb}_{k,\iota}^{X^L}[0]\}.
\end{eqnarray*}
\qed

\vspace{2mm}

\nd
{\bf Conflict of Interest}
The corresponding author states that there is no conflict of interest.

\vspace{2mm}

\nd
{\bf Data availability statement}
All data generated or analysed during this study are included in this published article.


\begin{thebibliography}{9}







\bibitem{BZ}
A.Berenstein, A.Zelevinsky, 
Tensor product multiplicities, canonical bases and totally positive varieties,
Invent. Math. 143, no. 1, 77--128 (2001).


\bibitem{GKS16} V.Genz, G.Koshevoy and B.Schumann, Combinatorics of canonical bases revisited: type  A, 
Selecta Math. (N.S.) 27, no. 4, Paper No. 67, 45 pp, (2021).


\bibitem{GP}
O.Gleizer, A.Postnikov,
Littlewood-{R}ichardson coefficients via {Y}ang-{B}axter equation.
Internat. Math. Res. Notices 14, 741--774 (2000).

\bibitem{Ha}
T.Hayashi,
$Q$-analogues of Clifford and Weyl algebras—spinor and oscillator representations of quantum enveloping algebras,
Comm. Math. Phys. 127, no. 1, 129--144 (1990).




\bibitem{HN2}D.Hernandez, H.Nakajima,
Level 0 monomial crystals,
Nagoya Math. J. 184, 85--153 (2006).

\bibitem{H1}
A.Hoshino,
Polyhedral realizations of crystal bases for quantum algebras of finite types,
J. Math. Phys. 46, no. 11, 113514,  31 pp, (2005). 

\bibitem{H2}
A.Hoshino,
Polyhedral realizations of crystal bases for quantum algebras of classical affine types,
J. Math. Phys. 54, no. 5, 053511, 28 pp, (2013).


\bibitem{HN}
A.Hoshino, K.Nakada,
Polyhedral realizations of crystal bases $B(\lambda)$ for quantum algebras of nonexceptional affine types,
J. Math. Phys. 60, no. 9, 56 pp (2019). 

\bibitem{JMMO}
M.Jimbo, K. C. Misra, T.Miwa, M.Okado,
Combinatorics of representations of  
$U_q(\widehat{\mathfrak{s}\mathfrak{l}}(n))$ at $q=0$,
Comm. Math. Phys. 136, no. 3, 543-–566 (1991).



\bibitem{Ka} Y.Kanakubo, Polyhedral realizations for $B(\infty)$
and extended Young diagrams, Young walls of type
${\rm A}^{(1)}_{n-1}$, ${\rm C}^{(1)}_{n-1}$, ${\rm A}^{(2)}_{2n-2}$, ${\rm D}^{(2)}_{n}$,
Algebr. Represent. Theor., https://doi.org/10.1007/s10468-022-10172-z (2022). 


\bibitem{KaN} Y.Kanakubo, T.Nakashima, Adapted sequence for polyhedral realization of crystal bases,
Comm. Algebra 48, no. 11, 4732--4766 (2020).


\bibitem{KaN2} Y.Kanakubo, T.Nakashima, Adapted sequences and polyhedral realizations of crystal bases for highest weight modules,
J. Algebra 574, 327--374 (2021).



\bibitem{Kang}
S.-J. Kang, Crystal bases for quantum affine algebras and combinatorics of Young walls,
Proc. London Math. Soc. (3) 86, no. 1, 29-–69 (2003).


\bibitem{KK}
S.-J. Kang, J.-H. Kwon, 
Crystal bases of the Fock space representations and string functions,
J. Algebra 280, no. 1, 313--349 (2004).



\bibitem{KMM}
S.-J. Kang, K. C. Misra, T.Miwa,
Fock space representations of the quantized universal enveloping algebras
$U_q(C_{\ell}^{(1)})$, $U_q(A_{2l}^{(2)})$, and $U_q(D_{l+1}^{(2)})$,
J. Algebra 155, no. 1, 238-–251 (1993).







\bibitem{K0} M.Kashiwara, Crystalling the $q$-analogue of universal 
              enveloping algebras, Comm. Math. Phys.,
            {\it 133}, 249--260 (1990).


\bibitem{K1} M.Kashiwara,
 On crystal bases of the $q$-analogue of universal enveloping algebras,
	Duke Math. J., {\it 63} (2), 465--516 (1991).

\bibitem{K3}M.Kashiwara, 
The crystal base and Littelmann's refined Demazure character formula,
Duke Math. J., 71, no 3,  839--858 (1993).

\bibitem{K}M.Kashiwara, Realizations of crystals, 
Combinatorial and geometric representation theory (Seoul, 2001), 133-–139,
Contemp. Math., 325, Amer. Math. Soc., Providence, RI, (2003).



\bibitem{KN}M.Kashiwara, T.Nakashima,
Crystal graphs for representations of the $q$-analogue of classical Lie algebras,
J. Algebra 165, no. 2, 295--345 (1994).


\bibitem{KS}\label{KS}
Kim, J.-A., Shin, D.-U., 
Monomial realization of crystal bases $B(\infty)$ for the quantum finite algebras,
Algebr. Represent. Theory 11, no. 1, 93--105 (2008).



\bibitem{Lit2}
P. Littelmann, Paths and root operators in representation theory, Annals
of Math. 142, 499--525 (1995).


\bibitem{Lit}
P.Littelmann, 
Cones, crystals, and patterns,
Transform. Groups 3, no. 2, 145--179 (1998).




\bibitem{L}
G.Lusztig, Canonical bases arising from quantized enveloping algebras, 
J. Amer. Math. Soc. 3, no. 2, 447--498 (1990). 



\bibitem{Nj}H.Nakajima, $t$-analogs of $q$-characters of quantum affine algebras of type $A_n$, $D_n$,
Combinatorial and geometric representation theory (Seoul, 2001), 141--160,
Contemp. Math., 325, Amer. Math. Soc., Providence, RI (2003).

\bibitem{N99}T.Nakashima, Polyhedral realizations of crystal bases for integrable highest weight modules,
J. Algebra, vol.219, no. 2, 571–-597, (1999).

\bibitem{NZ}
T.Nakashima,  A.Zelevinsky, Polyhedral realizations of 
crystal bases for quantized Kac-Moody algebras,
Adv. Math. {\bf 131}, no. 1, 253--278, (1997). 



\end{thebibliography}
\end{document}